\documentclass[11pt,oneside]{amsart}

\usepackage{tikz-cd}

\usepackage{fix-cm}

\DeclareMathAlphabet{\mathcal}{OMS}{cmsy}{m}{n}

\usepackage[arrow,curve,matrix,tips,2cell]{xy}

\usepackage[hmargin=2cm,vmargin=1.5cm]{geometry}

\usepackage{anyfontsize}
\usepackage{mathrsfs}
\usepackage{bm}

\usepackage{bbm}

\usepackage{fancyhdr}
\usepackage[obeyspaces]{url}

\usepackage{tensor}

\usepackage{mathtools}

\usepackage{accents}

\usepackage{amsmath}
\usepackage{amscd}
\usepackage{amssymb}
\usepackage{latexsym}
\usepackage{url}

\usepackage{graphicx}

\newtheorem{theorem}{Theorem}[section]
\newtheorem*{theorem*}{Theorem}
\newtheorem{lemma}[theorem]{Lemma}
\newtheorem*{lemma*}{Lemma}
\newtheorem{corollary}[theorem]{Corollary}
\newtheorem{proposition}[theorem]{Proposition}

\newtheorem{remark}[theorem]{Remark}
\newtheorem{definition}[theorem]{Definition}
\newtheorem*{definition*}{Definition}
\newtheorem{question}[theorem]{Question}
\newtheorem*{question*}{Question}

\newtheorem{example}[theorem]{Example}
\newtheorem{examples}[theorem]{Examples}

\makeatletter
\def\revddots{\mathinner{\mkern1mu\raise\p@
\vbox{\kern7\p@\hbox{.}}\mkern2mu
\raise4\p@\hbox{.}\mkern2mu\raise7\p@\hbox{.}\mkern1mu}}
\makeatother % ddots gespiegelt
\newcommand{\bgl}{\begin{equation}} %eine Gleichung mit Ziffer
\newcommand{\egl}{\end{equation}}
\newcommand{\bgloz}{\begin{equation*}} %eine Gleichung ohne Ziffer
\newcommand{\egloz}{\end{equation*}}
\newcommand{\bgln}{\begin{eqnarray}} %mehrere Gleichungen mit Ziffer
\newcommand{\egln}{\end{eqnarray}}
\newcommand{\bglnoz}{\begin{eqnarray*}} %mehrere Gleichungen ohne Ziffer
\newcommand{\eglnoz}{\end{eqnarray*}}
\newcommand{\btheo}{\begin{theorem}}
\newcommand{\etheo}{\end{theorem}}
\newcommand{\btheooz}{\begin{theorem*}}
\newcommand{\etheooz}{\end{theorem*}}

\newcommand{\blemma}{\begin{lemma}}
\newcommand{\elemma}{\end{lemma}}
\newcommand{\blemmaoz}{\begin{lemma*}}
\newcommand{\elemmaoz}{\end{lemma*}}
\newcommand{\bproof}{\begin{proof}}
\newcommand{\eproof}{\end{proof}}
\newcommand{\bbew}{\begin{beweis}}
\newcommand{\ebew}{\end{beweis}}
\newcommand{\bremark}{\begin{remark}\em}
\newcommand{\eremark}{\end{remark}}
\newcommand{\bdefin}{\begin{definition}}
\newcommand{\edefin}{\end{definition}}
\newcommand{\bdefinoz}{\begin{definition*}}
\newcommand{\edefinoz}{\end{definition*}}
\newcommand{\bex}{\begin{example}\em}
\newcommand{\eex}{\end{example}}
\newcommand{\bexs}{\begin{examples}}
\newcommand{\eexs}{\end{examples}}

\newcommand{\bprop}{\begin{proposition}}
\newcommand{\eprop}{\end{proposition}}
\newcommand{\bcor}{\begin{corollary}}
\newcommand{\ecor}{\end{corollary}}
\newcommand{\bfa}{\begin{cases}} %Fallunterscheidung
\newcommand{\efa}{\end{cases}}

\newcommand{\bquestion}{\begin{question}}
\newcommand{\equestion}{\end{question}}
\newcommand{\bquestionoz}{\begin{question*}}
\newcommand{\equestionoz}{\end{question*}}
\newtheorem{introtheorem}{Theorem}

\newtheorem{introcor}[introtheorem]{Corollary}

\newcommand{\cC}{\mathcal C}
\newcommand{\cD}{\mathcal D}
\newcommand{\cE}{\mathcal E}
\newcommand{\cF}{\mathcal F}
\newcommand{\cG}{\mathcal G}

\newcommand{\cJ}{\mathcal J}

\newcommand{\cL}{\mathcal L}
\newcommand{\cM}{\mathcal M}
\newcommand{\cN}{\mathcal N}
\newcommand{\cO}{\mathcal O}

\newcommand{\cS}{\mathcal S}

\newcommand{\cU}{\mathcal U}
\newcommand{\cV}{\mathcal V}

\def\Nz{\mathbb{N}}

\def\Rz{\mathbb{R}}

\def\Zz{\mathbb{Z}}

\newcommand{\fA}{\mathfrak A}

\newcommand{\fC}{\mathfrak C}
\newcommand{\fD}{\mathfrak D}

\newcommand{\fG}{\mathfrak G}

\newcommand{\fS}{\mathfrak S}

\newcommand{\fV}{\mathfrak V}
\newcommand{\fX}{\mathfrak X}

\newcommand{\mfd}{\mathfrak d}
\newcommand{\mfe}{\mathfrak e}
\newcommand{\mff}{\mathfrak f}
\newcommand{\mfg}{\mathfrak g}
\newcommand{\mfh}{\mathfrak h}

\newcommand{\mfs}{\mathfrak s}
\newcommand{\mft}{\mathfrak t}

\newcommand{\mfv}{\mathfrak v}
\newcommand{\mfw}{\mathfrak w}
\newcommand{\mfx}{\mathfrak x}
\newcommand{\an}[1]{``#1''} % Anfuehrungsstriche
\newcommand{\ti}{\tilde}

\newcommand{\ma}{\mapsto} % wird abgebildet auf
\newcommand{\into}{\hookrightarrow} % injektiv
\newcommand{\isom}{\xrightarrow{\raisebox{-1ex}[0ex][0ex]{$\sim$}}} % isom
\def\SEMI{\mbox{$\times\kern-2pt\vrule height5pt width.6pt \kern3pt $}}

\newcommand{\id}{{\rm id}}

\newcommand{\img}{{\rm im\,}}

\newcommand{\lcm}{{\rm lcm}} % kleinstes gemeinsames Vielfaches
\newcommand{\reg}{^\times} % regulaer
\newcommand{\defeq}{\mathrel{:=}} % per Definition
\newcommand{\dop}{\text{: }} % in Mengen
\newcommand{\lge}{\left\{} % links geschweift
\newcommand{\rge}{\right\}} % rechts geschweift
\newcommand{\lsp}{\left\langle} % links spitz
\newcommand{\rsp}{\right\rangle} % links spitz
\newcommand{\gekl}[1]{\lge #1 \rge} % geschweifte Klammer
\newcommand{\spkl}[1]{\lsp #1 \rsp} % spitze Klammer
\newcommand{\menge}[2]{\gekl{ #1 \dop #2 }} % Menge
\def\bf1{\mathbf{1}}
\newcommand{\dom}{{\rm dom\,}}
\newcommand{\ran}{{\rm ran\,}}
\newcommand{\bfB}{\bm{\mathfrak B}}
\newcommand{\bfC}{\bm{\mathfrak C}}
\newcommand{\bfD}{\bm{\mathfrak D}}

\newcommand{\bfS}{\bm{\mathfrak S}}
\newcommand{\bfX}{\bm{\mathfrak X}}
\newcommand{\mcm}{{\rm mcm}}

\newcommand{\bGamma}{\bm{\Gamma}}
\newcommand{\bma}{\bm{a}}
\newcommand{\bmb}{\bm{b}}
\newcommand{\bmc}{\bm{c}}
\newcommand{\bmd}{\bm{d}}
\newcommand{\bmg}{\bm{g}}
\newcommand{\bmm}{\bm{m}}
\newcommand{\bmt}{\bm{t}}
\newcommand{\bmu}{\bm{u}}

\newcommand{\bmE}{\bm{E}}
\newcommand{\bmF}{\bm{F}}
\newcommand{\bmU}{\bm{U}}

\newcommand{\bbd}{\mathbbm{d}}

\newcommand{\Div}{{\rm Div}}

\newcommand{\rmr}{\ensuremath{\mathrm{r}}}
\newcommand{\rms}{\ensuremath{\mathrm{s}}}

\newcommand{\oset}[2]{%
  \mathop{#2}\limits^{\vbox to -1.66ex{%
  \kern -1.4ex\hbox{$#1$}\vss}}}

\newcommand{\tipreceq}{\oset{\sim}{\preceq}}
\newcommand{\tiprec}{\oset{\sim}{\prec}}

\newcommand{\C}{\mathscr{C}}
\newcommand{\D}{\mathscr{D}}
\newcommand{\Q}{\mathscr{Q}}
\newcommand{\crS}{\mathscr{S}}

\newcommand{\pars}{\setlength{\parindent}{0cm} \setlength{\parskip}{0.5cm}}
\newcommand{\pari}{\setlength{\parindent}{0.5cm} \setlength{\parskip}{0cm}}
\newcommand{\nopar}{\setlength{\parindent}{0cm} \setlength{\parskip}{0cm}}

\usepackage{newtxtext}
\usepackage{newtxmath}

\begin{document}

\title[Left regular representations of Garside categories II]{Left regular representations of Garside categories II. \\ Finiteness properties of topological full groups}

\thispagestyle{fancy}

\author{Xin Li}

\address{Xin Li, School of Mathematics and Statistics, University of Glasgow, University Place, Glasgow G12 8QQ, United Kingdom}
\email{Xin.Li@glasgow.ac.uk}


\thanks{This project has received funding from the European Research Council (ERC) under the European Union's Horizon 2020 research and innovation programme (grant agreement No. 817597).}

\begin{abstract}
We study topological full groups attached to groupoid models for left regular representations of Garside categories. Groups arising in this way include Thompson's group $V$ and many groups of dynamical origin such as R{\"o}ver-Nekrashevych groups. Our key observation is that a Garside structure for the underlying small category induces a new Garside structure for a new small category of bisections, and that our topological full group coincides with the fundamental group of the enveloping groupoid of the new category. As a consequence, we solve the word problem and identify general criteria for establishing finiteness properties of our topological full groups. In particular, we show that topological full groups arising from products of shifts of finite type are of type ${\rm F}_{\infty}$, answering a natural question left open by Matui.
\end{abstract}

\maketitle

%\tableofcontents

\setlength{\parindent}{0cm} \setlength{\parskip}{0.5cm}

%$\cG$ {\'e}tale groupoid, ample; $X = \cG^{(0)}$ unit space; $\bm{F}(\cG)$ topological full group.
%
%$U, V, W \in \cC \cO$
%
%$\bm{\fB}$, $Q(\bm{\fB})$: SC, groupoid of general compact open bisections
%
%
%
%$\bm{\sigma} \in \bm{\fB}(\bm{U},\bm{V})$: morphism, $\bm{U} = (U_i)_i$, $\bm{V} = (V_j)_j$, $\sigma_i: \: U_i \to V_j$.
%
%
%
%$\fC$: LCSC, $a, b, c, d$: elements in $\fC$; $\fC^0$: objects; $\mfu, \mfv, \mfw \in \fC^0$; $\mfd$, $\mft$: domain, target; $\fA$: atoms
%
%$\fC^*$: invertible elements of $\fC$; $a \sim b$ $\Leftrightarrow$ $a \in b \fC^*$; $a \preceq b$ $\Leftrightarrow$ $b \in a \fC$, $a \prec b$ if $a \preceq b$ and $a \not\sim b$
%
%$\lcm$: right-lcm, $\mcm$: right-mcm, well-defined up to $\sim$
%
%$I_l$: left inverse hull (inverse semigroup); $E$: semilattice; $\Omega$: $\widehat{E}$ (spectrum)
%
%$\bfC$: SC of $\fC$-bisections; $\bmd$, $\bmt$: domain and target in $\bfC$; $\bma, \bmb, \bmc$: elements of $\bfC$
%
%$\widehat{E}(-;-)$, $\Omega(-;-)$, $X(-;-)$: basic compact open sets
%
%
%
%$\fA$: set of atoms in $\fC$
%
%
%
%$\fS$: Garside family, $\Vert \cdot \Vert$: $\fS$-length
%
%$\bfX_{\fS}$, $\bfX$: collection of domains and targets in $\bfC$.
%
%$\Gamma$: generator (subset of $\bm{\fC}$); $\alpha, \beta, \gamma$: elements in $\Gamma$
%
%$\cC$: sub-category of $\bm{\fC}$ generated by $\Gamma$ and units

\section{Introduction}

The story of this paper starts with Thompson's group $V$, which was introduced in unpublished notes by Thompson in 1965. It can be described as a group of certain right-continuous bijections of the unit interval $[0,1]$ and was one of the first examples of a finitely presented infinite simple group (see for instance \cite{CFP} for more details). Thompson's group $V$ has been studied from many different perspectives, and its construction has been generalized in several different ways (see for instance \cite{Hig,Ste,Bri04}). One perspective is based on the observation, going back to Nekrashevych \cite{Nek04}, that $V$ can be described as the topological full group of a groupoid naturally attached to the one-sided full shift on two symbols. This led to the construction of interesting classes of groups of dynamical origin, for instance R{\"o}ver-Nekrashevych groups, which are topological full groups of groupoids attached to expanding maps \cite{Rov,Nek04,Nek18a}. Moreover, replacing the full shift by shifts of finite type, Matui constructed and studied topological full groups arising from shift of finite type groupoids and their products \cite{Mat15,Mat16}. Roughly speaking, the general idea is that {\'e}tale groupoids encode generalized topological dynamical systems because they generate partial homeomorphisms on their underlying unit spaces, and the associated topological full groups consist of global symmetries (i.e., globally defined homeomorphisms) of the unit space which are pieced together from these partial homeomorphisms. These topological full groups turn out to have very interesting properties, leading to solutions of outstanding open problems in group theory. For instance, they provided the first examples of finitely generated, infinite, simple, amenable groups \cite{JM} and of finitely generated, infinite, simple groups of intermediate growth \cite{Nek18b}.

There is also an interesting connection to C*-algebra theory because {\'e}tale groupoids serve as models for C*-algebras \cite{Li20} and introduce ideas and techniques from topological dynamics to the study of C*-algebras and their structural properties. For instance, the groupoid attached to the one-sided full shift on two symbols mentioned above is nothing else but the canonical groupoid model of the Cuntz algebra $\cO_2$ \cite{Cun}, a classical example of a C*-algebra which plays a distinguished role in the Elliott classification programme. 

A unifying framework is given by left regular representations of left cancellative small categories and their groupoid models (see \cite{Spi20,Li21a}). The corresponding topological full groups generalize all of the above-mentioned examples. Indeed, the underlying groupoids of R{\"o}ver-Nekrashevych groups can be described as boundary groupoids attached to left regular representations of Zappa-Sz{\'e}p product monoids arising from self-similar groups, while the groups studied in \cite{Mat15,Mat16} are topological full groups of boundary groupoids associated with small categories obtained from shifts of finite type (and products of them). Generally speaking, topological full groups of boundary groupoids obtained from these left regular representations build a bridge between group theory, dynamical systems (in the form of topological groupoids) and C*-algebras and are the main objects of study of this paper. We will mainly focus on finiteness properties of this rich class of topological full groups. A group is of type ${\rm F}_n$ if it admits a classifying space with a compact $n$-skeleton. These finiteness properties play an important role in the study of group homology and reduce to familiar notions in low dimensions (a group is of type ${\rm F}_1$ if and only if it is finitely generated and of type ${\rm F}_2$ if and only if it is finitely presented). The study of finiteness properties of $V$ and its generalizations goes back to \cite{BG,Bro,Ste} (see also for instance \cite{FMWZ,SWZ} and the references therein). In \cite{Nek18a}, it is shown that topological full groups arising from classes of self-similar groups are finitely presented, i.e., of type ${\rm F}_2$. \cite{SWZ} established that topological full groups arising from classes of self-similar groups give rise to first examples of infinite simple groups which are of type ${\rm F}_{n-1}$ but not of type ${\rm F}_n$, for each $n$. In this context, a general framework for establishing finiteness properties has been developed in \cite{Wit19}. The key ingredient is the notion of a Garside category (see \cite{Deh15}), which originated from the study of Braid groups and allow us to carry over classical results and methods from Braid groups to more general groups, monoids or small categories. Applications of Garside structures have been developed in the context of the $K(\pi,1)$-conjecture \cite{Bes,Par,PS} as well as isomorphism conjectures such as the Farrell-Jones conjecture or the coarse Baum-Connes conjecture because of connections with Helly graphs \cite{HO}. Roughly speaking, \cite{Wit19} provides criteria when isotropy groups of Garside categories have prescribed finiteness properties. Using similar techniques as in \cite{BG,Bro,Ste}, it was shown in \cite{Mat15} that topological full groups arising from shifts of finite type are of type ${\rm F}_{\infty}$ (i.e., type ${\rm F}_n$ for all $n$). The natural question whether the same is true for topological full groups arising from products of shifts of finite type was left open in \cite[\S~5.3]{Mat16}. 

Our goal is to show that for topological full groups of groupoids arising from left regular representations of Garside categories, the Garside structure of the underlying small category induces a new Garside structure on the category of bisections of the groupoids, whose isotropy groups (or rather those of the enveloping groupoid) can be identified with the topological full groups we are interested in. The construction of this new Garside structure is a key novelty of this paper. Consequently, we succeed in solving the word problem and establishing finiteness properties for topological full groups of groupoids attached to left regular representations of several classes of Garside categories. For the second result, the technical difficulty we overcome is the identification of sufficient conditions for connectivity of certain simplicial complexes arising in our setting (see \S~\ref{s:Fn} for details). In particular, we answer the natural question left open in \cite[\S~5.3]{Mat16} and show that topological full groups arising from products of shifts of finite type are of type ${\rm F}_{\infty}$. This class of groups generalize higher dimensional analogues of Thompson's group $V$ such as the groups $nV$ from \cite{Bri04} and is an instructive example class illustrating the key ideas of this paper (see \S~\ref{ss:k-graphs} for details).

Let us now formulate our main results. Given a left cancellative small category $\fC$, we denote its set of objects by $\fC^0$ and its set of invertible elements by $\fC^*$. The groupoid model $I_l \ltimes \Omega$ attached to $\fC$ is discussed in \cite{Li21a} and is recalled in \S~\ref{ss:GPD_LCSC}. In the following, we will present our main theorems only in the special case of the boundary groupoid $I_l \ltimes \partial \Omega$ for simplicity, even though our results are much more general (they apply to groupoids of the form $(I_l \ltimes X)_Y^Y$ where $X$ is a closed invariant subspace of $\Omega_{\infty}$ and $Y$ is a closed subset of $X$ of a particular form). Starting with a Garside family for $\fC$ (the notion of Garside family is explained in \S~\ref{s:CatBisGar}), we first construct another small category $\C$ which again admits a Garside family such that the topological full group $\bmF(I_l \ltimes \partial \Omega)$ can be identified with an isotropy group of the enveloping groupoid of $\C$.
\begin{introtheorem}
\label{intro:Gars}
Suppose that $\fC$ is a left cancellative small category with finite $\fC^0$. Further assume that $\fC$ is right cancellative up to $=^*$, finitely aligned, right Noetherian and admits disjoint mcms, and that (F) holds. Let $\fS$ be a Garside family in $\fC$ which is locally finite, $=^*$-transverse with $\fS \cap \fC^* = \emptyset$. Moreover assume that for all $L \geq 1$, $(\fS^{\leq L})^{\sharp}$ is closed under left divisors.
\pari

Let $\C$ be the small category constructed in \S~\ref{ss:GPD_LCSC} and \S~\ref{s:CatBisGar}, with base object $*$ and enveloping groupoid $\Q$. Then there exists a right Garside map $\Delta$ for $\C$ such that $\Div_{\C}(\Delta)$ is a Garside family for $\C$, and we have $\Q(*,*) \cong \bmF(I_l \ltimes \partial \Omega)$. 
\end{introtheorem}
\nopar

$\C$ is constructed explicitly as a category of bisections of our groupoid (see \S~\ref{ss:GPD_LCSC} and \S~\ref{s:CatBisGar}). The construction of $\C$ is interesting on its own right as similar constructions led to a better understanding of group homology for topological full groups \cite{Li22}. Theorem~\ref{intro:Gars} is proved in \S~\ref{s:GarsCat-TFG} (see Corollary~\ref{cor:thmA}), where the reader will find more explanations and details. We restrict to the case of small categories $\fC$ with finite set of objects $\fC^0$ for convenience and because this is required for the applications we have in mind. The point about Theorem~\ref{intro:Gars} is that it allows us to apply tools for Garside categories as in \cite{Deh15} in the study of topological full groups. 
\pars

For example, Theorem~\ref{intro:Gars} allows us to solve the word problem for topological full groups of the form $\bmF(I_l \ltimes \partial \Omega)$. 
\begin{introcor}
\label{intro:WordProblem}
Assume that we are in the same situation as in Theorem~\ref{intro:Gars}. Suppose that there exists a computable $=^*$-map for $\fS^{\sharp}$. Then $\bmF(I_l \ltimes \partial \Omega)$ has decidable word problem.
\end{introcor}
\nopar

Here a $=^*$-map for $\fS^{\sharp} = \fS \fC^* \cup \fC^*$ is a partial map $E$ from $\fS^{\sharp} \times \fS^{\sharp}$ to $\fC^*$ with the property that $E(s,t)$ is defined if and only if $s =^* t$, and in that case $E(s,t) = u \in \fC^*$ with $su = t$. Such a $=^*$-map is called computable if it can be implemented on a Turing machine. Corollary~\ref{intro:WordProblem} is proved in \S~\ref{s:GarsCat-TFG} (see Corollary~\ref{cor:WordProblem}), where we construct a concrete algorithm to solve the word problem, i.e., to decide whether a given word in the generators represents the trivial element of our topological full groups.
\pars

As another application of Garside structures, we establish general criteria for finiteness properties of topological full groups.
\begin{introtheorem}
\label{intro:Fn}
Let $\fC$ and $\fS$ be as in Theorem~\ref{intro:Gars}. Assume that conditions (F), (St), (LCM) and ($\bmt < \bmd$) are satisfied. Then for all natural numbers $n$, $\bmF(I_l \ltimes \partial \Omega)$ is of type ${\rm F}_n$ if $\fC^*(\mfv,\mfv)$ is of type ${\rm F}_n$ for all $\mfv \in \fC^0$.
\end{introtheorem}
\nopar

The conditions (F), (St), (LCM) and ($\bmt < \bmd$) are introduced in \S~\ref{s:Fn}, which also contains the proof of Theorem~\ref{intro:Fn} (see Theorem~\ref{thm:Fn}). On our way of proving Theorem~\ref{intro:Fn}, we construct a concrete simplicial complex on which our topological full group acts and for which our conditions allow us to establish the desired connectivity properties. The construction of the complex is interesting on its own right as we expect a further analysis of it to reveal valuable information about topological full groups, for instance regarding concrete presentations (see for instance \cite{Mat15}).
\pars

Having established a general criterion for finiteness properties of topological full groups, we now turn to specific contexts where the conditions in Theorem~\ref{intro:Fn} are satisfied. In the presence of degree maps, we identify criteria for finiteness properties which are easy to check.

\begin{introcor}
\label{intro:deg}
Let $P$ be a left cancellative monoid with $P^* = \gekl{1}$. Assume that $P$ is right Noetherian, left reversible and admits conditional lcms. Suppose that $S_P \subseteq P$ is a finite Garside family in $P$ with $1 \notin S_P$ and assume that $(S_P^{\leq L})^{\sharp}$ is closed under left divisors for all $L \geq 1$. Let $\fC$ be a left cancellative small category with finite $\fC^0$ equipped with a $P$-valued degree map $\bbd$ such that $\mfv \dot{\bbd}^{-1}(p) < \infty$ for all $\mfv \in \fC^0$ and $p \in P$. Suppose that condition (F) holds.
\pari
 
If for all $\mfv \in \fC^0$ and $s \in S_P$, we have $\# \mfv \bbd^{-1}(s) \mfv \geq 2$, then for all natural numbers $n$, $\bmF(I_l \ltimes \partial \Omega)$ is of type ${\rm F}_n$ if $\fC^*(\mfv,\mfv)$ is of type ${\rm F}_n$ for all $\mfv \in \fC^0$.
\end{introcor}
\nopar

Degree maps are introduced in \S~\ref{s:Gars-deg}, where Corollary~\ref{intro:deg} is proved, too (see Theorem~\ref{thm:deg}). Condition (F) is for example satisfied if $\fC$ is cancellative. Corollary~\ref{intro:deg} applies for instance to small categories arising from higher rank graphs (where $P = \Zz_{\geq 0}^k$). In particular, we obtain the following special case.
\pars

\begin{introcor}
\label{intro:productgraphs}
Suppose that $\fC_j$, $1 \leq j \leq k$, are path categories of finite graphs. Assume that for all $1 \leq j \leq k$ and $\mfv \in \fC_j^0$, we have $\# \mfv \fC_j \mfv \geq 2$. Let $\cG_j \defeq I_l(\fC_j) \ltimes \partial \Omega_{\fC_j}$. Then $\bmF(\cG_1 \times \dotso \times \cG_k)$ is of type ${\rm F}_{\infty}$. 
\pari

In particular, if $\cG_j$, $1 \leq j \leq k$, are groupoids arising from irreducible one-sided shifts of finite type as in \cite{Mat15}, then $\bmF(\cG_1 \times \dotso \times \cG_k)$ is of type ${\rm F}_{\infty}$.
\end{introcor}
\nopar

This answers a natural question left open in \cite[\S~5.3]{Mat16}. For the proof of Corollary~\ref{intro:productgraphs}, see \S~\ref{ss:k-graphs} (Corollary~\ref{cor:ProdGraphs}). Corollary~\ref{intro:deg} also applies to topological full groups of groupoids arising from one vertex higher rank graphs (see Corollary~\ref{cor:OneVertex}). Note that we are not allowing arbitrary graphs in Corollary~\ref{intro:productgraphs} because of the condition that $\# \mfv \fC_j \mfv \geq 2$. A corresponding condition is required in the case of one vertex higher rank graphs. However, our results do cover all products of groupoids from irreducible one-sided shifts of finite type, as we explain in the proof of Corollary~\ref{cor:ProdGraphs}.
\pars

Next, we consider Zappa-Sz{\'e}p products of a small category and a groupoid. Such examples have been considered in various situations in \cite{Nek04,Nek18a,EP,Sta,LRRW14,LRRW18,LY,ABRW,BKQS,Wit19} and, in particular, give rise to R{\"o}ver-Nekrashevych groups as mentioned above. For our purposes, Zappa-Sz{\'e}p products allow us to adjoin invertible elements to a small category while keeping the key properties which are needed to establish finiteness properties for the topological full groups of the corresponding groupoid models.
\begin{introcor}
\label{intro:ZS}
Assume that $\fC$ and $\fS$ are as in Theorem~\ref{intro:Gars}. Let $\fG \curvearrowright \fC$ be a self-similar action and form $\fD \defeq \fC \bowtie \fG$. Assume that (Inv) holds and that $\fD$ is right cancellative up to $=^*$. Further suppose that conditions (St), (LCM) and ($\bmt < \bmd$) are satisfied for $\fC$ and that condition (F) is satisfied for $\fD$.
\pari

Then for all natural numbers $n$, $\bmF(I_l(\fD) \ltimes \partial \Omega_{\fD})$ is of type ${\rm F}_n$ if $\fD^*(\mfw,\mfw)$ is of type ${\rm F}_n$ for all $\mfw \in \fD^0$.
\end{introcor}
\nopar

Self-similar actions and all relevant related notions are introduced in \S~\ref{ss:ZS}, where the reader will also find the proof of Corollary~\ref{intro:ZS} (see Theorem~\ref{thm:ZS}). As explained in Example~\ref{ex:ZS}, Corollary~\ref{intro:ZS} generalizes \cite[Theorem~4.15]{SWZ} and covers self-similar groups as in \cite{Nek04,Nek18a}, self-similar actions on graphs as in \cite{EP,LRRW14,LRRW18} as well as self-similar actions on higher rank graphs as in \cite{LY,ABRW}. For instance, Corollary~\ref{intro:ZS} covers topological full groups of groupoid models for Katsura algebras as discussed in \cite[\S~18]{EP}. 
\pars

\section{Category of bisections}
\label{s:CatBisec}

First, we develop a general framework which allows us to describe topological full groups of classes of groupoids as isotropy groups of small categories.

\subsection{The setting of general groupoids}
\label{ss:GenGPD}

Let $\cG$ be an {\'e}tale groupoid with range and source maps $\rmr, \rms$. Let $X \defeq \cG^{(0)}$ denote the unit space of $\cG$, and assume that $X$ is totally disconnected, compact and Hausdorff. Let $\cC \cO$ be the set of compact open subsets of $X$. We set out to define the category of bisections of $\cG$.

\bdefin
Define a small category $\bfB$ as follows: 
\setlength{\parindent}{0.5cm} \setlength{\parskip}{0cm}

Objects of $\bfB$ are given by finite tuples $\bm{U} = (U_i)_{i \in I}$ with $U_i \in \cC \cO$ for all $i \in I$. 

Morphisms of $\bfB$ from $\bm{U} = (U_i)_{i \in I}$ to $\bm{V} = (V_j)_{j \in J}$ are of the form $\bm{a} = (a_i)_{i \in I}$, together with a map $I \to J, \, i \ma j(i)$, where for each $i \in I$, $a_i$ is a compact open bisection of $\cG$ with $\rms(a_i) = U_i$ and $\rmr(a_i) \subseteq V_{j(i)}$. Moreover, for each $j \in J$, if we set $I_j \defeq \menge{i \in I}{j(i) = j}$, then we require that $V_j = \coprod_{i \in I_j} r(a_i)$. For such a morphism $\bm{a}$, set $\bm{d}(\bm{a}) = \bm{U}$ and $\bm{t}(\bm{a}) = \bm{V}$, and the map $I \to J$ is called the base map. Let $\bfB(\bm{V},\bm{U})$ denote the set of morphisms with $\bm{d}(\bm{a}) = \bm{U}$ and $\bm{t}(\bm{a}) = \bm{V}$.

Let us now define composition in $\bfB$. Given $\bm{b} = (b_j) \in \bfB(\bm{W},\bm{V})$ and $\bm{a} = (a_i) \in \bfB(\bm{V},\bm{U})$, where $\bm{U} = (U_i)_{i \in I}$, $\bm{V} = (V_j)_{j \in J}$ and $\bm{W} = (W_k)_{k \in K}$, define $\bm{b} \bm{a} = ((b a)_i)_{i \in I} \in \bfB(\bm{W},\bm{U})$, where $(b a)_i \defeq b_{j(i)} a_i$ for each $i \in I$, and the base map for $\bm{b} \bm{a}$ is given by the composition $I \to J \to K, \, i \ma j(i) \ma k(j(i))$. 
\edefin
\nopar

Here and in the following, our index sets $I$ are finite ordered sets of the form $\gekl{1, \dotsc, m}$.
\setlength{\parindent}{0cm} \setlength{\parskip}{0.5cm}

To see that composition in $\bfB$ is well defined, observe that $\rms((b a)_i) = U_i$ and $\rmr((b a)_i) \subseteq W_{k(j(i))}$. Moreover, given $k \in K$, if we set $I_k \defeq \menge{i \in I}{k(j(i)) = k}$, then 
$$
 \coprod_{i \in I_k} r((b a)_i) = \coprod_{j \in J_k} \coprod_{i \in I_j} r(b_j a_i) = \coprod_{j \in J_k} r(b_j) = W_k.
$$

\blemma
\label{lem:bfBLeftRev}
$\bfB$ is left reversible, i.e., for all $\bm{a}, \bm{b} \in \bfB$ with $\bm{t}(\bm{a}) = \bm{t}(\bm{b})$, there exist $\bm{a}', \bm{b}' \in \bfB$ such that $\bm{a} \bm{a}' = \bm{b} \bm{b}'$ in $\bfB$. 
\elemma
\nopar

\bproof
Suppose that $\bm{a} = (a_h)_{h \in H}$ and $\bm{b} = (b_i)_{i \in I}$ with $\bm{t}(\bm{a}) = \bm{V} = \bm{t}(\bm{b})$, where $\bm{V} = (V_j)_{j \in J}$, and assume that $H \to J, \, h \ma j_{\bm{a}}(h)$ and $I \to J, \, i \ma j_{\bm{b}}(i)$ are the base maps for $\bm{a}$ and $\bm{b}$, respectively. First, we arrange by right multiplication by $\bm{a}', \bm{b}' \in \bfB$ that there exists a bijection $I \cong H, \, i \ma h_i$ such that $\rmr(a_{h_i}) = \rmr(b_i)$. Let $H \bullet I \defeq \menge{(h,i) \in H \times I}{\rmr(a_h) \cap \rmr(b_i) \neq \emptyset}$. For $(h,i) \in H \bullet I$, define $\dot{a}_{h,i} = \rms \left( (\rmr(a_h) \cap \rmr(b_i)) a_h \right)$ and $\dot{b}_{h,i} = \rms \left( (\rmr(a_h) \cap \rmr(b_i)) b_i \right)$. We obtain morphisms $\dot{\bm{a}} \defeq (\dot{a}_{h,i})_{(h,i) \in H \bullet I} \in \bfB(\bm{d}(\bm{a}),(\dot{a}_{h,i})_{(h,i) \in H \bullet I})$ with base map $H \bullet I \to H, \, (h,i) \ma h$ and $\dot{\bm{b}} \defeq (\dot{b}_{h,i})_{(h,i) \in H \bullet I} \in \bfB(\bm{d}(\bm{b}),(\dot{b}_{h,i})_{(h,i) \in H \bullet I})$ with base map $H \bullet I \to I, \, (h,i) \ma i$. We have $\bm{a} \dot{\bm{a}} = (a_h \dot{a}_{h,i})_{(h,i) \in H \bullet I} \in \bfB(\bm{V},\dot{a}_{h,i})_{(h,i) \in H \bullet I})$ with base map $H \bullet I \to H, \, (h,i) \ma j_{\bm{a}}(h)$ and $\bm{b} \dot{\bm{b}} = (b_i \dot{b}_{h,i})_{(h,i) \in H \bullet I} \in \bfB(\bm{V},\dot{b}_{h,i})_{(h,i) \in H \bullet I})$ with base map $H \bullet I \to H, \, (h,i) \ma j_{\bm{b}}(i)$. Then we obtain for each $(h,i) \in H \bullet I$ that $\rmr(a_h \dot{a}_{h,i}) = \rmr(a_h) \cap \rmr(b_i) = \rmr(b_i \dot{b}_{h,i})$. So without loss of generality, we may assume that $H = I$ and that $\rmr(a_i) = \rmr(b_i)$ for all $i$. In that situation, let $\bm{a}' = (a'_i)_{i \in I} \in \bfB((\rms(a_i))_{i \in I}, (\rms(b_i))_{i \in I})$ with base map $\id_I$ and $a'_i = a_i^{-1} b_i$. Then $\bm{a} \bm{a}' = \bm{b}$.
\eproof
\pars

Let us now describe the enveloping groupoid of $\bfB$. 
\blemma
The enveloping groupoid $Q(\bfB)$ of $\bfB$ is given as follows: 
\setlength{\parindent}{0.5cm} \setlength{\parskip}{0cm}

The set of objects of $Q(\bfB)$ coincides with the set of objects of $\bfB$. 

Morphisms of $Q(\bfB)$ from $\bm{U} = (U_i)_{i \in I}$ to $\bm{V} = (V_j)_{j \in J}$ are of the form $\bm{a} = (a_l)_{l \in L}$, where $L$ is a subset of $J \times I$, where for each $(j,i) \in L$, $a_{j,i}$ is a compact open bisection of $\cG$ with $\rms(a_{j,i}) \subseteq U_i$ and $\rmr(a_{j,i}) \subseteq V_j$. Moreover, for each $i \in I$, we require that $U_i = \coprod_{(j,i) \in L} \rms(a_{j,i})$, and for each $j \in J$, we require that $V_j = \coprod_{(j,i) \in L} \rmr(a_{j,i})$. 

Composition in $Q(\bfB)$ is given as follows: Suppose that $\bm{a} = (a_l)_{l \in L} \in Q(\bfB)(\bm{V},\bm{U})$ and $\bm{b} = (b_m)_{m \in M} \in Q(\bfB)(\bm{W},\bm{V})$, where $\bm{U} = (U_i)_{i \in I}$, $\bm{V} = (V_j)_{j \in J}$ and $\bm{W} = (W_k)_{k \in K}$. Let $N$ be the set of pairs $(k,i) \in K \times I$ for which there exists $j \in J$ such that $\rms(b_{k,j}) \cap \rmr(a_{j,i}) \neq \emptyset$. For $(k,i) \in N$, define $(b a)_{k,i} \defeq \coprod_j b_{k,j} a_{j,i}$, where the disjoint union is taken over all $j \in J$ with $\rms(b_{k,j}) \cap \rmr(a_{j,i}) \neq \emptyset$. $(b a)_{k,i}$ is again a compact open bisection, and we define the product of $\bm{b}$ and $\bm{a}$ as the morphism $\bm{b} \bm{a} = ((b a)_{k,i})_{(k,i) \in N} \in Q(\bfB)(\bm{W},\bm{U})$. 

We obtain an embedding $\bfB(\bm{V},\bm{U}) \into Q(\bfB)(\bm{V},\bm{U}), \, \bm{a} \ma \bm{a}'$ as follows: For $\bm{a} = (a_i)$, set $L \defeq \menge{(j(i),i)}{i \in I}$, and $a'_{j(i),i} \defeq a_i$, $\bm{a}' \defeq (\bm{a}_l)_{l \in L}$. 
\elemma
\nopar

Note that we are implicitly choosing a fixed bijection $L \cong \gekl{1, \dotsc, \# L}$.

\bproof
It is straightforward to check that $Q(\bfB)$ is well-defined and that the maps $\bfB(\bm{V},\bm{U}) \into Q(\bfB)(\bm{V},\bm{U})$ defined above indeed give rise to an embedding $\bfB \into Q(\bfB)$ of categories.
\pari

Inverses are given in $Q(\bfB)$ as follows: Given $\bm{a} = (a_l)_{l \in L} \in Q(\bfB)(\bm{V},\bm{U})$, let $M = \menge{(i,j) \in I \times J}{(j,i) \in L}$. Then the inverse of $\bm{a}$ is given by $\bm{a}^{-1} = ((a^{-1})_m)_{m \in M} \in Q(\bfB)(\bm{U},\bm{V})$, where $(a^{-1})_{i,j} = (a_{j,i})^{-1}$. 

Finally, to see that $Q(\bfB) = \bfB \bfB^{-1}$, take an arbitrary element $\bm{a} = (a_l)_{l \in L} \in Q(\bfB)(\bm{V},\bm{U})$. Define the elements $\dot{\bm{a}} = (a_l)_{l \in L} \in \bfB(\bm{V},(\rms(a_l))_{l \in L})$ and $\bm{b} = (b_l)_{l \in L} \in \bfB(\bm{U},(\rms(a_l))_{l \in L})$, where $b_l \defeq \rms(a_l)$ for all $l \in L$. Then $\bm{a} = \dot{\bm{a}} \bm{b}^{-1}$ in $Q(\bfB)$.
\eproof
\nopar

We will identify $\bfB$ as a subcategory of $Q(\bfB)$ via the embedding constructed above.  
\setlength{\parindent}{0cm} \setlength{\parskip}{0.5cm}

Recall that the topological full group of $\cG$ is given by the group of compact open bisections $a$ of $\cG$ which are global, in the sense that $\rmr(a) = X = \rms(a)$. The product in the topological full group is given by multiplication of bisections as subsets of the groupoid (see \cite{Nek19}). The following is an immediate consequence of our construction.
\blemma
For every compact open subset $Y \subseteq X$, we have $\bm{F}(\cG_Y^Y) = \bfB(Y,Y)$ and $\bm{F}(\cG_Y^Y) = Q(\bfB)(Y,Y)$. Here $Y$ denotes the object $(Y)$ of $\bfB$ or $Q(\bfB)$ consisting of the single element $Y \in \cC \cO$.
\elemma
\nopar

In the above, $\cG_Y^Y$ denotes the restriction of $\cG$ to $Y$, i.e., $\cG_Y^Y = \menge{\gamma \in \cG}{\rmr(\gamma), \rms(\gamma) \in Y}$.
\pars

\subsection{Groupoids arising from left cancellative small categories}
\label{ss:GPD_LCSC}

We use the same notation as in \cite{Li21a}, where the reader may find more details as well. Let $\fC$ be a small category. We will identify the category with its set of morphisms, again denoted by $\fC$, and write $\fC^0$ for its set of objects, which we view as identity morphisms and hence as a subset of $\fC$. Let $\mfd: \: \fC \to \fC^0$ and $\mft: \: \fC \to \fC^0$ be the domain and target maps. $\fC^*$ denotes the set of invertible elements of $\fC$. We assume that $\fC$ is left cancellative, i.e., for all $c, x, y \in \fC$ with $\mfd(c) = \mft(x) = \mft(y)$, $cx = cy$ implies $x = y$. Note that our convention is the same as the one in \cite{Li21a,Spi20,Wit19}, while it is opposite to the one used in \cite{Deh15}. Also note that $\fC^*$ is denoted by $\fC\reg$ in \cite{Deh15,Wit19}. 

Let $I_l$ be the left inverse hull of $\fC$, i.e., the inverse semigroup of partial bijections of $\fC$ generated by the left multiplication maps $c: \: \mfd(c) \fC \to c \fC, \, x \ma cx$. 
Here and in the sequel, for $c \in \fC$ and $S \subseteq \fC$, we use the notation $cS \defeq \menge{cs}{s \in S, \, \mft(s) = \mfd(c)}$. A general element $s$ of $I_l$ is of the form $s = d_n^{-1} c_n \dotso d_2^{-1} c_2 d_1^{-1} c_1$, where $\mft(d_i) = \mft(c_i)$ and $\mfd(d_i) = \mfd(c_{i+1})$. We denote the domain and image of $s$ by $\dom(s)$ and $\img(s)$. Let $\cJ$ be the semilattice of idempotents in $I_l$. We identify $\cJ$ with the semilattice $\menge{\dom(s)}{s \in I_l}$ of subsets of $\fC$.

The space of characters $\widehat{\cJ}$ is given by the set of non-zero multiplicative maps $\cJ \to \gekl{0,1}$, which send $0 \in \cJ$ to $0 \in \gekl{0,1}$ in case $I_l$ contains $0$. Here multiplication in $\gekl{0,1}$ is the usual one induced by multiplication in $\Rz$. The topology on $\widehat{\cJ}$ is given by point-wise convergence. A basis of compact open sets for the topology of $\widehat{\cJ}$ is given by sets of the form
$$
 \widehat{\cJ}(e;\mff) \defeq \big\lbrace \chi \in \widehat{\cJ} : \: \chi(e) = 1, \, \chi(f) = 0 \ \forall \ f \in \mff \big\rbrace,
$$
where $e \in \cJ$ and $\mff \subseteq \cJ$ is a finite subset. We can always assume that $f \leq e$ for all $f \in \mff$. We will also set $\widehat{\cJ}(e) \defeq \lbrace \chi \in \widehat{\cJ} : \: \chi(e) = 1 \rbrace$.

Consider the subspace $\Omega$ of $\widehat{\cJ}$ consisting of those $\chi \in \widehat{\cJ}$ with the following property: Suppose that $\mfv \in \fC^0$, and we are given $f \in E$, $e_i \in E$ ($1 \leq i \leq n$) with $f \subseteq \mfv \fC$, $e_i \subseteq \mfv \fC$ for all $i$. If $f = \bigcup_{i=1}^n e_i$ as subsets of $\fC$, then $\chi(f) = 1$ implies that $\chi(e_i) = 1$ for some $1 \leq i \leq n$ (compare \cite[Definition~6.1 -- 6.3]{Spi20}).

Now every $s \in I_l$ induces the partial bijection
$$
 \big\lbrace \chi \in \widehat{\cJ} : \: \chi(\dom(s)) = 1 \big\rbrace \isom \big\lbrace \chi \in \widehat{\cJ} : \: \chi(\ran(s)) = 1 \big\rbrace, \, \chi \ma s.\chi \defeq \chi(s^{-1} \sqcup s).
$$
By restriction, this yields an action $I_l \curvearrowright \Omega$. Now we can form the groupoid
$$
 I_l \ltimes \Omega \defeq \menge{(s,\chi) \in I_l \times \Omega}{\chi(s^{-1}s) = 1} / {}_{\sim}
$$
where we define $(s,\chi) \sim (t,\eta)$ if $\chi = \eta$ and there exists $e \in \cJ$ with $\chi(e) = 1$ such that $se = te$. The topology of $I_l \ltimes \Omega$ is described in \cite[\S~2.2]{Li21a}.

Similarly, if $X \subseteq \Omega$ is a closed invariant subspace, then we consider the groupoid $I_l \ltimes X = (I_l \ltimes \Omega)_X^X$. Here \an{invariant} always refers to the $I_l$-action.

Let us now introduce a special closed invariant subspace called the boundary. First of all, let $\widehat{\cJ}_{\max}$ be the set of characters $\chi \in \widehat{\cJ}$ for which $\chi^{-1}(1)$ is maximal among all characters $\chi \in \widehat{\cJ}$. As observed in \cite[\S~2]{Li21a}, we have $\widehat{\cJ}_{\max} \subseteq \Omega$. The boundary $\partial \Omega$ is given as the closure of $\widehat{\cJ}_{\max}$ in $\Omega$, i.e., $\partial \Omega \defeq \overline{\widehat{\cJ}_{\max}} \subseteq \Omega$. The boundary $\partial \Omega$ is $I_l$-invariant, so that we may form the boundary groupoid $I_l \ltimes \partial \Omega$.

Let us recall the notion of finite alignment (see \cite[Definition~3.2]{Spi20}), which will allow us to focus on principal right ideals as opposed to general elements of $\cJ$. $\fC$ is called finitely aligned if for all $a,b \in \fC$, there exists a finite subset $\gekl{c_i} \subseteq \fC$ such that $a \fC \cap b \fC = \bigcup_i c_i \fC$. Alternatively, this can be rephrased using the notion of mcms (minimal common right multiples), in the sense of \cite[Definition~2.38]{Deh15}. Recall that given $a, b, c \in \fC$, $c$ is called a mcm of $a$ and $b$ if $c \in a \fC \cap b \fC$ and no proper left divisor $d$ (i.e., an element $d \in \fC$ with $c \in d \fC$) satisfies $d \in a \fC \cap b \fC$. It follows from \cite[Lemma~3.3]{Spi20} that $\fC$ is finitely aligned if and only if for all $a, b \in \fC$, the set of mcms $\mcm(a,b)$ is non-empty and finite up to right multiplication by $\fC^*$. Given a subset $A \subseteq \fC$, we introduce the notation $\mcm(A) \defeq \menge{\mcm(a_1,a_2)}{a_1, a_2 \in A}$. Given subsets $A_1, A_2, \dotsc \subseteq \fC$, we set $\mcm(A_1, A_2, \dotsc) \defeq \mcm(\bigcup_i A_i)$.

We briefly explain the connection to C*-algebras and refer to \cite{Li21a} for details. Form the Hilbert space $\ell^2 \fC$, with canonical orthonormal basis given by $\delta_x(y) = 1$ if $x=y$ and $\delta_x(y) = 0$ if $x \neq y$. For each $c \in \fC$, the mapping 
$$
 \delta_x \mapsto
 \bfa
 \delta_{cx} & {\rm if} \ \mft(x) = \mfd(c),\\
 0 & {\rm else},
 \efa
$$
extends to a bounded linear operator on $\ell^2(\fC)$ which we denote by $\lambda_c$. Note that left cancellation is needed at this point to ensure boundedness of the extension of this mapping. Now the left reduced C*-algebra of $\fC$ is given by $C^*_{\lambda}(\fC) \defeq C^*(\menge{\lambda_c}{c \in \fC}) \subseteq \cL(\ell^2 \fC)$. It is shown in \cite[\S~11]{Spi20} that there is a canonical isomorphism $C^*_r(I_l \ltimes \Omega) \isom C^*_{\lambda}(\fC)$ if $\fC$ is finitely aligned.

$\Omega$ is compact if and only if $\fC^0$ is finite. Since we would like our groupoids to have compact unit spaces, we will later on assume that $\fC^0$ is finite.

Now fix a closed invariant subspace $X$ of $\Omega$. Let $\bfB$ be the category of compact open bisections of the groupoid $I_l \ltimes X$. Given $c \in \fC$, the compact open bisection $\menge{[c,\chi]}{\chi \in \Omega, \, \chi(\mfd(c) \fC) = 1}$ will be denoted by $c$ again.  
\bdefin
Let $\bfC \subseteq \bfB$ be the sub-category with the same set of objects and morphisms of the form $\bm{a} = (a_i U_i)_{i \in I} \in \bfB(\bm{V}, \bm{U})$, where $\bm{U} = (U_i)_{i \in I}$ and $a_i \in \fC$.
\edefin
\nopar

It is easy to see that this is indeed a sub-category.
\pars

Let us now explain the connection to topological full groups. In the following, the enveloping groupoid of $\bfC$ is denoted by $Q(\bfC)$.
\blemma
\label{lem:QbfC=bmF}
For every compact open subset $Y \subseteq X$, we have $Q(\bfC)(Y,Y) = \bm{F}( (I_l \ltimes X)_Y^Y )$.
\elemma
\nopar

\bproof
We have $Q(\bfC)(Y,Y) \subseteq Q(\bfB)(Y,Y) = \bm{F}( (I_l \ltimes X)_Y^Y )$. It remains to show \an{$\supseteq$}. Every element of $\bm{F}((I_l \ltimes X)_Y^Y)$ is a finite union of compact open bisections of the form $[s,U]$ for some $s \in I_l$ and $U \subseteq Y$. It is straightforward to see that we an always arrange this to be a disjoint union. Hence every element of $\bm{F}((I_l \ltimes X)_Y^Y)$ is of the form $\bmc_{\bm{V}} \bmc \bmc_{\bm{U}}^{-1}$, where $\bm{U} = (U_i)_{i \in I}$ with $Y = \coprod_{i \in I} U_i$, $\bmc_{\bm{U}} = (U_i)_{i \in I} \in \bfC(Y,\bm{U})$, $\bm{V} = (V_i)_{i \in I}$ with $Y = \coprod_{i \in I} V_i$, $\bmc_{\bm{V}} = (V_i)_{i \in I} \in \bfC(Y,\bm{V})$, and $\bmc = ([s_i, U_i])_{i \in I} \in Q(\bfC)(\bm{V},\bm{U})$ for some $s_i \in I_l$ with $\rmr([s_i, U_i]) = V_i$, where the base map is given by $\id_I$. Write $s_i = b_{i,N}^{-1} a_{i,N} \dotso b_{i,2}^{-1} a_{i,2} b_{i,1}^{-1} a_{i,1}$. Note that, by choosing $N$ big enough, the same $N$ works for all $i$. Then $[s_i, U_i] = \ti{b}_{i,N}^{-1} \ti{a}_{i,N} \dotso \ti{b}_{i,2}^{-1} \ti{a}_{i,2} \ti{b}_{i,1}^{-1} \ti{a}_{i,1}$, where $\ti{a}_{i,1} = [a_{i,1}, U_i] = a_{i,1} U_i$, $\ti{b}_{i,1} = [b_{i,1}, b_{i,1}^{-1} a_{i,1}. U_i] = b_{i,1} (b_{i,1}^{-1} a_{i,1}. U_i)$, .... Now define $\bma_1 \defeq (\ti{a}_{i,1})_{i \in I} \in \bfC$, $\bmb_1 \defeq (\ti{b}_{i,1})_{i \in I} \in \bfC$, ..., $\bma_N \defeq (\ti{a}_{i,N})_{i \in I} \in \bfC$, $\bmb_N \defeq (\ti{b}_{i,N})_{i \in I} \in \bfC$, where the base map is always $\id_I$. Then 
$$
 \bmc = \bmb_N^{-1} \bma_N \dotso \bmb_1^{-1} \bma_1
$$
lies in $Q(\bfC)$. Thus $\bmc_{\bm{V}} \bmc \bmc_{\bm{U}}^{-1}$ lies in $Q(\bfC)$.
\eproof
\pars

For a closed invariant subspace $X \subseteq \Omega$, we write $X(e;\mff) \defeq X \cap \widehat{\cJ}(e;\mff)$ and $X(e) \defeq X \cap \widehat{\cJ}(e)$.

\blemma
\label{lem:QbfCYY=QbfCYvYv}
We have $Q(\bfC)(Y,Y) \cong Q(\bfC)( (Y \cap X(\mfv \fC))_{\mfv \in \fC^0}, (Y \cap X(\mfv \fC))_{\mfv \in \fC^0} )$ for every compact open subset $Y \subseteq X$.
\elemma
\nopar

\bproof
Consider the morphism $\bmc = (c_{\mfv})_{\mfv} \in \bfC(Y,(Y \cap X(\mfv \fC))_{\mfv \in \fC^0})$ given by $c_{\mfv} = Y \cap X(\mfv \fC)$. Then the desired isomorphism is given by $\bmc \sqcup \bmc^{-1}: \: Q(\bfC)( (Y \cap X(\mfv \fC))_{\mfv \in \fC^0}, (Y \cap X(\mfv \fC))_{\mfv \in \fC^0} ) \isom Q(\bfC)(Y,Y)$.
\eproof
\pars

\bprop
\label{prop:bfCLeftRev}
If $\fC$ is finitely aligned, then $\bfC$ is left reversible.
\eprop
\nopar

\bproof
Suppose that $\bma, \bmb \in \bfC$ with $\bm{t}(\bma) = \bm{t}(\bmb)$. As we have seen in the proof of Lemma~\ref{lem:bfBLeftRev}, up to right multiplication with elements of $\bfC$, we may assume that $\bma = (a_i U_i)_{i \in I}$ and $\bmb = (b_i U_i)_{i \in I}$, where $\rmr(a_i U_i) = \rmr(b_i U_i)$ for all $i$. Now fix $i \in I$ and set $a \defeq a_i$, $b \defeq b_i$, $U \defeq U_i$ and set $V \defeq \rmr(a U) = \rmr(b U)$. As $\fC$ is finitely aligned, we can find a finite set $\gekl{c_n}_{1 \leq n \leq N} \subseteq \fC$ such that $a \fC \cap b \fC = \bigcup_n c_n \fC$. Write $c_n = a a_n = b b_n$. For $\chi \in V$, $\chi(a \fC) = 1 = \chi(b \fC)$ implies that there exists $n$ with $\chi(c_n \fC) = 1$. Define
$$
 V_n \defeq \menge{\chi \in V}{\chi(c_n \fC) = 1, \, \chi(c_m \fC) = 0 \text{ for all } m < n}.
$$
We have $\rms(a U) = a^{-1}.V = \coprod_n a^{-1}.V_n$ and $\rms(b U) = b^{-1}.V = \coprod_n b^{-1}.V_n$. Now we claim that $a^{-1}.V_n = a_n.c_n^{-1}.V_n$. Indeed, \an{$\supseteq$} is clear, and for \an{$\subseteq$}, take $\chi \in a^{-1}.V_n$. Then $a.\chi \in V_n$, so $a.\chi(a a_n \fC) = a.\chi(c_n \fC) = 1$, which implies $\chi(a_n \fC) = 1$ and hence $\chi = a_n.a_n^{-1}.\chi$. We conclude that
$$
 \chi = a_n.a_n^{-1}.\chi = a_n.a_n^{-1}.a^{-1}.(a.\chi) = a_n.c_n^{-1}.(a.\chi) \in a_n.c_n^{-1}.V.
$$
Similarly, $b^{-1}.V_n = b_n.c_n^{-1}.V_n$. Now set $\bma'_n \defeq a_n (c_n^{-1}.V)$ and $\bma'_i \defeq (\bma'_n)_n \in \bfC(\rms(\bma),(c_n^{-1}.V_n)_n)$, and $\bmb'_n \defeq b_n (c_n^{-1}.V)$ and $\bmb'_i \defeq (\bmb'_n)_n \in \bfC(\rms(\bmb),(c_n^{-1}.V_n)_n)$. Then $b b_n = a a_n$ for all $n$ implies that $(a U) \bma'_i = (b U) \bmb'_i$, where we view $a U$ as a morphism in $\bfC(\rmr(a U),\rms(a U))$ and $b U$ as a morphism in $\bfC(\rmr(b U),\rms(b U))$. Applying this construction to $a_i U_i$ for all $i \in I$, we obtain $\bma' \in \bfC$ with $\bmd(\bma) = \bmt(\bma')$ and $\bmb' \in \bfC$ with $\bmd(\bmb) = \bm{t}(\bmb')$ such that $\bma \bma' = \bmb \bmb'$.
\eproof
\pars

\bcor
The enveloping groupoid $Q(\bfC)$ is given by $\bfC \bfC^{-1}$ (formed inside $Q(\bfB)$, as constructed in \S~\ref{ss:GenGPD}).
\ecor

\subsection{Decompositions of compact open subsets}

In preparation for later applications, we show how to decompose compact open subsets into basic ones.
 
The following is stated in \cite[Lemma~4.1]{Li17}.
\blemma
\label{lem:CO=-;-}
Every compact open subset of $X$ is a disjoint union of finite sets of the form $X(e;\mff)$.
\elemma
\nopar

\bproof
Let $\cU$ be the collection of subsets of $X$ which can be written as disjoint unions of finite sets of the form $X(e;\mff)$. Our proof will be complete once we show that $\cU$ is closed under finite unions, because every compact open subset of $X$ is a finite union of sets of the form $X(e;\mff)$. For $e, g \subseteq v \fC$ ($v \in \fC^0$), we have $X(e;\mff) \cap X(g;\mfh) = X(eg;\mff \cup \mfh)$. Moreover, $X(eg;\mff \cup \mfh)^{\rm c} =  X(v \fC;eg) \cup \bigcup_{x \in \mff \cup \mfh}  X(x)$. So
$$
  X(g;\mfh) \cap (X(eg;\mff \cup \mfh))^{\rm c} = X(g; \gekl{eg} \cup \mfh) \amalg \bigcup_{x \in \mff \cup \mfh} X(xg;\mfh) = X(g; \gekl{eg} \cup \mfh) \amalg \bigcup_{x \in \mff} X(xg;\mfh),
$$
and
$$
 \bigcup_{x \in \mff} X(xg;\mfh) = X(x_0g;\mfh) \amalg \bigcup_{x \in \mff \setminus \gekl{x_0}} X(xg; \gekl{x_0 g} \cup \mfh).
$$
Proceeding inductively on $\# \mff$, we obtain that $\bigcup_{x \in \mff} X(xg;\mfh) \in \cU$. Hence 
$$
 X(e;\mff) \cup X(g;\mfh) = X(e;\mff) \amalg X(g;\mfh) \cap (X(e;\mff) \cap X(g;\mfh))^{\rm c} = X(e;\mff) \amalg X(g;\mfh) \cap (X(eg;\mff \cup \mfh))^{\rm c}
$$
lies in $\cU$, as desired.
\eproof
\pars

\bremark
\label{rem:bigcup=coprod}
The last computation shows that, given finite subsets $\mff, \mfh \subseteq \cJ$, we have
$$
 \bigcup_{f \in \mff} X(f;\mfh) = \coprod_{f \in \mff} X(f;\mfh_f)
$$
for some finite subsets $\mfh_f \subseteq \cJ$.
\eremark

Now suppose that $\fC$ is finitely aligned. Then $\cJ \setminus \gekl{\emptyset} = \menge{\bigcup_i a_i \fC}{a_i \in \fC}$ by \cite[Corollary~3.8]{Spi20}. By construction of $\Omega$, we have $\Omega(\bigcup_i a_i \fC) = \bigcup_i \Omega(a_i \fC)$ and
$$
 \Omega(\bigcup_i a_i \fC; \big\{ \bigcup_j \epsilon_j^k \fC \big\}_k) = \bigcup_i \Omega(a_i \fC; \big\{ \epsilon_j^k \fC \big\}_{j,k}) = \bigcup_i \Omega(a_i \fC; \mfe_i)
$$
for some finite subsets $\mfe_i \subseteq \menge{x \fC}{x \in \fC}$ by Remark~\ref{rem:bigcup=coprod}. A similar statement holds for $X$ in place of $\Omega$.

In the following, when there is no danger of confusion, we denote $a \fC$ by $a$, for $a \in \fC$. Then the above shows the following:
\nopar

\blemma
Suppose that $\fC$ is finitely aligned. Then every compact open subset of $X$ is a finite disjoint union of sets of the form $X(a;\mfe)$, where $a \in \fC$ and $\mfe$ is a finite subset of $\fC$. 
\elemma
\pars

\section{Categories of bisections attached to Garside categories}
\label{s:CatBisGar}

In this section, we restrict the set of objects of the category $\bfC$ of bisections from \S~\ref{ss:GPD_LCSC}. From now on, we assume that $\fC$ is a finitely aligned left cancellative small category. Finite alignment allows us to focus on principal ideals. As before, when the meaning is clear from the context, given $a \in \fC$, we will denote the principal ideal $a \fC$ by $a$ again. Since $\fC$ is finitely aligned, \cite[Lemma~3.3]{Spi20} implies that $\fC$ admits mcms. 

Let us introduce the following notation as in \cite{Li21a}: Given $a, b \in \fC$, we write $a \preceq b$ if $a$ is a left divisor of $b$, i.e., $b \in a \fC$. We write $a \prec b$ if $b \fC \subsetneq a \fC$. We write $a \tipreceq b$ if $a$ is a right divisor of $b$, i.e., $b \in \fC a$. We write $a \tiprec b$ if $\fC b \subsetneq \fC a$. We write $a =^* b$ if $a \in b \fC^*$ (which is equivalent to $a \fC = b \fC$). 

$\fC$ is right Noetherian in the sense of \cite[Chapter~II, Definition~2.26]{Deh15} if there is no infinite descending chain $\dotso \tiprec a_3 \tiprec a_2 \tiprec a_1$ in $\fC$.

From now on, let as assume that $\fC^0$ is finite, so that $\Omega$ is compact, and let $X$ be a closed invariant subspace of $\Omega$.

We introduce the following terminology and record a related observation for later use.
\nopar

\bdefin
$\fC$ is called right cancellative up to $=^*$ if, for all $a, b, x \in \fC$, $ax=bx$ implies $a =^* b$.
\edefin
\pars

\blemma
\label{lem:rNrc-->rN}
\begin{enumerate}
\item[(i)] If $\fC$ is right cancellative up to $=^*$, then, for all $a, b \in \fC$ and $\chi \in \Omega$, $[a,\chi] = [b,\chi]$ implies $a = ^* b$. In particular, for all $u \in \fC$ and $\chi \in \Omega$, $[u,\chi] = \chi$ implies $u \in \fC^*$.
\item[(ii)] If $\fC$ is right Noetherian and right cancellative up to $=^*$, then the category of bisections $\bfC$ from \S~\ref{ss:GPD_LCSC} is right Noetherian.
\end{enumerate}
\elemma
\nopar

\bproof
(i) By definition of the equivalence relation, $[a,\chi] = [b,\chi]$ implies that there exists $x \in \fC$ with $\chi(x) = 1$ and $ax=bx$. Right cancellation up to $=^*$ implies $a =^* b$, as desired.
\pars

(ii) Now suppose we have $\bma \tipreceq \bmb$ in $\bfC$, i.e., $\bmb = \bma' \bma$. If $\bmt(\bmb) = (V_j)_{j \in J}$ and $\bmt(\bma) = (W_k)_{k \in K}$, then $\# J \leq \# K$. If $\bmd(\bmb) = (U_i)_{i \in I}$, then $\# J \leq \# K \leq \# I$. Hence given a chain $\dotso \tipreceq \bma_3 \tipreceq \bma_2 \tipreceq \bma_1$ in $\bfC$, we may assume that $\bmt(\bma_n) = (V_{n,j})_{j \in J_n}$ for all $n$, and $\# J_n = \# J_1$ for all $n$. Now write $\bma'_n \bma_n = \bma_{n-1}$ and $\bma'_n = (a'_{n,j} V_{n,j})_j$. Let $\bmd(\bma_1) = (U_i)_{i \in I}$ and write $\bma_n = (a_{n,i} U_i)_{i \in I}$. Then $[a_{1,i},U_i] = [a'_{2,j(i)} a_{2,i}, U_i]$ implies by (i) that $a_{1,i} =^* a'_{2,j(i)} a_{2,i}$ and thus $a_{2,i} \tipreceq a_{1,i}$. Continuing this way, we obtain $\dotso \tipreceq a_{3,i} \tipreceq a_{2,i} \tipreceq a_{1,i}$, for all $i$. Since $\fC$ is right Noetherian, we deduce that there must exist $n_0$ such that for all $n \geq n_0$, we have $a'_{n,j} \in \fC^*$ for all $j$ and thus $\bma'_n \in \bfC^*$.
\eproof
\pars

Now let us specialize to the class of Garside categories. The idea behind the concept of Garside structures originated from the study of Braid groups and monoids, and of the more general Artin-Tits groups and monoids. Roughly speaking, Garside structures axiomatize the structures which are needed to carry over classical results and methods from Braid groups and monoids to more general groups, monoids or small categories. One important feature of Garside categories is that elements admit normal forms. This was used in \cite{Li21a} to derive general results about groupoids and C*-algebras attached to left regular representations of Garside categories. In \cite{Wit19}, criteria have been established for finiteness properties of isotropy groups of Garside categories, based on the existence of head functions given by maximal left divisors from a given Garside family. Now our goal is to show that a Garside family in $\fC$ induces Garside families in subcategories of $\bfC$, and based on this, to derive criteria when these Garside structures fit into the general framework developed in \cite{Wit19} which then allows us to establish finiteness properties for isotropy groups and hence topological full groups.

To explain the notion of Garside family, we need some terminology. First of all, a finite sequence $s_1, s_2, \dotsc$ in $\fC$ is called a path if $\mfd(s_k) = \mft(s_{k+1})$ for all $k$. Such a path will be denoted by $s_1 s_2 \dotsm$.
\nopar

\bdefin
A subset $\fS \subseteq \fS$ is closed under right comultiples if for all $r, s \in \fS$ and $a \in \fC$ with $r \preceq a$, $s \preceq a$, there exists $t \in \fS$ with $r \preceq t$, $s \preceq t$ and $t \preceq a$.
\edefin

\bdefin
Suppose $\fS \subseteq \fC$ is closed under right comultiples such that $\fS \cup \fC^*$ generates $\fC$ and $\fS^{\sharp} \defeq \fS \fC^* \cup \fC^*$ is closed under right divisors. 
\pari

A path $s_1, \dotsc, s_l \in \fS^{\sharp}$ is called normal if for all $1 \leq k \leq l-1$ and $r \in \fS$, if $r \preceq s_k s_{k+1}$ then $r \preceq s_k$.

For $a \in \fC$, a normal decomposition or normal form of $a$ is given by a normal path $s_1 \dotsc s_l$ in $\fS^{\sharp}$ with $a = s_1 \dotsm s_l$.

$\fS$ is called a Garside family if every element in $\fC$ admits a normal decomposition.
\edefin
\pars

Suppose that $\fS$ is a Garside family of $\fC$. We can always arrange (see \cite[Chapter~III, Corollary~1.34]{Deh15}) that $\fS$ is $=^*$-transverse, i.e., for all $s_1, s_2 \in \fS$, $s_1 =^* s_2$ implies $s_1 = s_2$. We can also assume without loss of generality that $\fS \cap \fC^* = \emptyset$ (see \cite[Chapter~III, Corollary~1.34]{Deh15}).

Given $a \in \fS$, we define $\Vert a \Vert \defeq 0$ if $a \in \fC^*$ and $\Vert a \Vert \defeq l$ if $s_1 \dotso s_l$ is the unique normal decomposition of $a \in \fC \setminus \fC^*$ (which exists by \cite[Chapter~III, Corollary~1.27]{Deh15}).

If $\fC$ is right Noetherian and admits mcms, then $\fS \subseteq \fC$ is a Garside family if and only if $\fS \cup \fC^*$ generates $\fC$ and $\fS^{\sharp}$ is closed under mcms and right divisors (see \cite[Chapter~IV, Proposition~2.25]{Deh15}). The reader may consult \cite{Deh15} for more information about Garside families.

Now assume that $\fS$ is locally finite (i.e., $\# \mfv \fS < \infty$ for all $\mfv \in \fC^0$). 

\cite[Chapter~III, Corollary~1.37]{Deh15} implies that for all $L \geq 1$, $\fS^{\leq L} = \fS \cup \fS^2 \cup \dotso \cup \fS^L$ is a Garside family again. \cite[Chapter~IV, Proposition~2.25]{Deh15} implies that $\fS^{\leq L}$ is closed under mcms. 
\blemma
\label{lem:SL_No_Inf_Chain}
There are no infinite chains $a \prec a_1 \prec a_2 \prec \dotso$ in $(\fS^{\leq L})^{\sharp}$.
\elemma
\nopar

\bproof
Since $\fS$ is locally finite, $\fS^l$ is again locally finite for all $1 \leq l \leq L$. 
\eproof
\pars

In addition, we assume that for all $L \geq 1$, $(\fS^{\leq L})^{\sharp}$ is closed under left divisors. This is for instance the case if the assumptions of \cite[Chapter~III, Proposition~1.62]{Deh15} are satisfied. 

\blemma
Let $a$ and $b$ be elements of $\fC$ with $\Vert a \Vert = L = \Vert b \Vert$. Then every $c \in \mcm(a,b)$ satisfies $\Vert c \Vert = L$. 
\elemma
\nopar

\bproof
As $\fS^{\leq L}$ is closed under mcms, we must have $\Vert c \Vert \leq L$. If $\Vert c \Vert \leq L-1$, then, since $(\fS^{\leq L-1})^{\sharp}$ is closed under left divisors, it would follow that $\Vert a \Vert \leq L-1$ and $\Vert b \Vert \leq L-1$. But this is a contradiction.
\eproof
\pars

\blemma
\label{lem:aB=aiBi}
Given $a_i \in \fC$ and a finite set $\mfe \subseteq \fC$, there exist finite subsets $\mfe_i \subseteq \mcm(\gekl{a_i},\mfe)$ such that 
$$
  \bigcup_i X(a_i; \mfe) = \coprod_i X(a_i; \mfe_i).
$$
\elemma
\nopar

\bproof
We have
$$
 \bigcup_i X(a_i; \mfe) = \bigcup_i X(a_i; \mcm(a_i, \mfe)) = X(a_1; \mcm(a_1, \mfe)) \amalg \bigcup_{i > 1} X(a_i; \mcm(a_i,a_1), \mcm(a_i, \mfe)).
$$
Now proceed inductively on $\# \gekl{a_i}$.
\eproof
\pars

\blemma
\label{lem:aBaai}
Every compact open subset of the form $X(a;\mfe)$ is a finite disjoint union of sets of the form $X(b;b \mff)$ for some $b \in \fC$ and some finite set $\mff \subseteq \mfd(b) \fS$.
\elemma
\nopar

\bproof
Write $\mfe = \gekl{a s_i \epsilon_i}$ for some $s_i \in \fS$. Choose $L$ big enough so that $\mfe \subseteq (\fS^{\leq L})^{\sharp}$. Then
$$
 X(a;\mfe) = X(a;\gekl{a s_i}) \amalg \bigcup_i X(a s_i; \mfe).
$$
By Lemma~\ref{lem:aB=aiBi}, we can further decompose
$$
 \bigcup_i X(a s_i; \mfe) = \coprod_i X(a s_i; \mfe_i),
$$
where $\mfe_i \subseteq \mcm(a s_i, \mfe)$ are finite subsets of $(\fS^{\leq L})^{\sharp}$ since $a s_i \in (\fS^{\leq L})^{\sharp}$ and $\mfe \subseteq (\fS^{\leq L})^{\sharp}$ because $(\fS^{\leq L})^{\sharp}$ is closed under left divisors and mcms.
\pars

Without loss of generality we may assume that $X(a s_i; \mfe_i) \neq \emptyset$. So if $\mfe_i \neq \emptyset$, then we must have $a s_i \prec \epsilon$ for all $\epsilon \in \mfe_i$. If there exists $i$ with $\mfe_i \neq \emptyset$, then let $a^{(1)} \defeq a$, $a^{(2)} \defeq a s_i$, and continue the above process with $a s_i$ in place of $a$. In this way, we obtain a chain $a^{(1)} \prec a^{(2)} \prec a^{(3)} \prec \dotso$ in $(\fS^{\leq L})^{\sharp}$. As observed above, every such chain in $(\fS^{\leq L})^{\sharp}$ must be finite. Hence this process must end. This means that at some point, we obtain $\mfe_i = \emptyset$ for all $i$.
\eproof
\pars

\bdefin
Given objects $\bmU = (U_i)_{i \in I}$ and $\bmU' = (U'_{i'})_{i' \in I'}$ in $\bfC$, let $K = I \amalg I'$ and define the object $\bmU \amalg \bmU' \defeq (\ti{U}_k)_{k \in K}$ of $\bfC$ by $\ti{U}_k \defeq U_k$ if $k \in I$ and $\ti{U}_k \defeq U'_k$ if $k \in I'$.
\edefin

\bdefin
We define $\bfX_{\fS}$ as the smallest subset of objects of $\bfC$ which contains $\menge{X(\mfv;\mfe)}{\mfv \in \fC^0, \, \mfe \subseteq \mfv \fS}$ and which is closed under the operation $\amalg$.
\edefin

\bdefin
We set
$$
 \bfC_{\bfX_{\fS}} \defeq \menge{\bma \in \bfC}{\bmt(\bma), \, \bmd(\bma) \in \bfX_{\fS}}.
$$
\edefin

\bprop
\label{prop:QdXX=F}
Suppose that $\fC$ is finitely aligned and that $\fC^0$ is finite. Assume that $\fS$ is a Garside family in $\fC$ which is locally finite, $=^*$-transverse with $\fS \cap \fC^* = \emptyset$. Moreover assume that for all $L \geq 1$, $(\fS^{\leq L})^{\sharp}$ is closed under left divisors.
\pari

For all $\bm{U} = (U_i)_{i \in I}$, where $U_i$ are compact open subsets of $X$, there exists $\bmc \in \bfC$ with $\bmt(\bmc) = \bm{U}$ and $\bmd(\bmc) \in \bfX_{\fS}$.

$\bfC_{\bfX_{\fS}}$ is left reversible. 

The enveloping groupoid $Q(\bfC_{\bfX_{\fS}})$ of $\bfC_{\bfX_{\fS}}$ satisfies $Q(\bfC_{\bfX_{\fS}})(*,*) = Q(\bfC)(*,*)$ for all $* \in \bfX_{\fS}$.
\eprop
\nopar

\bproof
For all $i \in I$, Lemma~\ref{lem:CO=-;-} and Lemma~\ref{lem:aBaai} imply that there exists a decomposition $U_i = \coprod_k X(c_{ik};c_{ik} \mfe_{ik})$ with $\mfe_{ik} \subseteq \fS$. Then $\bmc = (c_{ik} X(\mfd(c_{ik});\mfe_{ik}))_{(i,k)}$ defines a morphism in $\bfC(\bm{U}, (X(\mfd(c_{ik});\mfe_{ik}))_{(i,k)})$, where the base map sends $(i,k)$ to $i$. 
\pari

To see that $\bfC_{\bfX_{\fS}}$ is left reversible, take $\bma, \bmb \in \bfC_{\bfX_{\fS}}$. By Proposition~\ref{prop:bfCLeftRev}, there exist $\bma', \bmb' \in \bfC$ with $\bma \bma' = \bmb \bmb'$. We have just seen that there exists $\bmc \in \bfC$ with $\bmt(\bmc) = \bmd(\bma') = \bmd(\bmb')$ and $\bmd(\bmc) \in \bfX_{\fS}$. Hence $\bma' \bmc, \bmb' \bmc \in \bfC_{\bfX_{\fS}}$ and $\bma (\bma' \bmc) = \bmb (\bmb' \bmc)$. 

Finally, every element of $Q(\bfC)(*,*)$ is of the form $\bma \bmb^{-1}$ for some $\bma, \bmb \in \bfC$. As we have seen above, there exists $\bmc \in \bfC$ with $\bmt(\bmc) = \bmd(\bma) = \bmd(\bmb)$ and $\bmd(\bmc) \in \bfX_{\fS}$. Then $\bma \bmb^{-1} = (\bma \bmc) (\bmc \bmb)^{-1} \in Q(\bfC_{\bfX_{\fS}})(*,*)$.
\eproof
\pars

\section{Garside subcategories describing topological full groups}
\label{s:GarsCat-TFG}

First, we identify criteria when subcategories of $\bfC$ or $\bfC_{\bfX_{\fS}}$ still allow for descriptions of topological full groups as isotropy groups (of enveloping groupoids). Then we turn to the question when such subcategories admit right Garside maps.

Throughout this section, let $\fC$ be a finitely aligned left cancellative small category with finite $\fC^0$. Assume that $\fS$ is a Garside family in $\fC$ which is locally finite, $=^*$-transverse with $\fS \cap \fC^* = \emptyset$. Moreover assume that for all $L \geq 1$, $(\fS^{\leq L})^{\sharp}$ is closed under left divisors.

As $\fC$ is generated by $\fS$ and $\fC^*$, we know that $\fC \setminus \fC^* = \fS \fC$. It follows that, since every $s \in \fS$ is contained in some $a \fC$, where $a \in \fA$, we must have $\fC \setminus \fC^* = \fA \fC$. Here $\fA$ denotes the set of atoms in $\fC$, i.e., elements of $\fC$ which do not admit proper divisors. For $x \in \fC$, let $\chi_x$ be the element of $\Omega$ with $\chi_x(e) = 1$ if and only if $x \fC \subseteq e$ (as defined above). Now suppose $\mfv \in \fC^0$. Our assumption that $\fS$ is locally finite implies that $\mfv \fA$ must be finite. Since we have $\gekl{\chi_{\mfv}} = \Omega(\mfv;\mfv\fA)$, it follows that $\gekl{\chi_{\mfv}}$ is clopen. Hence $\menge{\chi_x}{x \in \fC}$ is open in $\Omega$, so that $\Omega_{\infty} \defeq \Omega \setminus \menge{\chi_x}{x \in \fC}$ is closed. It is easy to see that $\Omega_{\infty}$ is invariant as well. From now on, assume that $X \subseteq \Omega_{\infty}$ is a closed invariant subspace of $\Omega_{\infty}$.

\subsection{Topological full groups as isotropy groups of subcategories}

Now suppose that $\mfv \in \fC^0$ and $\mfe \subseteq \mfv \fS$. Let $\mfs \subseteq \mfv \fC \setminus \fC^*$ be closed under mcms, i.e., $\mcm(s,t) \subseteq \mfs$ for all $s, t \in \mfs$. Further assume that the canonical projection $\fC \isom \fC/{ }_{=^*}$ restricts to an injective map on $\mfs$. 

\bdefin
For $s \in \gekl{\mfv} \cup \mfs$, let $\mff_s$ be a $=^*$-transverse subset of $\fC$ consisting up to $=^*$ of all minimal elements of
$$
 \menge{f \in \fC \setminus \fC^*}{sf \in \mfs} \cup \menge{f \in \fC \setminus \fC^*}{sf \in \mcm(s,\mfe)}.
$$
\edefin
Here \an{minimal} refers to $\prec$. Note that $\mff_s$ is only well-defined up to $=^*$. But since we will only use $\mff_s$ to construct $X(\mfd(s); \mff_s)$, different choices of $\mff_s$ will lead to the same set $X(\mfd(s); \mff_s)$. Later on, we will only need the case where $\mff_s \subseteq \fS$, so that $\mff_s$ will be automatically $=^*$-transverse.

For $L \geq 1$, let $\fS_L \subseteq \mfv \fC$ be such that the canonical projection $\fC \isom \fC/{ }_{=^*}$ restricts to a bijection $\fS_L \isom \menge{a \in \fC}{\Vert a \Vert = L} / { }_{=^*}$. 
\blemma
\label{lem:gamma}
\begin{enumerate}
\item[(i)] $X(\mfv; \mfe) \cap (\bigcup_{s \in \mfs} sX) = \coprod_{s \in \gekl{\mfv} \cup \mfs} s.X(\mfd(s); \mff_s)$.
\item[(ii)] $\gamma(\mfe, \mfs) \defeq (sX(\mfd(s); \mff_s))_{s \in \gekl{\mfv} \cup \mfs}$ defines an element of $\bfC(X(\mfv; \mfe), (X(\mfd(s); \mff_s))_{s \in \gekl{\mfv} \cup \mfs})$.
\item[(iii)] For $\mfs = \fS_L$, we have $\bigcup_{s \in \fS_L} sX = X$ and $\mff_s \subseteq \fA$, so that we can arrange $\mff_s \subseteq \fS$, and $\mff_s$ coincides with the set of minimal elements of $\menge{f \in \fS}{sf \in \mfs}$. In particular, $\gamma(\mfe, \fS_L) \in \bfC_{\bfX_{\fS}}$.
\end{enumerate}
\elemma
\nopar

\bproof
(i) Given $s, s' \in \mfs$, write $\mcm(s,s') = \gekl{s s_i} = \{s' s'_j\}$. We have $\bigcup_i s s_i \subseteq \bigcup s \mff_s$ unless $\mcm(s,s') = \gekl{s}$. Similarly, $\bigcup_j s' s'_j \subseteq \bigcup s' \mff_{s'}$ unless $\mcm(s,s') = \gekl{s'}$. So $X(s; s \mff_s)$ and $X(s'; s' \mff_{s'})$ are disjoint unless $s =^* s'$ (which implies $s = s'$ by assumption on $\mfs$).
\pari

Moreover, by construction, we have $s X(\mfd(s); \mff_s) \subseteq X(s; \mcm(s, \mfe)) = X(\mfv; \mfe) \cap sX$. Given $\chi \in X(\mfv; \mfe) \cap (\bigcup_{s \in \mfs} sX)$, we have $\menge{s \in \mfs}{\chi(s) = 1} \neq \emptyset$. Hence there exists $t \in \mfs$ maximal with respect to $\prec$ such that $\chi(t) = 1$. By maximality of $t$ and since $\chi(\mfe) = 0$ implies $\chi(tf) = 0$ for all $f \in \fC$ with $tf \in \mcm(t,\mfe)$, we deduce that $\chi \in t X(\mfd(t); \mff_t)$.
\pars

(ii) The bisection of $\gamma(\mfe, \mfs)$ corresponding to $\mfv$ is given by $X(\mfv; \mfs \cup \mfe)$, which is disjoint from $sX$ for all $s \in \mfs$, hence from $s X(\mfd(s); \mff_s)$ for all $s \in \mfs$. Moreover, (i) implies that
$$
 X(\mfv; \mfs \cup \mfe) \amalg (\coprod_{s \in \mfs} X(\mfd(s); \mff_s)) = X(\mfv; \mfs) \cap X(\mfv; \mfe) \amalg X(\mfv; \mfe) \cap (\bigcup_{s \in \mfs} sX) 
 = X(\mfv; \mfe) \cap (X(\mfv; \mfs) \cup (\bigcup_{s \in \mfs} sX))
 = X(\mfv; \mfe).
$$

(iii) We have $\bigcup_{s \in \fS_L} sX = X$ because $X \subseteq \Omega_{\infty}$. Now take $s \in \fS_L$ and $f \in \fC$ such that $sf \in \fS_L$. It then follows that given any left divisor $f'$ of $f$, we must have $sf' \in \fS_L$ up to $=^*$. Indeed, $sf' \in (\fS^{\leq L})^{\sharp}$ as $(\fS^{\leq L})^{\sharp}$ is closed under left divisors. If $sf' \in (\fS^{\leq L-1})^{\sharp}$, then we would deduce $s \in (\fS^{\leq L-1})^{\sharp}$, which would contradict $s \in \fS_L$. Hence $\Vert sf \Vert = L$. This shows that the minimal elements of $\menge{f \in \fC \setminus \fC^*}{sf \in \mfs}$ are indeed contained in $\fA$. Moreover, for every $s \in \fS_L$, $\mcm(s,\mfe) \subseteq \fS_L$ because $\mfe \subseteq \fS$ implies $\mcm(s,\mfe) \subseteq (\fS^{\leq L})^{\sharp}$, and no element of $\mcm(s,\mfe)$ can lie in $(\fS^{\leq L-1})^{\sharp}$ as $s$ lies in $\fS_L$. It follows that $\menge{f \in \fC \setminus \fC^*}{sf \in \mcm(s,\mfe)} \subseteq \menge{f \in \fC \setminus \fC^*}{sf \in \fS_L}$, as desired.
\eproof

\blemma
\label{lem:char-gammaS}
Let $\gamma = (s X(\mfd(s); \mff'_s))_{s \in \fS_L} \in \bfC_{\bfX_{\fS}}(X(\mfv; \mfe), (X(\mfd(s);\mff'_s))_{s \in \fS_L})$ with $\mff'_s \subseteq \fS$ such that $\Vert s f' \Vert = L$ for all $f' \in \mff'_s$. Then $\gamma = \gamma(\mfe, \fS_L)$.
\elemma
\nopar

\bproof
Our assumption implies $\mff'_s \subseteq \mff_s$ for all $s \in \fS_L$, which in turn yields $X(\mfd(s);\mff'_s) \supseteq X(\mfd(s);\mff_s)$ and hence $X(s; s \mff'_s) \supseteq X(s; s \mff_s)$. Let $\ti{s}$ be maximal with respect to $\prec$ such that $X(\mfd(\ti{s});\mff'_{\ti{s}}) \supsetneq X(\mfd(\ti{s});\mff_{\ti{s}})$, which implies $\mff'_{\ti{s}} \subsetneq \mff_{\ti{s}}$, and $X(\mfd(t);\mff'_t) = X(\mfd(t);\mff_t)$ for all $\ti{s} \prec t$. Take $f \in \mff_{\ti{s}} \setminus \mff'_{\ti{s}}$ such that $X(\mfd(\ti{s});\mff'_{\ti{s}}) \cap f X \neq \emptyset$, which implies $X(\ti{s}; \ti{s} \mff'_{\ti{s}}) \cap \ti{s} f X \neq \emptyset$. At the same time,
$$
 \ti{s} f X = \coprod_{\ti{s} f \preceq t} X(t; t \mff_t). 
$$
This is a contradiction because $X(\ti{s}; \ti{s} \mff'_{\ti{s}})$ and $\coprod_{\ti{s} f \preceq t} X(t; t \mff_t) = \coprod_{\ti{s} f \preceq t} X(t; t \mff'_t)$ are disjoint.
\eproof
\pars

\blemma
\label{lem:gammaS=gammaS'}
Given $\fS_L$, $\fS'_L$ such that the canonical projection restricts to bijections $\fS_L \isom \menge{a \in \fC}{\Vert a \Vert = L} / { }_{=^*}$ and $\fS'_L \isom \menge{a \in \fC}{\Vert a \Vert = L} / { }_{=^*}$, there exists $\bmu \in \bfC_{\bfX_{\fS}}^*$ with $\bmu \in \bfC_{\bfX_{\fS}}(\bmd(\gamma(\mfe, \fS_L)), \bmd(\gamma(\mfe, \fS'_L)))$ such that $\gamma(\mfe, \fS_L) \bmu = \gamma(\mfe, \fS'_L)$.
\elemma
\nopar

\bproof
Suppose that we have $s =^* s' \in \fS'_L$ for all $s \in \fS_L$, say $s' = s u_{s'}$ for $u_{s'} \in \fC^*$. Then $s \fC = s' \fC$ and $\menge{s f \fC}{f \in \mff_s} = \menge{s' f' \fC}{f' \in \mff_{s'}} = \menge{s u_{s'} f' \fC}{f' \in \mff_{s'}}$. It follows that $\menge{f \fC}{f \in \mff_s} = u_{s'} \menge{f' \fC}{f' \in \mff_{s'}}$ and thus $u_{s'}. X(\mfd(s'); \mff_{s'}) = X(\mfd(s); \mff_s)$. Define $\bmu = (u_{s'} X(\mfd(s'); \mff_{s'}))_{s' \in \fS'_L} \in \bfC_{\bfX_{\fS}}(\bmd(\gamma(\mfe, \fS_L)), \bmd(\gamma(\mfe, \fS'_L)))$. This morphism has the desired properties.
\eproof
\pars

We need the following operations.
\nopar
\bdefin
Suppose that $\bma = (a_i)_{i \in I} \in \bfC_{\bfX_{\fS}}((V_j)_{j \in J}, (U_i)_{i \in I})$ and $\bma' = (a'_{i'})_{i' \in I'} \in \bfC_{\bfX_{\fS}}((V'_{j'})_{j' \in J'}, (U'_{i'})_{i' \in I'})$. Let $K = I \amalg I'$, $L = J \amalg J'$ (equipped with some ordering) and define $\bma \amalg \bma' \defeq (a \amalg a')_k)_{k \in K} \in \bfC_{\bfX_{\fS}}((\ti{V}_l)_{l \in L}, (\ti{U}_k)_{k \in K})$ by $\ti{U}_k \defeq U_k$ if $k \in I$, $\ti{U}_k \defeq U'_k$ if $k \in I'$, $\ti{V}_l \defeq V_l$ if $l \in J$, $\ti{V}_l \defeq V'_l$ if $l \in J'$, $(a \amalg a')_k \defeq a_k$ if $k \in I$, $(a \amalg a')_k \defeq a'_k$ if $k \in I'$, and the base map $K \to L$ is given as the disjoint union of the base maps of $\bma$ and $\bma'$.
\pari

Given $\bfX \subseteq \bfX_{\fS}$, a subset $\bGamma$ of $\bfC_{\bfX}$ is called $\amalg_{\fX}$-closed if for all $\bma \in \bGamma$ and $x \in \bfX$, $\bma \amalg x$ lies in $\bGamma$.
\edefin
\pars

Now consider $\fX \subseteq \menge{X(\mfv;\mfe)}{\mfv \in \fC^0, \, \mfe \subseteq \mfv \fS}$ and let $\bfX$ be the smallest $\amalg$-closed subset of $\bfX_{\fS}$ containing $\fX$. Further let $\Gamma$ be a subset of $\menge{\bma \in \bfC_{\bfX}}{\bmt(\bma) \in \fX}$ and $\bGamma$ the smallest $\amalg_{\fX}$-closed subset of $\bfC_{\bfX}$ containing $\Gamma$. Let $\Pi$ be the smallest subcategory of $\bfC$ containing $\bGamma$ and $\bfC_{\bfX}^*$. Since both $\bGamma$ and $\bfC_{\bfX}^*$ are closed under the operation $\amalg$, so is $\Pi$.

\bdefin
We say that condition (St) is satisfied if the following conditions hold.
\nopar

\begin{enumerate}
\item[(1$_{\fX}$)] For all $u \in \fC^*$, $\mfe \subseteq \mfd(u) \fS$ such that $X(\mfd(u); \mfe) \in \fX$, we have $X(\mft(u); u \mfe) \in \fX$.
\item[(2$_{\fX}$)] For all $\mfv \in \fC^0$ and $\mfe \subseteq \mfv \fS$, there exists $\bma \in \bfC$ with $\bmt(\bma) = X(\mfv; \mfe)$ and $\bmd(\bma) \in \fX$.
\item[(1$_{\Gamma}$)] For all $\mfv \in \fC^0$ and $\mfe \subseteq \mfv \fS$ with $X(\mfv; \mfe) \in \fX$, we have $\gamma(\mfe, \fS) \in \Gamma$.
\item[(2$_{\Gamma}$)] For all $L \geq 1$, $\mfv \in \fC^0$, $c \in \fC \mfv$ with $\Vert c \Vert = L$, $\mfs \defeq \menge{s \in \mfv \fS}{\Vert cs \Vert = L}$ and $\mfe \subseteq \mfv \fS$ such that $X(\mfv; \mfe) \in \fX$, we have $\gamma(\mfe, \mfs) \in \Gamma$.
\end{enumerate}
\edefin
\nopar

Note that (1$_{\Gamma}$) includes the condition that $\bmd(\gamma(\mfe, \fS)) \in \bfX$ and (2$_{\Gamma}$) includes the condition that $\bmd(\gamma(\mfe, \mfs)) \in \bfX$ for those $\mfs$ which appear in (2$_{\Gamma}$).
\pars

\bprop
\label{prop:ab=gamma}
Assume that (1$_{\fX}$), (1$_{\Gamma}$) and (2$_{\Gamma}$) are satisfied. Let 
$$
 \bma = (a_i X(\mfd(a_i); \mff_i))_{i \in I} \in \bfC_{\bfX}(X(\mfv; \mfe), (X(\mfd(a_i); \mff_i))_{i \in I}),
$$
and suppose that $L \geq 1$ is such that
\begin{equation}
\label{e:norm<=L}
 \Vert a_i f \Vert \leq L \text{ for all } i \in I \text{ and } f \in \mff_i.
\end{equation}
Then there exist $\fS_L \subseteq \fC$ such that the canonical projection restricts to a bijection $\fS_L \isom \menge{a \in \fC}{\Vert a \Vert = L} / { }_{=^*}$ and $\bmb \in \Pi$ with $\bma \bmb = \gamma(\mfe, \fS_L)$.
\eprop
\nopar

\bproof
Fix $\fS_L \subseteq \fC$ such that the canonical projection restricts to a bijection $\fS_L \isom \menge{a \in \fC}{\Vert a \Vert = L} / { }_{=^*}$. If $\bma \neq \gamma(\mfe, \fS_L)$, then by right multiplication with elements in $\C$, we can arrange
\begin{enumerate}
\item[(i)] $\Vert a_i \Vert = L$ for all $i \in I$;
\item[(ii)] There exists $\fS'_L \subseteq \fC$ such that the canonical projection restricts to a bijection $\fS'_L \isom \menge{a \in \fC}{\Vert a \Vert = L} / { }_{=^*}$, and $\gekl{a_i} \subseteq \fS'_L$;
\item[(iii)] $\mff_i = \mff_{a_i}$ for all $i \in I$;
\item[(iv)] $\gekl{a_i} = \fS'_L$ (where $\fS'_L$ is from (ii)).
\end{enumerate}
\pars

(i) Assume that there exists $a_j \in \gekl{c_i}$ with $\Vert a_j \Vert < L$. Consider the morphism $\bmb \defeq \coprod_{i \in I} \bmb_i$ with $\bmb_i = X(\bmd(a_i); \mff_i)$ for $i \neq j$ and $\bmb_j = \gamma(\mff_j, \fS)$. Then $\bma \bmb$ still satisfies \eqref{e:norm<=L}. This process must stop after finitely many steps, so that we indeed can arrange $\Vert a_i \Vert = L$ for all $i \in I$.

(ii) Assume that there exist $a_k, a_l \in \gekl{a_i}$ such that $a_k =^* a_l$, say $a_l = a_k u$. Because of $u \fS \subseteq \fS \fC^*$, there exists $\mff_l \subseteq \fS$ so that $u \mff_l =^* \mfg_l$ (i.e., for all $f \in \mff_l$ and $g \in \mfg_l$, we have $uf =^* g$). It follows that $u^{-1}.X(\mfd(a_k); \mfg_l) = X(\mfd(a_l); \mff_l)$. Define $\bmb \defeq (b_i)_{i \in I} \in \bfC_{\bfX}((X(\mfd(a_i); \mff_i))_{i \in I}, (U_i)_{i \in I})$, where $U_i \defeq X(\mfd(a_i); \mff_i)$ for $i \neq l$ and $U_l \defeq X(\mfd(a_k); \mfg_l)$, $b_i \defeq X(\mfd(a_i); \mff_i)$ for $i \neq l$ and $b_l \defeq u^{-1} X(\mfd(a_k); \mfg_l)$, and the base map is given by $\id_I$. Note that we need condition (1$_{\fX}$) to ensure that $X(\mfd(a_k); \mfg_l) \in \fX$. Then by multiplying $\bma$ from the right with $\bmb$, we can replace $a_l$ by $a_k = a_l u^{-1}$. Continuing this way, we can make sure that for all $a_k, a_l \in \gekl{a_i}$, $a_k =^* a_l$ implies $a_l = a_k$.

(iii) Assume that there exists $j \in I$ with $\mff_j \neq \mff_{a_j}$. Define
$$
 \mfs \defeq \menge{s \in \fS}{\Vert a_j s \Vert = L}.
$$
Now define $\bmb \defeq \coprod_{i \in I} \bmb_i$, where $\bmb_i \defeq X(\mfd(a_i); \mff_i)$ for $i \neq j$ and $\bmb_j \defeq \gamma(\mff_j, \mfs)$. By multiplying $\bma$ from the right with $\bmb$, we replace $\mff_j$ by $\mff_{a_j}$ because $\gamma(\mff_j, \mfs)$ contains the compact open bisection $X(\mfd(a_j); \mff_{a_j})$.

(iv) Take $\fS_L$ as in (ii). Assume that there is $s \in \fS_L$ with $s \neq^* a_i$ for all $i \in I$. We know that $X(s; s \mff_s)$ is disjoint to $X(a_i; a_i \mff_i)$ for all $i \in I$ by Lemma~\ref{lem:gamma}~(i), and that $X(\mfv; \mfe) = \coprod_{i \in I} X(a_i; a_i \mff_i)$ since $\bma$ is a morphism. It follows that $X(s; s \mff_s) = \emptyset$ and thus $X(\mfd(s); \mff_s) = \emptyset$.

Having gone through (i) -- (iv), we arrive at $\gamma(\mfe, \fS_L)$ because of Lemmas~\ref{lem:char-gammaS} and \ref{lem:gammaS=gammaS'}.
\eproof

Now fix an object $* \in \bfX$. We consider the component of $*$ in $\Pi$ as follows:
\nopar

\bdefin
\label{def:component}
Let $\bfX(*) \defeq \menge{x \in \bfX}{\exists \ \alpha \in \Pi \text{ with } \bmt(\alpha) = *, \, \bmd(\alpha) = x}$, $\fX(*) \defeq \menge{U \in \fX}{\exists \ x \in \bfX(*) \text{ with } U \in x}$ and $\Gamma(*) \defeq \menge{\gamma \in \Gamma}{\bmt(\gamma) \in \fX(*)}$. The component of $*$ in $\Pi$ is given by $\C \defeq \menge{\alpha \in \Pi}{\bmt(\alpha) \in \bfX(*)}$. Let $\Q$ denote the enveloping groupoid of $\C$.
\edefin
Clearly, $\C$ is a subcategory of $\Pi$. We set out to study properties of $\C$. 
\pars

\bcor
\label{cor:ab=gamma}
Assume that (St) holds. Then $\gamma(\mfe, \fS_L) \in \Pi$ for all $\mfv \in \fC^0$, $\mfe \subseteq \mfv \fS$ with $X(\mfv; \mfe) \in \fX$ and all positive integers $L$. Moreover, $\Pi$ and $\C$ are left reversible, and we have $\Q(*,*) = Q(\bfC)(*,*)$. 
\ecor
\nopar

\bproof
The first claim is a direct consequence of Proposition~\ref{prop:ab=gamma}. For the second claim, take $\alpha, \alpha' \in \Pi$ with $\bmt(\alpha) = \bmt(\alpha') = \coprod_{i \in I} X(\mfv_i;\mfe_i)$. By Proposition~\ref{prop:ab=gamma}, there are $\beta, \beta' \in \Pi$ such that $\alpha \beta = \coprod_{i \in I} \gamma(\mfe_i,\fS_L) = \alpha' \beta'$ for sufficiently big $L$. If $\alpha$ and $\alpha'$ both lie in $\C$, then $\beta$ and $\beta'$ lie in $\C$ as well. To prove the last claim, we know from Lemma~\ref{prop:QdXX=F} that every element of $Q(\bfC)(*,*)$ is of the form $\bma \bmb^{-1}$ for some $\bma, \bmb \in \bfC_{\bfX_{\fS}}$ with $\bmt(a) = * = \bmt(b)$. Suppose that $* = \coprod_i X(\mfv_i;\mfe_i)$. Because of condition (2$_{\fX}$), there exists $\bmc \in \bfC$ with $\bmt(\bmc) = \bmd(\bma) = \bmd(\bmb)$ and $\bmd(\bmc) \in \bfX$. Applying Proposition~\ref{prop:ab=gamma}, we can find $\bm{\xi} \in \Pi$ with $\bma \bmc \bm{\xi} = \coprod_i \gamma(\mfe_i, \fS_L)$ for some $L \geq 1$. Applying Proposition~\ref{prop:ab=gamma} again, we can find $\bm{\eta} \in \C$ with $\bmb \bmc \bm{\xi} \bm{\eta} = \coprod_i \gamma(\mfe_i, \fS_{L'})$ for some $L' \geq 1$. Thus we conclude
\begin{align*}
 \bma \bmb^{-1} &= (\bma \bmc \bm{\xi}) (\bmb \bmc \bm{\xi})^{-1} = (\coprod_i \gamma(\mfe_i, \fS_L)) (\bmb \bmc \bm{\xi})^{-1} = (\coprod_i \gamma(\mfe_i, \fS_L)) \bm{\eta} (\bmb \bmc \bm{\xi} \bm{\eta})^{-1}\\
 &= (\coprod_i \gamma(\mfe_i, \fS_L)) \bm{\eta} (\coprod_i \gamma(\mfe_i, \fS_{L'}))^{-1} \in \Q(*,*).
\end{align*}
\eproof
\pars

\subsection{Existence of least common multiples}

Let $\fX$, $\bfX$, $\Gamma$, $\bGamma$, $\Pi$ and $\C$ be as above.

\bdefin
Given $\alpha, \beta \in \bGamma$, we define $\lcm_{\Gamma \subseteq \bfC}(\alpha,\beta)$ to be an element $\gamma \in \Pi$ such that $\gamma \in \alpha \bGamma \cap \beta \bGamma$ and for all $\chi \in \Pi$, $\chi \in \alpha \bfC \cap \beta \bfC$ implies that $\chi \in \gamma \bfC$.
\edefin

Let us now introduce the following condition:
\bdefin
We say that condition (LCM) holds if the following are satisfied:
\nopar

\begin{enumerate}
\item[(3$_{\Gamma}$)] $\Pi^* \bGamma \subseteq \bGamma \Pi^*$;
\item[(4$_{\Gamma}$)] $\bGamma \Pi^*$ is closed under non-invertible right divisors in $\Pi$, i.e., if $\alpha \omega \in \bGamma \Pi^*$ with $\alpha \in \bGamma$ and $\omega \in \Pi \setminus \Pi^*$, then $\omega \in \bGamma \Pi^*$;
\item[(5$_{\Gamma}$)] $\bfC$ admits lcms (denoted by $\lcm_{\bfC}$) and for all $\alpha, \beta \in \bGamma$, $\lcm_{\bfC}(\alpha, \beta) \in \alpha \bGamma \cap \beta \bGamma$.
\end{enumerate}
\edefin
\pars

\bprop
\label{prop:RightDivInC}
Suppose that $\fC$ is right cancellative up to $=^*$ and right Noetherian. Assume that (3$_{\Gamma}$) and (4$_{\Gamma}$) hold. If $\lcm_{\bGamma \subseteq \bfC}(\alpha,\beta)$ exists for all $\alpha, \beta \in \bGamma$, then given $\chi \in \Pi$ with $\chi = \alpha x_1$ for some $\alpha \in \bGamma$ and $x_1 \in \bfC$, we must have $x_1 \in \Pi$.
\eprop
\nopar

\bproof
Because of condition (3$_{\Gamma}$), $\chi \in \Pi$ implies that $\chi = \beta \chi_1$ for some $\beta \in \bGamma$ and $\chi_1 \in \Pi$. So $\chi = \alpha x_1 = \beta \chi_1$. Let $\gamma = \lcm_{\Gamma \subseteq \bfC}(\alpha,\beta)$. Write $\gamma = \alpha \alpha' = \beta \alpha_1$. Then $\chi = \gamma x_2$.
\pars

Thus $x_1 = \alpha' x_2$, $\chi_1 = \alpha_1 x_2$. Moreover, $\chi_1 \in \Pi$ implies that $\chi_1 = \beta_1 \chi_2$ for some $\beta_1 \in \bGamma$ and $\chi_2 \in \Pi$. Let $\gamma_1 = \lcm_{\Gamma \subseteq \C \subseteq \bfC}(\alpha_1,\beta_1)$. Write $\gamma_1 = \alpha_1 \alpha'_1 = \beta_1 \alpha_2$. Then $x_2 = \alpha'_1 x_3$, $\chi_2 = \alpha_2 x_3$ and $\chi_2 = \beta_2 \chi_3$. 

Continue in this way to obtain $x_i = \alpha'_{i-1} x_{i+1}$, $\chi_i = \alpha_i x_{i+1}$ and $\chi_i = \beta_i \chi_{i+1}$.

So $\chi = \beta \chi_1 = \beta \beta_1 \chi_2 = \beta \beta_1 \beta_2 \chi_3 = \dotso = \beta \beta_1 \beta_2 \dotsm \beta_i \chi_{i+1}$. By Lemma~\ref{lem:rNrc-->rN}, $\bfC$ is right Noetherian. Hence we must arrive at an element $\chi_{i+1}$ with $\chi_{i+1} \in \C^*$, so that $\alpha_i x_{i+1} \in \beta_i \Pi^* \subseteq \bGamma \Pi^*$. Hence $x_{i+1} \in \bfC_{\bfX}^* = \Pi^*$, or condition (4$_{\Gamma}$) implies $x_{i+1} \in \bGamma \Pi^* \subseteq \Pi$. In both cases, we deduce $x_{i+1} \in \Pi$ and thus $x_i \in \Pi$, etc., and finally $x_1 \in \Pi$.
\eproof
\pars

\bcor
In the situation of Proposition~\ref{prop:RightDivInC}, we have $\lcm_{\Gamma \subseteq \bfC}(\alpha,\beta) = \lcm_{\Pi}(\alpha,\beta)$.
\ecor
\nopar

\bproof
Suppose that $\alpha \zeta = \beta \eta$ in $\Pi$. Then, if $\gamma = \lcm_{\Gamma \subseteq \bfC}(\alpha,\beta)$, then $\alpha \zeta = \gamma z$ for some $z \in \bfC$. But then Proposition~\ref{prop:RightDivInC} implies $z \in \Pi$, so that $\alpha \zeta = \beta \eta \in \gamma \Pi$. 
\eproof
\pars

\blemma
Suppose we have condition (3$_{\Gamma}$). Assume that for all $\alpha, \beta \in \bGamma$, $\lcm_{\Pi}(\alpha,\beta)$ exists and $\lcm_{\Pi}(\alpha,\beta) \in \alpha \bGamma \cap \beta \bGamma$. Then $\Pi$ admits lcms.
\elemma
\nopar

\bproof
Because of condition (3$_{\Gamma}$), it suffices to show that lcms exist for all elements of $\Pi$ which are finite products of generators in $\bGamma$. Let $\spkl{\bGamma}$ denote the set of such finite products. For $\alpha \in \spkl{\bGamma}$, let $\ell_{\Gamma}(\alpha)$ be the minimal number of factors in $\bGamma$ needed to express $\alpha$ as a product.
\pars

First of all, we show for all $\sigma \in \bGamma$ and $\chi \in \spkl{\bGamma}$, $\lcm_{\Pi}(\sigma, \chi)$ exists and is of the form $\sigma \eta$ with $\ell_{\Gamma}(\eta) \leq \ell_{\Gamma}(\chi)$. We proceed inductively on $\ell_{\Gamma}(\chi)$. The case $\ell_{\Gamma}(\chi) = 1$ holds by assumption. Now take $\omega \chi \in \spkl{\bGamma}$ with $\omega \in \bGamma$ and $\ell_{\Gamma}(\omega \chi) = 1 + \ell_{\Gamma}(\chi)$. By assumption, $\lcm_{\Pi}(\sigma, \omega) = \omega \zeta = \sigma \tau$ for some $\zeta, \tau \in \bGamma$. Thus we conclude that
$$
 \lcm_{\Pi}(\sigma, \omega \chi) = \lcm_{\Pi}(\lcm_{\Pi}(\sigma, \omega), \omega \chi) = \lcm_{\Pi}(\omega \zeta, \omega \chi) = \omega \lcm_{\Pi}(\zeta, \chi) = \omega \zeta \xi = \sigma \tau \xi,
$$
where $\ell_{\Gamma}(\xi) \leq \ell_{\Gamma}(\chi)$ by induction hypothesis. Hence $\lcm_{\Pi}(\sigma, \omega \chi) = \sigma (\tau \xi)$, and $\ell_{\Gamma}(\tau \xi) \leq 1 + \ell_{\Gamma}(\xi) \leq 1 + \ell_{\Gamma}(\chi) = \ell_{\Gamma}(\omega \chi)$.

Now we prove that $\lcm_{\Pi}(\sigma, \chi)$ exists for all $\sigma, \chi \in \spkl{\bGamma}$ inductively on $\max(\ell_{\Gamma}(\sigma), \ell_{\Gamma}(\chi))$. Write $\sigma = \alpha \beta$, $\chi = \delta \varepsilon$ with $\alpha, \delta \in \bGamma$ and $\ell_{\Gamma}(\beta) < \ell_{\Gamma}(\sigma)$ and $\ell_{\Gamma}(\varepsilon) < \ell_{\Gamma}(\chi)$. Write $\lcm_{\Pi}(\alpha, \delta) = \alpha \alpha' = \delta \delta'$ for some $\alpha', \delta' \in \bGamma$. Then
\begin{align*}
 \lcm_{\Pi}(\alpha \beta, \delta \varepsilon) &= \lcm_{\Pi}(\lcm_{\Pi}(\alpha, \delta), \alpha \beta, \delta \varepsilon) = \lcm_{\Pi}(\lcm_{\Pi}(\alpha \alpha', \alpha \beta), \lcm_{\Pi}(\delta \delta', \delta \varepsilon))\\
 &= \lcm_{\Pi}(\alpha \, \lcm_{\Pi}(\alpha', \beta), \delta \, \lcm_{\Pi}(\delta', \varepsilon))
 = \lcm_{\Pi}(\alpha \alpha' \mu,  \delta \delta' \nu) = \alpha \alpha' \lcm_{\Pi}(\mu, \nu),
\end{align*}
where $\ell_{\Gamma}(\mu) \leq \ell_{\Gamma}(\beta)$ and $\ell_{\Gamma}(\nu) \leq \ell_{\Gamma}(\varepsilon)$ as shown above, so that $\lcm_{\Pi}(\mu, \nu)$ exists by induction hypothesis.
\eproof
\pars

Putting everything together, we arrive at the following conclusion.
\nopar

\bcor
\label{cor:lcm}
Suppose that $\fC$ is right cancellative up to $=^*$ and right Noetherian. Assume that conditions (3$_{\Gamma}$) and (4$_{\Gamma}$) are satisfied. If $\lcm_{\Gamma \subseteq \bfC}(\alpha,\beta)$ exists for all $\alpha, \beta \in \bGamma$, then $\Pi$ admits lcms. In particular, if (5$_{\Gamma}$) holds, then $\Pi$ admits lcms.
\ecor
\pars

Define a map $\Delta$ on $\bfX_{\fS}$ by $\Delta((X(\mfv; \mfe_i)_i) \defeq \coprod_i \gamma(\mfe_i, \fS)$. Write $\Div_{\Pi}(\Delta)$ and $\widetilde{\Div}_{\Pi}(\Delta)$ for the set of all left divisors and right divisors in $\Pi$ of $\Delta(x)$ for all $x \in \bfX$, respectively. Define $\Div_{\C}(\Delta)$ and $\widetilde{\Div}_{\C}(\Delta)$ analogously.

\bprop
\label{prop:GarsideTFG}
Suppose that $\C$ admits lcms. Then $\Div_{\C}(\Delta)$ is closed under lcms, $\widetilde{\Div}_{\C}(\Delta) \subseteq \Div_{\C}(\Delta)$, $\Delta$ restricted to the objects of $\C$ is a right Garside map for $\C$ and $\Div_{\C}(\Delta)$ is a Garside family for $\C$.
\eprop
\nopar

Recall from \cite[Chapter~V, Definition~1.15]{Deh15} that a map $\Delta$ from the objects of $\Pi$ to $\Pi$ is called a right Garside map if $\bmt(\Delta(x)) = x$ for every object $x$, $\Div_{\Pi}(\Delta)$ generates $\Pi$, $\widetilde{\Div}_{\Pi}(\Delta) \subseteq \Div_{\Pi}(\Delta)$ and for every $\bmg \in \Pi$ with $\bmt(\bmg) = x$, the elements $\bmg$ and $\Delta(x)$ admit a left gcd. The reader may consult \cite[Chapter~V]{Deh15} for more information about right Garside maps. 
\pars

Note that, by Corollary~\ref{cor:lcm}, $\Pi$ admits lcms if condition (LCM) holds.
\nopar

\bproof
Corollary~\ref{cor:ab=gamma} implies that $\Delta(x) \in \C$ for all $x \in \bfX(*)$. Given $\bma, \bmb \in \Div_{\C}(\Delta)$ with $x \defeq \bmt(\bma) = \bmt(b)$, we have $\lcm_{\C}(\bma, \bmb) \preceq_{\C} \Delta(x)$, so that $\lcm_{\C}(\bma, \bmb) \in \Div_{\C}(\Delta)$. This shows that $\Div_{\C}(\Delta)$ is closed under lcms. Let $\bma \in \widetilde{\Div}_{\C}(\Delta)$, $\bma = (a_i U_i) \in \bfC((V_j), (U_i))$ where $U_i = X(\mfd(a_i); \mff_i)$. $\bma \in \widetilde{\Div}_{\C}(\Delta)$ implies that $\Vert a_i f_i \Vert \leq 1$ for all $i$ and all $f_i \in \mff_i$. Hence Proposition~\ref{prop:ab=gamma} implies that $\bma \in \Div_{\C}(\Delta)$. This shows $\widetilde{\Div}_{\C}(\Delta) \subseteq \Div_{\C}(\Delta)$, as desired.
\pars

The third claim follows from \cite[Chapter~V, Proposition~1.20]{Deh15} (see \cite[Chapter~V, Definition~1.15]{Deh15}). Indeed, given $\bmg \in \C$ with $\bmt(\bmg) = x$, $\lcm_{\C}(\menge{\bma \in \C}{\bma \preceq \bmg, \, \bma \preceq \Delta(x)} \in \Div_{\C}(\Delta)$ is the left gcd of $\bmg$ and $\Delta(x)$.
\eproof
\pars

Combining Proposition~\ref{prop:GarsideTFG} with Corollary~\ref{cor:lcm}, we obtain the following:
\nopar

\bcor
\label{cor:GarsideTFG}
Suppose that $\fC$ is right cancellative up to $=^*$ and right Noetherian. Assume that condition (LCM) is satisfied. Then $\C$ admits lcms, $\Div_{\C}(\Delta)$ is closed under lcms, $\widetilde{\Div}_{\C}(\Delta) \subseteq \Div_{\C}(\Delta)$, $\Delta$ restricted to the objects of $\C$ is a right Garside map for $\C$ and $\Div_{\C}(\Delta)$ is a Garside family for $\C$.
\pari

If, in addition, condition (St) is satisfied, then $\Q(*,*) = Q(\bfC)(*,*)$.
\ecor 
\pars

Thus we are naturally led to the following question: When does $\bfC$ admit lcms? To discuss a sufficient criterion, we introduce the following condition.
\nopar

\bdefin
We say that condition (F) holds if for all $\mfv \in \fC^0$, $a, b \in \fC^*(\mfv,\mfv)$ and $U = X(\mfv;\mfe) \in \fX$, $[a,U] = [b,U]$ in $\bfC$ implies $a=b$. 
\edefin
Note that (F) holds if, for instance, $\fC$ is right cancellative.
\pars

\blemma
\label{lem:aU=bU}
Suppose that $\fC$ is right cancellative up to $=^*$ and that condition (F) holds. Given $a, b \in \fC$ and $U = X(\mfv;\mfe) \in \fX$, $[a,U] = [b,U]$ in $\bfC$ implies $a=b$.
\elemma
\nopar

\bproof
By Lemma~\ref{lem:rNrc-->rN}~(i), $[a,U] = [b,U]$ implies $a =^* b$. Hence $b = au$ for some $u \in \fC^*$. It follows that $[a,U] = [a,u.U] [u,U]$ and thus $U = u.U$. So $[a,U] = [a,U] [u,U]$ and hence $U = [u,U]$. Condition (F) now implies $u \in \fC^0$, i.e., $a=b$.
\eproof
\pars

Suppose that $\fC$ is right cancellative up to $=^*$, right Noetherian and that condition (F) holds. Further assume that $\fC$ has disjoint mcms, in the sense that for all $a, b \in \fC$, we have $a \fC \cap b \fC = \coprod_k c_k \fC$ for some finite collection $\gekl{c_k} \subseteq \fC$. Note that the case $a \fC \cap b \fC = \emptyset$ is allowed. Our assumption is for instance satisfied if $\fC$ admits lcms, i.e., for all $a, b \in \fC$, we have $a \fC \cap b \fC = \emptyset$ or $a \fC \cap b \fC = c \fC$ for some $c \in \fC$. 

We now construct lcms in $\bfC$. Suppose that we are given $\bma = (a_i U_i)_{i \in I} \in \bfC(W,(U_i))$ and $\bmb = (b_j V_j)_{j \in J} \in \bfC(W,(V_j))$. Write $a_i \fC \cap b_j \fC = \coprod_{k \in K_{ij}} c_k \fC$ and $a_i a'_{ik} = c_k = b_j b'_{jk}$. Set
$$
 O_k \defeq c_k^{-1}. ((a_i.U_i \cap b_j.V_j) \cap c_k.X).
$$
Now define 
$$
 \bma' \defeq (a'_{ik} O_k)_{i \in I, k \in K_{ij}, j \in J} \in \bfC((U_i), (O_k)_{i, k \in K_{ij}, j \in J}),
$$
where the base map is given by $(i,k) \mapsto i$, and
$$
 \bmb' \defeq (b'_{jk} O_k)_{j \in J, k \in K_{ij}, i \in I} \in \bfC((V_j), (O_k)_{j \in J, k \in K_{ij}, i \in I}),
$$
where the base map is given by $(j,k) \mapsto j$.

\blemma
\label{lem:lcm_bfC}
$\bma'$ and $\bmb'$ are well-defined. We have $\bma \bma' = \bmb \bmb' = \lcm_{\bfC}(\bma, \bmb)$.
\pari

In general, given $\bar{\bma} = \coprod \bma_n$ and $\bar{\bmb} = \coprod \bmb_n$ with $\bmt(\bma_n) = \bmt(\bmb_n)$, then $\lcm_{\bfC}(\bar{\bma},\bar{\bmb}) = \coprod \lcm_{\bfC}(\bma_n,\bmb_n)$.
\elemma
\nopar

\bproof
It is straightforward to check that $\bma'$ and $\bmb'$ are well-defined. The equation $\bma \bma' = \bmb \bmb'$ follows from $a_i a'_{ik} = b_j b'_{jk}$.
\pars

Let us now verify the lcm property. Assume that $\alpha, \beta \in \bfC$ satisfy $\bma \alpha = \bmb \beta$. Suppose that $\alpha = (\alpha_l Z_l)_{l \in L}$ and the base map for $\alpha$ sends $l \in L_i \subseteq L$ to $i$, and that $\beta = (\beta_l Z_l)_{l \in L}$ and the base map for $\beta$ sends $l \in L_j \subseteq L$ to $j$. For $l \in L_i \cap L_j$, we conclude that $\alpha_l. Z_l \subseteq U_i$ and $\beta_l. Z_l \subseteq V_j$. Hence $a_i \alpha_l Z_l = b_j \beta_l Z_l$. It follows that $a_i \alpha_l = b_j \beta_l$ because of Lemma~\ref{lem:aU=bU}. So there exists $k \in K_{ij}$ such that $a_i \alpha_l = b_j \beta_l \in c_k \fC = a_i a'_{ik} \fC$. It follows that $\alpha_l = a'_{ik} \zeta_l$ and $\beta_l = b'_{jk} \zeta_l$ for some $\zeta_l \in \fC$. $a_i \alpha_l. Z_l \subseteq a_i.U_i \cap b_j.V_j \cap c_k.X = c_k.O_k$ implies that $\alpha_l.Z_l \subseteq a'_{ik}.O_k$ and thus $\zeta_l.Z_l \subseteq O_k$. The decompositions $U_i = \coprod_{k \in K_{ij}, j \in J} a'_{ik} O_k$ and $U_i = \coprod_{l \in L_i} \alpha_l.Z_l$ imply $O_k = \coprod_{l \in L_k} \zeta_l.Z_l$, where $L_k = \menge{l \in L}{\zeta_l.Z_l \subseteq O_k}$. It also follows that $\coprod_{k \in K_{ij}, i \in I, j \in J} L_k = L$. Hence $\zeta \defeq (\zeta_l Z_l)_{l \in L}$ defines an element in $\bfC((O_k)_{k \in K_{ij}, i \in I, j \in J}, (Z_l)_{l \in L})$, where the base map of $\zeta$ sends $l \in L_k$ to $k$. By construction, we have $\bma' \zeta = \alpha$ and $\bmb' \zeta = \beta$, as desired.

The last claim is straightforward to see.
\eproof
\pars

We record the following consequence.
\nopar

\bcor
\label{cor:5Gamma}
Suppose that $\fC$ is right cancellative up to $=^*$, right Noetherian and has disjoint mcms, and assume that (F) holds. Then condition (5$_{\Gamma}$) holds if for all $\bma, \bmb \in \Gamma$, the elements $\bma'$, $\bmb'$ constructed above lie in $\bGamma$.
\ecor
\pars

Now let us write $\Div_{\bfC_{\bfX_{\fS}}}(\Delta)$ for the set of all left divisors in $\bfC_{\bfX_{\fS}}$ of $\Delta(x)$ for all $x \in \bfX_{\fS}$. We now show that Corollary~\ref{cor:GarsideTFG} applies to $\bfX = \bfX_{\fS}$ and $\bGamma = \bfS$, where
$$
 \bfS = \menge{\bma = (a_i X(\mfd(a_i); \mff_i))_i \in \bfC_{\bfX_{\fS}}}{a_i \in \fS \ \forall i; \, \Vert a_i f_i \Vert \leq 1 \ \forall i, \, f_i \in \mff_i}.
$$
Define $\fX$ and $\Gamma$ correspondingly.

\blemma
\label{lem:maxXG:allCOND}
For $\fX$ and $\Gamma$, conditions (1$_{\fX}$), (2$_{\fX}$), (1$_{\Gamma}$) -- (5$_{\Gamma}$) are satisfied.
\elemma
\nopar

\bproof
Clearly, (1$_{\fX}$) and (2$_{\fX}$) hold. Moreover, Proposition~\ref{prop:ab=gamma} implies that $\bGamma \bfC_{\bfX_{\fS}}^* = \Div_{\bfC_{\bfX_{\fS}}}(\Delta)$. As an immediate consequence, we deduce that (1$_{\Gamma}$) -- (4$_{\Gamma}$) are satisfied.
\pars

Finally, it remains to verify (5$_{\Gamma}$), i.e., for all $\alpha, \beta \in \bGamma$, we have $\lcm_{\bfC}(\alpha, \beta) \in \alpha \bGamma \cap \beta \bGamma$. In other words, we have to show that the elements $\bma'$ and $\bmb'$ constructed before Lemma~\ref{lem:lcm_bfC} (for $\bma = \alpha$ and $\bmb = \beta$) lie in $\bGamma$. Since $\fS$ is closed under mcms and right divisors, it follows that all $a'_{ik}$ and $b'_{jk}$ lie in $\fS$. It remains to show that all $O_k$ lie in $\fX_{\fS}$. Recall that
$$
 O_k = c_k^{-1}. (X(a_i; a_i \mfe_i) \cap X(b_j; b_j \mff_j) \cap c_k. X),
$$
if $\bmd(\alpha) = (X(\mfd(a_i); \mfe_i))_i$ and $\bmd(\beta) = (X(\mfd(b_j); \mff_j))_j$. Now
$$
 X(a_i; a_i \mfe_i) \cap X(b_j; b_j \mff_j) = \coprod_k X(c_k; \mcm(c_k, a_i \mfe_i), \mcm(c_k, b_j \mff_j)).
$$
As $c_k \in \fS$, $a_i \mfe_i, b_j \mff_j \subseteq \fS$, we conclude that $\mcm(c_k, a_i \mfe_i), \mcm(c_k, b_j \mff_j) \subseteq \fS$. Hence $O_k$ lies in $\fX_{\fS}$, as desired.
\eproof
\pars

\bcor
\label{cor:thmA}
Suppose that $\fC$ is a left cancellative small category with finite $\fC^0$. Further assume that $\fC$ is right cancellative up to $=^*$, finitely aligned, right Noetherian and admits disjoint mcms, and that (F) holds. Let $\fS$ be a Garside family in $\fC$ which is locally finite, $=^*$-transverse with $\fS \cap \fC^* = \emptyset$. Moreover assume that for all $L \geq 1$, $(\fS^{\leq L})^{\sharp}$ is closed under left divisors.
\pari

Given $* \in \bfX_{\fS}$, let $\fX$, $\Gamma$ be as above and $\C$, $\Q$ as in Definition~\ref{def:component}. Then $\Delta$ restricted to the objects of $\C$ is a right Garside map for $\C$, $\Div_{\C}(\Delta)$ is a Garside family for $\C$, and $\Q(*,*) = Q(\bfC)(*,*)$. If $* = (X(\mfv;\mfe))_{\mfv \in \fV}$ for some $\fV \in \fC^0$, $Y = \coprod_{\mfv \in \fC^0} Y_{\mfv}$, where $Y_{\mfv} = \emptyset$ if $\mfv \notin \fV$ and $Y_{\mfv} = X(\mfv;\mfe) \in \fX$ if $\mfv \in \fV$, then $\Q(*,*) \cong \bmF((I_l \ltimes X)_Y^Y)$. 
\ecor
\pars

As a consequence, we solve the word problem for $\Q(*,*)$. More precisely, given a word in $\bfS$, $\bfS^{-1}$ and $\C^*$, we will construct an algorithm which decides whether our word represents the trivial element of $\Q(*,*)$. We need a $=^*$-map for $\fS^{\sharp}$, i.e., a partial map $E$ from $\fS^{\sharp} \times \fS^{\sharp}$ to $\fC^*$ with the property that $E(s,t)$ is defined if and only if $s =^* t$, and in that case $E(s,t) = u \in \fC^*$ with $su = t$. Such a $=^*$-map is called computable if it can be implemented on a Turing machine. Our goal is to establish the following:
\nopar

\bcor
\label{cor:WordProblem}
Assume that we are in the same situation as in Corollary~\ref{cor:5Gamma}, Lemma~\ref{lem:maxXG:allCOND} and Corollary~\ref{cor:thmA}. Suppose that there exists a computable $=^*$-map for $\fS^{\sharp}$. Given $* \in \bfX_{\fS}$, let $\C$ and $\Q$ be as in Definition~\ref{def:component}. Then $\Q(*,*) = Q(\bfC)(*,*)$ has decidable word problem. If $* = (X(\mfv;\mfe))_{\mfv \in \fV}$ for some $\fV \in \fC^0$, $Y = \coprod_{\mfv \in \fC^0} Y_{\mfv}$, where $Y_{\mfv} = \emptyset$ if $\mfv \notin \fV$ and $Y_{\mfv} = X(\mfv;\mfe) \in \fX$ if $\mfv \in \fV$, then $\Q(*,*) \cong \bmF((I_l \ltimes X)_Y^Y)$ has decidable word problem.
\ecor
\pars

For the proof, we need the following observation.
\nopar

\blemma
\label{lem:bE-computable}
If there exists a computable $=^*$-map for $\fS^{\sharp}$, then there exists a computable $=^*$-map for $\bfS^{\sharp}$.
\elemma
\nopar

The latter means a partial map $\bmE$ from $\bfS^{\sharp} \times \bfS^{\sharp}$ to $\C^*$ with the property that $\bmE(\alpha,\beta)$ is defined if and only if $\alpha =^* \beta$ in $\C$, and in that case $\bmE(\alpha,\beta) = \bmu \in \C^*$ with $\alpha \bmu = \beta$.
\bproof
Suppose that $\alpha = ([a_i,V_i])_{i \in I}$ and $\beta = ([b_j,W_j])_{j \in J}$. We can only have $\alpha =^* \beta$ if $I$ and $J$ agree up to permutation, so that we may assume $I = J$. Then $\alpha =^* \beta$ if and only if $[a_i,V_i] =^* [b_i,W_i]$ for all $i$. Now $[a_i,V_i] =^* [b_i,W_i]$ holds if and only if $a_i u_i = b_i$ for some $u_i \in \fC^*$ and $u_i.W_i = V_i$. With the help of the given $=^*$-map $E$, we can decide whether the first condition is satisfied, and $E$ also produces the value for $u_i$. For the second condition, suppose that $W_i = X(\mfw;\mff)$ and $V_i = X(\mfv;\mfe)$. Then $u_i. W_i = V_i$ if and only if $\mfd(u_i) = \mfw$, $\mft(u_i) = \mfv$ and $u_i \mff =^* \mfe$. The first two conditions are decidable as $\fC^0$ is finite, and the last condition is decidable again using $E$. This shows that it is decidable whether $([a_i,V_i])_i =^* ([b_i,W_i])_i$ holds in $\C$, and if that is the case, define $\bmE(([a_i,V_i])_i, ([b_i,W_i])_i)$ as $([u_i,W_i])_i$, which is computable because $E$ is computable.
\eproof
\pars

\bproof[Proof of Corollary~\ref{cor:WordProblem}]
First, we need an algorithm which starts with a word $\omega$ in $\bfS$, $\bfS^{-1}$ and $\C^*$ which represents an element of $\Q(*,*)$ and produces a word of the form $\alpha_1 \dotsm \alpha_m \beta_1^{-1} \dotsm \beta_n^{-1}$ for some $\alpha_i \in \bfS \cup \C^*$ and $\beta_j \in \bfS \cup \C^*$ such that $\omega = \alpha_1 \dotsm \alpha_m \beta_1^{-1} \dotsm \beta_n^{-1}$ in $\C$. This algorithm runs through the letters of $\omega$, starting from the right, and replaces a subword of the form $\alpha^{-1} \beta$ by a word of the form $\alpha' u (\beta')^{-1}$ for some $\alpha', \beta' \in \bfS$ and $u \in \C^*$ (some of these words could be empty). It suffices to explain the replacements in the case where $\alpha, \beta \in \fX \bfS$ and $\alpha \in \fX \bfS$, $\beta \in \C^*$. Extending this in a $\amalg$-compatible way, this will explain the replacements whenever $\alpha \notin \C^*$. In case $\alpha, \beta \in \fX \bfS$, define $\alpha', \beta', u$ such that $\lcm_{\Pi}(\alpha,\beta) = \alpha \alpha' u = \beta \beta'$. In case $\alpha \in \fX \bfS$, $\beta \in \C^*$, define $\beta'$ as the element of $\bmd(\beta) \bfS$ which is uniquely determined by the property that $\bmE(\alpha,\beta \beta')$ is defined. $\beta'$ is computable as $\bmd(\beta) \bfS$ is finite and $\bmE$ is computable. Furthermore, define $u \defeq \bmE(\alpha,\beta \beta')$, and let $\alpha'$ be the empty word. It is now straightforward to check that this algorithm always terminates and that it produces the desired output.

Now suppose that we are given a word $\omega$ in $\bfS$, $\bfS^{-1}$ and $\C^*$ which represents an element of $\Q(*,*)$, and that our algorithm produces a word of the form $\alpha_1 \dotsm \alpha_m \beta_1^{-1} \dotsm \beta_n^{-1}$. Deciding whether $\omega$ represents the trivial element of $\Q(*,*)$ amounts to deciding whether $\alpha_1 \dotsm \alpha_m = \beta_1 \dotsm \beta_n$ in $\C$. Now Corollary~\ref{cor:thmA} implies that $\bfS$ is a Garside family in $\C$. Hence, as explained in \cite[Chapter~III, Proposition~3.59 and Corollary~3.60]{Deh15}, the word problem in $\C$ is decidable if there exists a computable $\Box$-witness $\bmF$ for $\bfS^{\sharp}$, i.e., a function defined on pairs $(\alpha,\beta)$, where $\alpha, \beta \in \bfS^{\sharp}$ such that $\alpha \beta$ is a path, with the property that $\bmF(\alpha,\beta) = (\sigma,\tau) \in \bfS^{\sharp} \times \bfS^{\sharp}$ such that $\sigma \tau$ is a normal form of $\alpha \beta$.

To define $\bmF$ and to see that it is computable, take a path $\alpha \beta$ in $\bfS^{\sharp}$. For $\xi \in \bfS$, we have $\xi \preceq_{\Pi} \alpha \beta$ if and only if our algorithm above, applied to $\xi^{-1} \alpha \beta$, produces a word in $\bfS$ and $\C^*$ (but containing no letters in $\bfS^{-1}$). Similarly, given $\eta, \zeta \in \bfS$, we have $\eta \preceq_{\Pi} \zeta$ if and only if our algorithm above, applied to $\xi^{-1} \zeta$, produces a word in $\bfS$ and $\C^*$ (but containing no letters in $\bfS^{-1}$). Hence, applying our algorithm finitely many times, we can find the unique element $\xi \in \bfS$ with $\xi \preceq_{\Pi} \alpha \beta$ and which is $\preceq_{\Pi}$-maximal with that property. Now we know that there exist $\eta \in \bfS$ and $u \in \C^*$ such that a normal form of $\alpha \beta$ is given by $\xi (\eta u)$. Our algorithm, applied to $\xi^{-1} \alpha \beta$, yields a word $\chi_1 \chi_2 \dotsm$ in $\bfS \cup \C^*$. Applying our algorithm finitely many times, we can find $\eta \in \bfS$ such that our algorithm, applied to $\eta^{-1} \chi_1 \chi_2 \dotsm$, produces a word with letters only in $\C^*$. In this way, we compute $\eta$ and $u$. Now define $\bmF(\alpha,\beta) \defeq (\xi,\eta u)$. This shows that there exists a computable $\Box$-witness $\bmF$ for $\bfS^{\sharp}$, and thus, as explained above, the proof of Corollary~\ref{cor:WordProblem} is complete.
\eproof

\bremark
The first algorithm we construct in the proof of Corollary~\ref{cor:WordProblem} is closely related to the method of right-reversing (see for instance \cite[\S~4 in Chapter~II]{Deh15}).
\eremark

\section{Finiteness properties of topological full groups}
\label{s:Fn}

Throughout this section, let $\fC$ be a finitely aligned left cancellative small category with finite $\fC^0$ which is right cancellative up to $=^*$. Assume that $\fS$ is a Garside family in $\fC$ which is locally finite, $=^*$-transverse with $\fS \cap \fC^* = \emptyset$, and suppose that for all $L \geq 1$, $(\fS^{\leq L})^{\sharp}$ is closed under left divisors. Let $X \subseteq \Omega_{\infty}$ be a closed invariant subspace. Let $\fX$, $\bfX$, $\Gamma$, $\bGamma$, $* \in \bfX$, $\C$ and $\Q$ be as in \S~\ref{s:GarsCat-TFG}. We need the following condition:
\nopar

\bdefin
We say that condition ($\bmt < \bmd$) holds if for all $\gamma \in \Gamma(*)$ with $\bmt(\gamma) = X(\mfv;\mfe)$, there exists $u \in \C^*$ such that $\bmd(\gamma u) = (U_i)_{i \in I}$ and there exist $i_1, i_2 \in I$ such that $i_1 \neq i_2$ and $U_{i_1} = X(\mfv;\mfe) = U_{i_2}$.
\edefin
\pars

The goal now is to prove the following result.
\btheo
\label{thm:Fn}
Assume that condition (F) holds, that $\fX$ and $\Gamma$ satisfy condition (St) and that $\C$ admits lcms. If condition ($\bmt < \bmd$) is satisfied, then for all natural numbers $n$, $\Q(*,*)$ is of type ${\rm F}_n$ if $\fC^*(\mfv,\mfv)$ is of type ${\rm F}_n$ for all $\mfv \in \fC^0$. 
\pari

In particular, if $* = (X(\mfv;\mfe))_{\mfv \in \fV}$ for some $\fV \in \fC^0$ and $Y = \coprod_{\mfv \in \fC^0} Y_{\mfv}$, where $Y_{\mfv} = \emptyset$ if $\mfv \notin \fV$ and $Y_{\mfv} = X(\mfv;\mfe) \in \fX$ if $\mfv \in \fV$, then for all natural numbers $n$, $\bmF((I_l \ltimes X)_Y^Y)$ is of type ${\rm F}_n$ if $\fC^*(\mfv,\mfv)$ is of type ${\rm F}_n$ for all $\mfv \in \fC^0$.
\etheo
\nopar

Note that $\C$ admits lcms if, for example, in addition to our assumptions above, $\fC$ is right Noetherian, admits disjoint mcms, and condition (LCM) holds.
\pars

\subsection{The strategy}
\label{subsec:strategy}

We briefly summarize the strategy developed in \cite{Wit19} for establishing finiteness properties for fundamental groups of left reversible cancellative categories with Garside families closed under factors. We will be able to apply this strategy in our setting because Corollary~\ref{cor:ab=gamma} and Proposition~\ref{prop:GarsideTFG} imply that our category $\C$ is cancellative and left reversible, and $\crS = \Div_{\C}(\Delta)$ is a Garside family if (St) holds and $\C$ admits lcms.

First, assume that $\rho: \: \bfX(*) \to \Zz_{\geq 0}$ is a height function, i.e., $\rho(x) = \rho(y)$ if $\C(x,y) \cap \C^* \neq \emptyset$ and $\rho(x) < \rho(y)$ if $\C(x,y) \neq \emptyset$ and $\C(x,y) \cap \C^* = \emptyset$. Moreover, suppose that for all $R \in \Nz$, $\# \menge{x \in \bfX(*)}{\rho(x) \leq R} < \infty$. For each $x \in \bfX(*)$, we construct the poset $E(x)$ as follows: Its underlying set is given by classes $[\alpha]$, for $\alpha \in \crS(-,x) \setminus \C^*(-,x)$, where we define $\alpha \sim \alpha'$ if $\bm{u} \alpha = \alpha'$ for some $\bm{u} \in \C^*$. Moreover, we define $[\alpha] \leq [\beta]$ if $\varepsilon \alpha = \beta$ for some $\varepsilon \in \crS$. Here $\crS(-,x)$ and $\C^*(-,x)$ denote the elements of $\crS$ and $\C^*$, respectively, with domain equal to $x$. In this situation, Witzel established the following criterion:
\btheo[{\cite[Theorem~3.12]{Wit19}}]
\label{thm:Wit}
Assume that $\C^*(x,x)$ is of type ${\rm F}_n$ for all $x \in \bfX(*)$ and that there exists $N \in \Nz$ such that the order complex $\vert E(x) \vert$ is ($n-1$)-connected for all $x \in \bfX(*)$ with $\rho(x) \geq N$. Then $\Q(*,*)$ is of type ${\rm F}_n$.
\etheo

\subsection{Establishing finiteness properties}

We follow the strategy from \S~\ref{subsec:strategy}. As explained above, because of our assumptions and by Corollary~\ref{cor:ab=gamma} and Proposition~\ref{prop:GarsideTFG}, our category $\C$ and $\crS = \Div_{\C}(\Delta)$ have the desired properties. As a first step, we need to define $\rho$. For $\bmU = (U_i)_{i \in I} \in \bfX(*)$, define $\bmm_{\bmU}: \: \fX(*) \to \Nz$ by $\bmm_{\bmU}(X(\mfv; \mff)) \defeq \# \menge{i \in I}{U_i = X(\mfv; \mff)}$. Moreover, for $\gamma \in \Gamma$, set $\fX(\bmd(\gamma)) \defeq \menge{X(\mfv; \mff) \in \fX}{X(\mfv; \mff) \in \bmd(\gamma)}$. Define
$$
 \rho'(\bmU) \defeq \max \Big\lbrace \sum_{\gamma} \rho_{\gamma} : \: \rho_{\gamma} \in \Nz \ \forall \ \gamma \in \Gamma, \, \sum_{\gamma} \rho_{\gamma} \bm{1}_{\fX(\bmd(\gamma))} \leq \bmm_{\bmU} \Big\rbrace.
$$
Now set, for $x \in \bfX(*)$, $\rho(x) \defeq \max \menge{\rho'(\bmU)}{\bmU \in \bfX(*), \, \C^*(x,\bmU) \neq \emptyset}$.

\blemma
$\rho$ is a height function. For all $R \in \Nz$, we have $\# \menge{x \in \bfX(*)}{\rho(x) \leq R} < \infty$.
\elemma
\nopar

\bproof
By construction, we have $\rho(x) = \rho(y)$ if $\C^*(x,y) \neq \emptyset$. Given $x, y \in \bfX(*)$ with $\C(x,y) \neq \emptyset$ and $\C^*(x,y) = \emptyset$, we may assume because of (3$_{\Gamma}$) that $\rho(x) = \rho'(x)$. So there exists $\alpha \in \spkl{\bGamma}$ and $u \in \C^*$ with $x = \bmt(\alpha u)$ and $y = \bmd(\alpha u)$. Without loss of generality, we may assume $\alpha \in \bGamma$. By construction, $\rho(y) = \rho(\bmd(\alpha))$. Now condition ($\bmt < \bmd$) implies that $\rho'(x) < \rho'(\bmd(\alpha))$. Hence we deduce $\rho(x) = \rho'(x) < \rho'(\bmd(\alpha)) \leq \rho(\bmd(\alpha)) = \rho(y)$, as desired.
\pars

For the last claim, take $x \in \bfX(*)$. Then there exists $\alpha \in \spkl{\bGamma}$ and $u \in \C^*$ with $\bmt(\alpha u) = *$ and $\bmd(\alpha u) = x$. Condition ($\bmt < \bmd$) implies that, if $\ell_{\Gamma}(\alpha) > R$, then $\rho(x) > R$. Here $\ell_{\Gamma}(\alpha)$ is the minimal number of factors in $\bGamma$ needed to express $\alpha$ as a product. In other words, $\rho(\bmd(\alpha u)) \leq R$ implies $\ell_{\Gamma}(\alpha) \leq R$. Now suppose that $\alpha$ is a product of at most $R$ factors in $\bGamma$. For each factor, there are only finitely many possibilities because $\bmt(\alpha) = *$, $\# \fC^0 < \infty$ and for all $\mfv \in \fC^0$, there are only finitely many possibilities for subsets $\mfe, \, \mfs \subseteq \mfv \fS$ because $\fS$ is locally finite. It follows that
$$
 \# \menge{\bmd(\alpha u)}{\alpha \in \spkl{\bGamma}, \, \ell_{\Gamma}(\alpha) \leq R, \, u \in \C^*} < \infty,
$$
and thus $\# \menge{x \in \bfX(*)}{\rho(x) \leq R} < \infty$, as desired.
\eproof
\pars

Our next aim is to construct another height function $h$. To do that, let us fix an object $x$ in $\C$. Let $\cL$ denote the set of all proper left divisors in $\C$ of $\Delta(U)$ with domain contained in $x$, where $U$ runs through all elements of $\fX(*)$. Further set $\cL_+ \defeq \menge{\lambda \in \cL}{\lambda \text{ not an atom}}$, $\cL_a \defeq \menge{\lambda \in \cL}{\lambda \text{ is an atom}}$ and $\cL_0 \defeq \bigcup_{U \in \fX(*), \, U \in x} \C^*(-,U)$. Set $\cE(x) \defeq \crS(-,x) \setminus \C^*(-,x)$, and for $\alpha, \beta \in \cE(x)$, write $\alpha \leq \beta$ if $\varepsilon \alpha = \beta$ for some $\varepsilon \in \crS$. Every $\mu \in \cE(x)$ is of the form $\mu = \mu_+ \amalg \mu_a \amalg \mu_0$, where $\mu_+ = \coprod \mu_+^k$, $\mu_a = \coprod \mu_a^j$ and $\mu_0 = \coprod \mu_0^i$, for some $\mu_+^k \in \cL_+$, $\mu_a^j \in \cL_a$ and $\mu_0^i \in \cL_0$. For such a $\mu$, define $n_+([\mu]) : \: [\cL_+] \to \Nz, \, n_+([\mu])([\lambda]) \defeq \# \menge{k}{[\mu_+^k] = [\lambda]}$. We set
$$
 \vec{h}: \: E(x) \to \Nz^{[\cL_+]} \times \Nz, \, [\mu] \ma ((n_+([\mu])(\lambda))_{[\lambda]}, \rho(\mu_0)).
$$
It is easy to see that this is well-defined. Now introduce a total order $\dot{<}$ on $[\cL_+]$ so that $[\lambda_1] < [\lambda_2]$ implies $[\lambda_1] \dot{<} [\lambda_2]$. This order induces the lexicographical order $\dot{<}$ on $\Nz^{[\cL_+]} \times \Nz$, where the last copy of $\Nz$ is considered last. Finally, define $h: \: E(x) \to \Nz$ such that $\vec{h}([\mu_2]) \dot{<} \vec{h}([\mu_1])$ implies $h([\mu_2]) < h([\mu_1])$.

By construction, if $[\mu] < [\nu]$, then either $n_+([\mu]) <  n_+([\nu])$ or $n_+([\mu]) =  n_+([\nu])$ and $\rho(\mu_0) > \rho(\nu_0)$. In the first case, we have $h([\mu]) < h([\nu])$, and in the second case, we have $h([\mu]) > h([\nu])$.

Now let $E(x)^0$ be the full subcomplex of $E(x)$ supported on vertices $[\mu] \in E(x)$ with $n_+([\mu]) = 0$. Let us show that $E(x)^0$ is highly connected, and that, for every $[\mu] \in E(x)$, the descending link $\menge{[\nu] \in E(x)}{h([\mu]) > h([\nu])}$, viewed as a subcomplex, is highly connected (or rather the corresponding order complex). To achieve this, we split the descending link into up-link 
$$\menge{[\nu] \in E(x)}{h([\mu]) > h([\nu]), \, [\mu] < [\nu]}$$
and down-link
$$\menge{[\nu] \in E(x)}{h([\mu]) > h([\nu]), \, [\mu] > [\nu]}.$$
The descending link is the join of the up-link and the down-link (see \cite[Observation~3.6]{FMWZ}), so that it suffices to show that the up-link or the down-link is highly connected, as long as $\rho(x)$ is sufficiently big. To achieve this, it suffices to treat the case that $\rho'(x)$ is sufficiently big, because if there exists $u \in \C^*(x,x')$, then we obtain a poset isomorphism $E(x) \cong E(x'), \, [\mu] \ma [\mu u]$ with inverse $[\mu' u^{-1}] \mapsfrom [\mu']$. Thus our goal is to show the following: 

\bprop
\label{prop:link:n-conn}
For each $n \in \Nz$, there exists $R \in \Nz$ such that for all $x \in \bfX(*)$ with $\rho(x) \geq R$ and for all $[\mu] \in E(x)$, the up-link or down-link of $[\mu]$ is ($n-1$)-connected.
\eprop

We proceed in several steps. 
\blemma
\label{lem:link:n-conn_1}
For all $n \in \Nz$, $x \in \bfX(*)$ and $[\mu] \in E(x)$ with $\# \gekl{\mu_+^k} + \# \big\lbrace \mu_a^j \big\rbrace \geq n+1$, the down-link of $[\mu]$ is ($n-1$)-connected.
\elemma
\nopar

\bproof
If $[\nu]$ is in the down-link $\menge{[\nu] \in E(x)}{h([\mu]) > h([\nu]), \, [\mu] > [\nu]}$, then we must have $n_+([\mu]) > n_+([\nu])$. Write $\mu = (\coprod_k \mu_+^k) \amalg \mu_a \amalg \mu_0$. Define $\omega(k,\omega^k) = \coprod_l \omega(k,\omega^k)^l$ by $\omega(k,\omega^k)^l = \mu_+^l$ if $l \neq k$ and $\omega(k,\omega^k)^k = \omega^k$, where $[\omega^k]$ is maximal with $[\omega^k] < [\mu_+^k]$, i.e., we have $\mu_+^k = a \omega^k$ where $a$ is an atom. Then the down-link is given by $\menge{[\nu]}{[\nu] \leq [\omega(k,\omega^k) \amalg \mu_a \amalg \mu_0] \text{ for some } k}$. So if we define the subcomplex $\cS(k,\omega^k) \defeq \menge{[\nu] \in E(x)}{[\nu] \leq [\omega(k,\omega^k) \amalg \mu_a \amalg \mu_0]}$, then the down-link is given by $\bigcup_k \cS(k,\omega^k)$.
\pars

To see the last claim, observe that if $\mu = a \nu$ and $[\nu]$ is maximal, then $a$ must be an atom. So if $\mu = (\coprod_k \mu_+^k) \amalg \mu_a \amalg \mu_0$, then $a$ must respect this $\coprod$-decomposition. $[\nu_+] < [\mu_+]$ implies that $a = a_+^k$ for some $k$.

(a) Finite intersections of $\cS(k,\omega^k)$ are empty or contractible: Given a finite collection $\gekl{\cS(k_p,\omega^{k_p})}_p$, let $\chi^l \defeq \gcd( \gekl{\omega(k_p,\omega^{k_p})^l}_p )$. Here we need the existence of gcds for right divisors of left divisors of $\Delta$. This follows from existence of lcms for left divisors of $\Delta$:

If $s = a_k t_k$, and $a = \lcm \gekl{a_k}$ in $\crS$, then $a = a_k b_k$. So $s = a t$ for some $t \in \crS$. Then (i) $t \leq t_k$ for all $k$ and (ii) if $\eta \leq t_k$ for all $k$, then $\eta \leq t$. Indeed, to see (i), we have $a_k b_k t = a_k t_k$, which implies $b_k t = t_k$, and $\crS$ is closed under divisors. To see (ii), assume $\chi_k \eta = t_k$. Then $s = a_k \chi_k \eta$, so that $a_k \chi_k \eta = a_l \chi_l \eta$ for all $k, l$. This implies $a_k \chi_k = a_l \chi_l$ and hence $a_k \chi_k = a \zeta = a_k b_k \zeta$ and hence $\chi_k = b_k \zeta$. So $s = a_k b_k \zeta \eta = a t$, so that $\zeta \eta = t$. All this takes place in $\crS$ (compare also \cite[Chapter~I, Lemma~2.37]{Deh15}).

Then $\bigcap_p \cS(k_p,\omega^{k_p}) = \menge{[\nu] \in E(x)}{[\nu] \leq [(\coprod_l \chi^l) \amalg \mu_a \amalg \mu_0]}$, which is contractible.

(b) For all $m < n$, all $m$-fold intersections of $\cS(k,\omega^k)$ (i.e., $\bigcap_{p \in P} \cS(k_p,\omega^{k_p})$, where $\# P \leq m$) are non-empty: This is clear if $\# \big\lbrace \mu_a^j \big\rbrace \geq 1$. So we may assume $\# \big\lbrace \mu_a^j \big\rbrace = 0$. Then $\# \gekl{\mu_+^k} + \# \big\lbrace \mu_a^j \big\rbrace \geq n+1$ implies $\# \gekl{\mu_+^k} \geq n+1$. Since $m < n$, we conclude that there exists $l$ such that $\omega(k_p,\omega^{k_p})^l = \mu^l$ for all $p \in P$, because $\# P < n$. Define $\zeta^k \defeq \mu^l$ if $k = l$ and $\zeta^k \defeq \bmd(\mu^k)$ if $k \neq l$. Then $[(\coprod_k \zeta^k) \amalg \mu_a \amalg \mu_0] \in \bigcap_{p \in P} \cS(k_p,\omega^{k_p})$.

Now \cite[Theorem~6]{Bjo03} implies that, if $\cD \cL$ denotes the down-link, $\pi_{\bullet}(\cD \cL) \cong \pi_{\bullet}(\cN(\cS(k,\omega^k)))$ for all $\bullet < n$, where $\cN(\cS(k,\omega^k))$ stands for the nerve simplex of $\gekl{\cS(k,\omega^k)}$ in the sense of \cite[\S~4]{Bjo03}. (a) and (b) together imply that $(\cN(\cS(k,\omega^k)))^m \cong (\cF \cS)^m$, where $\cF \cS$ stands for full simplex. Now apply \cite[Corollary~4.12, p.~360]{Hat} to obtain that $(\cN(\cS(k,\omega^k)))^m \to \cN(\cS(k,\omega^k))$ is a $\pi_{\bullet}$-isomorphism for all $\bullet < m$, and also that $(\cF \cS)^m \to \cF \cS$ is a $\pi_{\bullet}$-isomorphism for all $\bullet < m$. It follows that $\pi_{\bullet}(\cD \cL) \cong \pi_{\bullet}(\cF \cS) \cong (0)$.
\eproof
\pars

Hence, to prove Proposition~\ref{prop:link:n-conn}, it suffices to show the following.
\blemma
\label{lem:link:n-conn_2}
For each $n \in \Nz$, there exists $R \in \Nz$ such that for all $x \in \bfX(*)$ with $\rho(x) \geq R$ and for all $[\mu] \in E(x)$ with $\# \gekl{\mu_+^k} + \# \big\lbrace \mu_a^j \big\rbrace < n+1$, the up-link $[\mu]$ is ($n-1$)-connected.
\elemma

Fix $x \in \bfX(*)$ and $[\mu] \in E(x)$. Consider the up-link $\menge{[\nu] \in E(x)}{h([\mu]) > h([\nu]), \, [\mu] < [\nu]}$. If $[\nu]$ is in the up-link, then we must have $n_+([\mu]) = n_+([\nu])$ and $n_0([\mu]) > n_0([\nu])$. So $[\nu_+] = [\mu_+]$. Let $y = \bm{d}(\mu_0)$. Then $\nu = \nu_+ \amalg \nu_a \amalg \nu_0 = \mu_+ \amalg \mu_a \amalg \nu'$ with $\bm{d}(\nu') = y$. Let $\cM_{\mu}$ be the complex with set of vertices $\cV \defeq \menge{[\alpha]}{\alpha \text{ atom in } \crS \text{ with } \bm{d}(\alpha) \subseteq y}$, and $\gekl{[\alpha_i]}$ is a simplex in $\cM_{\mu}$ if $\bm{d}(\alpha_i)$ are pairwise disjoint and for $z \subseteq y$ such that $(\coprod_i \bmd(\alpha_i)) \amalg z = y$, there exists $[\nu]$ in the up-link of $[\mu]$ such that $\nu' = (\coprod_i \alpha_i) \amalg z$, i.e., $\bmt(\mu_+ \amalg \mu_a \amalg (\coprod_i \alpha_i) \amalg z)$ lies in $\bfX(*)$. By construction, the up-link of $[\mu]$ is isomorphic to the face poset of $\cM_{\mu}$ via $[\nu] \ma \menge{[\alpha]}{\alpha \in \cV, \, \alpha \leq \nu'}$. So connectivity of the up-link reduces to connectivity of $\cM_{\mu}$ (see \cite[\S~12.4, p.~1860--1861]{Bjo95}). Thus we set out to prove the following:

\blemma
\label{lem:link:n-conn_3}
For all $L, m \in \Nz$, there exists $R \in \Nz$ such that for all $x \in \bfX(*)$ with $\rho'(x) \geq R$ and for all $[\mu] \in E(x)$ with $\# \gekl{\mu_+^k} + \# \big\lbrace \mu_a^j \big\rbrace < n+1$, and for all $m$-simplices $\sigma_1, \dotsc, \sigma_L$ in $\cM_{\mu}$, there exists a vertex $[\alpha]$ of $\cM_{\mu}$ such that $\sigma_l \cup \gekl{[\alpha]}$ is an ($m+1$)-simplex of $\cM_{\mu}$ for all $1 \leq l \leq L$.
\elemma
\nopar

\bproof
Assume that the statement is not true. Then there exist $L, m \in \Nz$ and a sequence $\mu_p = \mu_{p,+} \amalg \mu_{p,a} \amalg \mu_{p,0} \in \cE(x_p)$ with $\# \gekl{\mu_{p,+}^k} + \# \big\lbrace \mu_{p,a}^j \big\rbrace < n+1$ and $\rho'(x_p) \nearrow \infty$ as $p \to \infty$, $m$-simplices $\sigma_{p,l}$ in $\cM_{\mu_p}$ such that there exists no vertex $[\alpha]$ in $\cM_{\mu_p}$ such that $\sigma_{p,l} \cup \gekl{[\alpha]}$ is an ($m+1$)-simplex of $\cM_{\mu_p}$. Our goal is to derive a contradiction. 
\pars

Assume that under the identification between the up-link of $[\mu_p]$ and the face poset of $\cM_{\mu_p}$ from above, $\sigma_{p,l}$ corresponds to $[\nu_{p,l}]$, where $\nu_{p,l} \in \cE(x_p)$ is of the form $\nu_{p,l} = \mu_{p,+} \amalg \mu_{p,a} \amalg \theta_{p,l} \amalg z_{p,l}$, where $\theta_{p,l}$ is a disjoint union of $m$ atoms, and $z_{p,l} \in \bfX(*)$. For $\alpha \in \C$, define $\bmm_{\alpha} \defeq \bmm_{\bmd(\alpha)} - \bmm_{\bmt(\alpha)}$. By passing to a subsequence if necessary, we may assume that $\bmm_{\nu_{p,l}}$ is independent of $p$. Since $\bmt(\nu_{p,l}) \in \bfX(*)$, there exist $\zeta_{p,l}, \xi_{p,l} \in \C$ with $\bmt(\zeta_{p,l} \xi_{p,l} \nu_{p,l}) = *$ and such that, by passing to a subsequence if necessary, we may assume that $\bmm_{\zeta_{p,l}}$ is independent of $p$, and for all $U \in \fX(*)$, either $\bmm_{\xi_{p,l}}(U) = 0$ or $\bmm_{\xi_{p,l}}(U)$ is strictly increasing as $p \to \infty$. Moreover, we have $\bmm_{\zeta_{p,l} \xi_{p,l} \nu_{p,l}} = \bmm_{x_i} - \bmm_{*}$, hence independent of $l$. We conclude that 
$$
 \bmm_{\zeta_{p,l}} + \bmm_{\xi_{p,l}} + \bmm_{\nu_{p,l}} = \bmm_{\zeta_{p,l'}} + \bmm_{\xi_{p,l'}} + \bmm_{\nu_{p,l'}},
$$
so that
$$
 \bmm_{\xi_{p,l}} - \bmm_{\xi_{p,l'}} = \bmm_{\zeta_{p,l'}} + \bmm_{\nu_{p,l'}} - \bmm_{\zeta_{p,l}} - \bmm_{\nu_{p,l}}
$$
is independent of $p$.

For sufficiently big $p$ (and after passing to a subsequence if necessary), there exists $\theta \in \C$, independent of $l$, such that $\bmd(\theta) \subseteq z_{p+2,l}$ and $\bmm_{\theta} = \bmm_{\xi_{p+1,l}} - \bmm_{\xi_{p,l}}$ for all $l$. Define
$$
 \ti{\nu}_{p+2,l} \defeq \mu_{p+2,+} \amalg \mu_{p+2,a} \amalg \theta_{p+2,l} \amalg \theta \amalg \ti{z}_{p+2,l},
$$
where $\bmd(\theta) \amalg \ti{z}_{p+2,l} = z_{p+2,l}$. Then, again for sufficiently big $p$ (and after passing to a subsequence if necessary), there exists $\eta_{p+2,l} \in \C$ such that $\bmm_{\eta_{p+2,l}} + \bmm_{\theta} = \bmm_{\xi_{p+2,l}} - \bmm_{\xi_{p,l}}$ and $\bmd(\eta_{p+2,l}) = \bmt(\ti{\nu}_{p+2,l})$. We claim that $\bmm_{\bmt(\eta_{p+2,l} \ti{\nu}_{p+2,l})} = \bmm_{\bmd(\xi_{p,l})}$. Indeed, we have
\begin{align*}
 \bmm_{x_{p+2}} - \bmm_{\bmt(\eta_{p+2,l} \ti{\nu}_{p+2,l})} &= \bmm_{\eta_{p+2,l} \ti{\nu}_{p+2,l}} = \bmm_{\eta_{p+2,l}} + \bmm_{\ti{\nu}_{p+2,l}} + \bmm_{\theta} = \bmm_{\xi_{p+2,l}} - \bmm_{\xi_{p,l}} + \bmm_{\nu_{p+2,l}}\\
 &= \bmm_{\bmd(\xi_{p+2,l)}} - \bmm_{\bmd(\xi_{p,l})} + \bmm_{x_{p+2}} - \bmm_* - \bmm_{\zeta_{p+2,l} \xi_{p+2,l}}\\
 &= \bmm_{\bmd(\xi_{p+2,l)}} - \bmm_{\bmd(\xi_{p,l})} + \bmm_{x_{p+2}} - \bmm_* + \bmm_* - \bmm_{\bmd(\xi_{p+2,l})}\\
 &= \bmm_{x_{p+2}} - \bmm_{\bmd(\xi_{p,l})}.
\end{align*}
We deduce that $\bmt(\eta_{p+2,l}) \in \bfX(*)$ and thus $\bmt(\ti{\nu}_{p+2,l}) \in \bfX(*)$, hence $\ti{\nu}_{p+2,l} \in \C$. 

Now take any atom $\alpha$ with $\alpha < \theta$. Then $[\alpha]$ is a vertex of $\cM_{\mu_{p+2}}$ and $\theta_{p+2,l} \cup \gekl{[\alpha]}$ is an ($m+1$)-simplex of $\cM_{\mu_{p+2}}$. This is the desired contradiction.
\eproof
\pars

\bproof[Proof of Lemma~\ref{lem:link:n-conn_2}]
Lemma~\ref{lem:link:n-conn_3} implies Lemma~\ref{lem:link:n-conn_2} following the argument for \cite[Lemma~6.18]{Mat15}, which in turn is based on the proof of \cite[Lemma~4.20]{Bro}.
\eproof

\bproof[Proof of Proposition~\ref{prop:link:n-conn}]
Proposition~\ref{prop:link:n-conn} follows from Lemmas~\ref{lem:link:n-conn_1} and \ref{lem:link:n-conn_2}.
\eproof

\bcor
\label{cor:link:n-conn}
For each $n \in \Nz$, there exists $R \in \Nz$ such that for all $x \in \bfX(*)$ with $\rho(x) \geq R$, $\vert E(x) \vert$ is ($n-1$)-connected.
\ecor
\nopar

\bproof
$E(x)^0$ is the up-link for $[\mu]$, where $\mu = x$, i.e., $\mu_+ = \emptyset$ and $\mu_a = \emptyset$. So Lemma~\ref{lem:link:n-conn_2} implies that $E(x)^0$ is ($n-1$)-connected. Now apply Morse theory, in the form of \cite[Lemma~3.11]{Wit19} ($X$ in \cite{Wit19} is our $E(x)$ and $\rho$ in \cite{Wit19} is our $h$), to deduce that the canonical inclusion $E(x)^0 \into E(x)$ induces a $\pi_{\bullet}$-isomorphism for all $\bullet < n$. Hence $\vert E(x) \vert$ is ($n-1$)-connected.
\eproof
\pars

The final ingredient is the following observation:
\nopar

\blemma
\label{lem:Stab}
Assume that $\fC$ is right cancellative up to $=^*$ and that condition (F) holds. If $\fC^*(\mfv,\mfv)$ is of type ${\rm F}_n$ for all $\mfv \in \fC^0$, then $\C^*(x,x)$ is of type ${\rm F}_n$.
\elemma
\bproof
First, consider $x = (U_i)_{i \in I}$. We have a short exact sequence 
$$
 1 \to \prod_{i \in I} \C^*(U_i,U_i) \to \C^*(x,x) \to {\rm Perm}(I) \to 1.
$$
As ${\rm Perm}(I)$ is finite, \cite[Corollary~7.2.4]{Geo} implies that $\C^*(x,x)$ is of type ${\rm F}_n$ if and only if $\prod_{i \in I} \C^*(U_i,U_i)$ is of type ${\rm F}_n$. Furthermore, if $\C^*(U_i,U_i)$ is of type ${\rm F}_n$ for all $i \in I$, then $\prod_{i \in I} \C^*(U_i,U_i)$ is of type ${\rm F}_n$ (see for instance \cite[Exercise~1 in \S~7.2]{Geo}). Thus it suffices to show that $\C^*(U,U)$ is of type ${\rm F}_n$, where $U = X(\mfv;\mfe) \in \fX(*)$. Now suppose that $[a,U]$ lies in $\C^*(U,U)$. Then there exists $a' \in \fC$ such that $[a',a.U] [a,U] = [a'a,U] = U$. Lemma~\ref{lem:rNrc-->rN}~(i) implies that $a'a \in \fC^*$ and thus $a \in \fC^*$, hence $a \in \fC^*(\mfv,\mfv)$. This, in combination with condition (F), implies that we obtain an embedding $\C^*(U,U) \into \fC^*(\mfv,\mfv), \, [a,U] \ma a$. In addition, consider the short exact sequence 
$$
 1 \to {\rm Ker} \to \fC^*(\mfv,\mfv) \to {\rm Perm}(\menge{X(\mfv;\mfe')}{\mfe' \in \mfv \fS}) \to 1,
$$
where the second map sends $a \in \fC^*(\mfv,\mfv)$ to the permutation $V \ma a.V$ and ${\rm Ker}$ denotes the kernel of this map. Since ${\rm Perm}(\menge{X(\mfv;\mfe')}{\mfe' \in \mfv \fS})$ is finite, ${\rm Ker}$ is a finite index subgroup of $\fC^*(\mfv,\mfv)$, and because ${\rm Ker} \subseteq \C^*(U,U) \subseteq \fC^*(\mfv,\mfv)$, this shows that $\C^*(U,U)$ is a finite index subgroup of $\fC^*(\mfv,\mfv)$ as well. So if $\fC^*(\mfv,\mfv)$ is of type ${\rm F}_n$, then $\C^*(U,U)$ is of type ${\rm F}_n$ because finiteness properties are inherited by finite index subgroups.
\eproof
\pars

\bproof[Proof of Theorem~\ref{thm:Fn}]
The first part of the theorem follows from Theorem~\ref{thm:Wit} using Corollary~\ref{prop:link:n-conn} and Lemma~\ref{lem:Stab}. 

The second part of the theorem follows from the first, using Lemmas~\ref{lem:QbfC=bmF}, \ref{lem:QbfCYY=QbfCYvYv} and Proposition~\ref{prop:GarsideTFG}.
\eproof

\section{Garside families arising from degree maps}
\label{s:Gars-deg}

We set out to describe a class of small categories where the general criteria from \S~\ref{s:Fn} apply and allow us to establish finiteness properties for topological full groups of the groupoid models attached to our small categories.

Let $P$ be a left cancellative monoid with identity element $1$ and $P^* = \gekl{1}$. Suppose that $S_P \subseteq P$ is a finite Garside family in $P$ with $1 \notin S_P$. Let $\fC$ be a left cancellative small category equipped with a $P$-valued degree map, i.e., a functor $\bbd: \: \fC \to P$ with $\bbd^{-1}(1) = \fC^*$ such that the following unique factorization property holds:
\begin{enumerate}
\item[(UFP*)] For all $c \in \fC$ with $\bbd(c) = pq$, there exist $a, b \in \fC$ with $c = ab$ and $\bbd(a) = p$, $\bbd(b) = q$. If we have $c = a' b'$ for some $a', b' \in \fC$ with $\bbd(a') = p$, $\bbd(b') = q$, then there exists $u \in \fC^*$ such that $a' = au$ and $b' = u^{-1}b$.
\end{enumerate}
\nopar

This type of factorization property has been considered in a special context in \cite{LV}. Note that (UFP*) implies cancellation up to $=^*$.
\pars

Let $\fS$ be a $=^*$-transversal for $\bbd^{-1}(S_P)$, i.e., $\fS$ is a subset of $\fC$ such that the canonical projection restricts to a bijection $\fS \isom \bbd^{-1}(S_P) / {}_{=^*}$. As before, $\fS^{\sharp} = \fC^* \cup \fS \fC^*$.

\blemma
$\fS^{\sharp}$ is closed under right divisors. We have $(\fS^{\leq L})^{\sharp} = \bbd^{-1}(S_P^{\leq L})$. For all $L \geq 1$, if $(S_P^{\leq L})^{\sharp}$ is closed under left divisors, then so is $(\fS^{\leq L})^{\sharp}$.
\elemma
\nopar

\bproof
Given $s \in \fS$ and $a, b \in \fC$ with $s = ab$, it follows that $\bbd(s) = \bbd(a) \bbd(b)$ and hence $\bbd(b) \in S_P$. It follows that $b \in \fS$.
\pars

We have $\fS^{\leq L} \subseteq \bbd^{-1}(S_P^{\leq L})$ because $\bbd$ is a functor. To see \an{$\supseteq$}, assume that $\bbd(a) = s_1 \dotsm s_l$ for some $s_1, \dotsc, s_l \in S_P$. Then (UFP*) allows us to find $a_1, \dotsc, a_l \in \fC$ with $\bbd(a_k) = s_k$ for all $1 \leq k \leq l$ such that $a = a_1 \dotsm a_l$. It follows that $a_k \in \bbd^{-1}(S_P) = \fS^{\sharp}$ for all $1 \leq k \leq l$ and thus $a \in (\fS^{\leq L})^{\sharp}$.

If $s \in \fS^{\leq L}$ and $s = ab$ for some $a, b \in \fC$, then $\bbd(s) = \bbd(a) \bbd(b)$. Since $S_P^{\leq L}$ is closed under left divisors, we deduce that $\bbd(a) \in S_P^{\leq L}$. Hence it follows that $a \in (\fS^{\leq L})^{\sharp}$ by what we just proved.
\eproof
\pars

\blemma
If $P$ is left Noetherian, then so is $\fC$. Similarly, if $P$ is right Noetherian, then so is $\fC$.
\elemma
\nopar

\bproof
A chain $\dotso \prec a_2 \prec a_1$ in $\fC$ leads to a chain $\dotso \prec \bbd(a_2) \prec \bbd(a_1)$ in $P$. This shows our claim for \an{left Noetherian}. The argument for \an{right Noetherian} is analogous.
\eproof
\pars

The degree functor $\bbd$ induces a map $\dot{\bbd}: \: \fC / { }_{=^*} \to P$ because $\bbd(\fC^*) = \gekl{1}$.
\blemma
If $\mfv \dot{\bbd}^{-1}(s) < \infty$ for all $\mfv \in \fC^0$ and $s \in S_P$, then $\mfv \fS < \infty$.
\elemma
\nopar

\bproof
This follows because $\mfv \fS \cong \mfv \dot{\bbd}^{-1}(S_P)$, and $S_P$ is finite by assumption.
\eproof
\pars

\blemma
\label{lem:mcmCmcmP}
Given $a, b \in \fC$, we have $\mcm(a,b) = \bbd^{-1}(\mcm(\bbd(a),\bbd(b))) \cap (a \fC \cap b \fC)$ if $\mcm(\bbd(a),\bbd(b))$ exists.
\pari

Suppose $M \subseteq P$ is a subset closed under mcms. Then $\bbd^{-1}(M)$ is also closed under mcms. 
\elemma
\nopar

\bproof
Assume that $a \ti{a} = b \ti{b}$. Then $\bbd(a) \bbd(\ti{a}) = \bbd(b) \bbd(\ti{b})$. Hence, if $\mcm(\bbd(a),\bbd(b))$ exists, there exists $m \in \mcm(\bbd(a),\bbd(b))$ and $m' \in P$ with $\bbd(a) \bbd(\ti{a}) = m m' = \bbd(b) \bbd(\ti{b})$. Write $m = \bbd(a) p = \bbd(b) q$. Then $\bbd(\ti{a}) = p m'$ and $\bbd(\ti{b}) = q m'$. Now (UFP*) implies that there exist $a', a'', b', b'' \in \fC$ with $\ti{a} = a' a''$, $\ti{b} = b' b''$ such that $\bbd(a') = p$, $\bbd(b') = q$, $\bbd(a'') = m'$, $\bbd(b'') = m'$. Consider the identity $(aa') a'' = (bb') b''$. As $\bbd(aa') = \bbd(bb')$ and $\bbd(a'') = \bbd(b'')$, (UFP*) implies that there exists $u \in \fC^*$ with $bb' = aa' u$ and $b'' = u^{-1} a''$. Hence it follows that $aa' \in a \fC \cap b \fC$ and $a \ti{a} \in aa' \fC$. If $aa' = cz$ for some $c, z \in \fC$ with $c \in a \fC \cap b \fC$, then $m = \bbd(aa') = \bbd(c) \bbd(z)$ and $\bbd(c) \in \bbd(a) P \cap \bbd(b) P$ implies that $\bbd(z) = 1$ by maximality of $m$. It follows that $z \in \fC^*$, so that $aa' =^* c$. 
\eproof
\pars

The following are immediate consequences.
\nopar

\bcor
If $S_P$ is closed under mcms, then $\fS$ is closed under mcms.
\pari
If $P$ is right Noetherian and admits mcms, then $\fS$ is a Garside family.

If $P$ is finitely aligned and $\mfv \dot{\bbd}^{-1}(p) < \infty$ for all $p \in P$, then $\fC$ is finitely aligned.
\ecor
\pars

\blemma
\label{lem:deg_lcm-->mcm}
Suppose that $P$ admits conditional lcms. Then $\fC$ has disjoint mcms.
\elemma
\nopar

\bproof
Given $\mfv \in \fC^0$ and $a, b \in \mfv \fC$, take $C \subseteq \mfv \fC$ such that the canonical projection $\fC \to \fC / { }_{=^*}$ induces a bijection $C \isom (\bbd^{-1}(\lcm(a,b)) \cap (a \fC \cap b \fC)) / { }_{=^*}$. Then we claim that
$$
 a \fC \cap b \fC = \coprod_{c \in C} c \fC.
$$
To see that the sets $c \fC$ are pairwise disjoint, take $c_1, c_2 \in C$. Then $\bbd(c_1) = \bbd(c_2)$ and $c_1, c_2 \in a \fC \cap b \fC$. If $c_1 z_1 = c_2 z_2$, then $\bbd(c_1) \bbd(z_1) = \bbd(c_2) \bbd(z_2)$. $\bbd(c_1) = \bbd(c_2)$ implies $\bbd(z_1) = \bbd(z_2)$. Hence (UFP*) implies that $c_1 =^* c_2$, and by construction, we deduce $c_1 = c_2$. 
\pars

Lemma~\ref{lem:mcmCmcmP} implies $C = \mcm(a,b)$, hence $a \fC \cap b \fC = \coprod_{c \in C} c \fC$.
\eproof
\pars

Now assume that $P$ is right Noetherian and admits conditional lcms. Further assume that $(S_P^{\leq L})^{\sharp}$ is closed under left divisors for all $L \geq 1$, and that $\mfv \dot{\bbd}^{-1}(p) < \infty$ for all $\mfv \in \fC^0$ and $p \in P$.

Let $\mfv \in \fC^0$. A subset $\mfe \subseteq \mfv \fS$ is called saturated if $\mfe = \mfv \bbd^{-1}(\bbd(\mfe)) \cap \fS$. Set
\begin{eqnarray*}
 \fX &\defeq& \menge{X(\mfv; \mfe)}{\mfv \in \fC^0, \, \mfe \subseteq \mfv \fS \text{ saturated}},\\
 \Gamma &\defeq& \menge{\gamma(\mfe, \mfs)}{\mfv \in \fC^0; \; \mfe, \, \mfs \subseteq \mfv \fS \text{ saturated and } \mcm \text{-closed}},
\end{eqnarray*}
and let $\bfX$, $\bGamma$, $* \in \bfX$, $\C$ and $\Q$ be as in \S~\ref{s:GarsCat-TFG}. Our goal now is to check the conditions needed to apply Theorem~\ref{thm:Fn}.

\blemma
Condition (1$_{\fX}$) is satisfied.
\elemma
\nopar

\bproof
Take $u \in \fC^*$ and $\mfe = \mfv \bbd^{-1}(\bbd(\mfe)) \cap \fS$. Suppose that $\mff \subseteq \fS$ satisfies $u \mfe =^* \mff$. We claim that $\mff = \mft(u) \bbd^{-1}(\bbd(\mfe)) \cap \fS$. Indeed, given $f \in \mff$, we have $f v = u e$ for some $e \in \mfe$ and $v \in \fC^*$. It follows that $\bbd(f) = \bbd(fv) = \bbd(ue) = \bbd(e) \in \bbd(\mfe)$. Conversely, given $f' \in \mft(u) \bbd^{-1}(\bbd(\mfe)) \cap \fS$, then $u^{-1} f' \in \mfv \fC$ and $\bbd(u^{-1} f') = \bbd(f') \in \bbd(\mfe)$, so that $u^{-1} f' \in \mfe \fC^*$ and thus $f' \in u \mfe \fC^*$.
\eproof
\pars

\blemma
Condition (2$_{\fX}$) is satisfied.
\elemma
\nopar

\bproof
Given a finite subset $\mfe \subseteq \mfv \fS$, we want to construct $\bma \in \bfC$ with $\bmt(\bma) = X(\mfv; \mfe)$ and $\bmd(\bma) \in \bfX$. Let $\mfe_P$ be the lcm-closure of $\bbd(\mfe)$. Enumerate $\mfe_P = \gekl{e_1, e_2, \dotsc}$ so that $e_i \prec e_j$ implies $i < j$. Starting with the biggest index $i_{\max}$, we inductively construct $\bma_i \in \bfC$ such that $\bmt(\bma_i) = X(\mfv; \mfe_{i+1})$ and $\bmd(\bma_i) = X(\mfv; \mfe_i) \amalg x_i$ for some $x_i \in \bfX$, where $\mfe_{i_{\max}+1} \defeq \mfe$ and $\mfe_i \defeq \mfe \cup \bigcup_{j \geq i} \mfv \bbd^{-1}(e_j) \cap \fS$.
\pars

To construct $\bma_{i_{\max}}$, set $J \defeq (\mfv \bbd^{-1}(e_{i_{\max}}) \cap \fS \cup \gekl{\mfv}) \setminus \menge{j \in \fS}{e \preceq j \text{ for some } e \in \mfe}$. Now set $a_{i_{\max}, j} \defeq j$ if $j \neq \mfv$ and $a_{i_{\max}, \mfv} \defeq X(\mfv; e_{i_{\max}})$. Then $\bma_{i_{\max}} \defeq (a_{i_{\max}, j})_{j \in J}$ lies in $\bfC(X(\mfv; \mfe), (U_j)_{j \in J})$, where $U_j = X$ if $j \neq \mfv$ and $U_{\mfv} = X(\mfv; \mfe_{i_{\max}})$.

Now suppose that $\bma_{i+1}$ has been constructed. To construct $\bma_i$, set 
$$
 J \defeq (\mfv \bbd^{-1}(e_i) \cap \fS \cup \gekl{\mfv}) \setminus \menge{j \in \fS}{e \preceq j \text{ for some } e \in \mfe_{i+1}}.
$$
Define $a_{i, j} \defeq j X(\mfd(j); \mff_j)$ if $j \neq \mfv$, where $\mff_j$ is the set of minimal elements of
$$
 \bigcup \menge{\mfd(j) \bbd^{-1}(q) \cap \fS}{q \neq 1; \; \lcm(\bbd(e),e_i) = e_i q \text{ for some } e \in \mfe_{i+1}}.
$$
and $a_{i, \mfv} \defeq X(\mfv; \mfe_i)$, i.e., $\mff_{\mfv} \defeq \mfe_i$.

Now we claim that $\bma_i \defeq (a_{i,j})_{j \in J}$ lies in $\bfC(X(\mfv; \mfe_{i+1}), (X(\mfd(j); \mff_j))_j)$. To show this, it suffices to prove for all $j \neq \mfv$ that
$$
 X(j; j \mff_j) = X(\mfv; \mfe_{i+1}) \cap jX.
$$
We have $X(\mfe; \mfe_{i+1}) \cap jX = X(j; \mcm(\mfe_{i+1}, j))$. For $e \in \mfe_{i+1}$, $\lcm(\bbd(e),e_i) = e_i q$, so that $\mcm(e,j) \subseteq j \bbd^{-1}(q)$. This shows \an{$\subseteq$}. To see \an{$\supseteq$}, suppose $f \in \mfd(j) \fS$ satisfies $\bbd(f) = q$, and $\bbd(j)q = e_i q = \lcm(\bbd(e),e_i)$ for some $e \in \mfe_{i+1}$. It follows that $e_i \prec \lcm(\bbd(e),e_i))$. Hence $jf \in \mfe_{i+1}$. This shows that $j \mff_j \subseteq \mfe_{i+1}$, as desired. 
\eproof
\pars

\blemma
\label{lem:deg--desc:d}
When $\gamma(\mfe,\mfs)$ is constructed for saturated sets of the form $\mfe = \mfv \bbd^{-1}(\bbd(\mfe)) \cap \fS$ and $\mfs = \mfv \bbd^{-1}(\bbd(\mfs)) \cap \fS$, then 
$\mff_s
 =
 \mfd(s) \bbd^{-1}(\min \menge{t \in S_P}{\bbd(s)t \in \bbd(\mfs)}) \cap \fS \cup \mfd(s) \bbd^{-1}(\bbd(s)^{-1} \lcm(\bbd(s),\bbd(\mfe))) \cap \fS
$.
For every $s \in \mfs$, $X(\mfd(s); \mff_s) \in \bfX$.
\elemma
\nopar

\bproof
It is clear that \an{$\subseteq$} holds. To see \an{$\supseteq$}, if $\bbd(s) \bbd(f) \in \bbd(\mfs)$, then $sf \in \mfs$ and $f$ is minimal with that property. If $\bbd(s) \bbd(f) = \lcm(\bbd(s),\bbd(\mfe))$, then $\bbd(sf) = \bbd(e) q$, so that (UFP*) implies that $sf = e't$ for some $e' \in \mfe$ and $sf \in \mcm(s,\mfe)$.
\eproof
\pars

\blemma
\label{lem:deg_1-4Gamma}
Conditions (1$_{\Gamma}$) -- (4$_{\Gamma}$) are satisfied.
\elemma
\nopar

\bproof
Clearly, (1$_{\Gamma}$) holds.
\pars

To verify (2$_{\Gamma}$), let $c \in \fC \mfv$ with $\Vert c \Vert = L$ and $\mfs = \menge{s \in \mfv \fS}{\Vert cs \Vert = L}$. It is straightforward to see that $\mfs$ is closed under mcms, and we claim that $\mfs = \mfv \bbd^{-1}(\bbd(\mfs)) \cap \fS$. Indeed, (UFP*) implies that for $s \in \mfv \fS$, $\Vert cs \Vert = L$ if and only if $\Vert \bbd(c) \bbd(s) \Vert = L$. This implies our claim, i.e., $\mfs$ is saturated.

For (3$_{\Gamma}$), suppose that $\mfe = \mfv \bbd^{-1}(\bbd(\mfe)) \cap \fS$ and $\mfs = \mfv \bbd^{-1}(\bbd(\mfs)) \cap \fS$. If $u \in \fC^*$ and $\mff \subseteq \fS$ satisfy $u \mfe =^* \mff$, then $\mff = \mft(u) \bbd^{-1}(\bbd(\mfe)) \cap \fS$ and $u \mfs =^* \mft(u) \bbd^{-1}(\bbd(\mfs)) \cap \fS$, so that $u \gamma(\mfe,\mfs) \in \gamma(\mff, \mft(u) \bbd^{-1}(\bbd(\mfs)) \cap \fS) \C^*$.

To see (4$_{\Gamma}$), suppose that $\gamma(\mfe,\mfs) \zeta = \gamma(\mfe,\mfs')$. Then $\zeta = \coprod_{s \in \mfs} \zeta_s$. We claim that $\zeta_s = \gamma(\mff_s,\mft_s)$, where $\mft_s = \menge{t \in \fS}{st \in \mfs'}$. Indeed, all bisections appearing in $\gamma(\mfe,\mfs')$ also appear in $\gamma(\mfe,\mfs) (\coprod_{s \in \mfs} \gamma(\mff_s,\mft_s))$. Moreover, proceeding inductively, starting with maximal elements $s' \in \mfs'$, it is straightforward to check that the domains of $\gamma(\mfe,\mfs')$ and $\gamma(\mfe,\mfs) (\coprod_{s \in \mfs} \gamma(\mff_s,\mft_s))$ agree as well. Induction works because $\mff_{s'}$ only depends on $\menge{s'' \in \mfs'}{s' \prec s''}$. Now a similar argument as for Lemma~\ref{lem:char-gammaS} completes the proof of the claim. Furthermore, it is straightforward to check that $\mft_s$ is saturated and $\mcm$-closed.
\eproof
\pars

\blemma
\label{lem:deg_5Gamma}
Condition (5$_{\Gamma}$) is satisfied.
\elemma
\nopar

\bproof
We establish the criterion in Corollary~\ref{cor:5Gamma}. Let $\gamma(\mfe,\mfs), \gamma(\mfe,\mft) \in \Gamma$. We follow the construction of $\lcm_{\bfC}$ before Lemma~\ref{lem:lcm_bfC} and show that 
$$
 \lcm_{\bfC}(\gamma(\mfe,\mfs), \gamma(\mfe,\mft)) \in \gamma(\mfe,\mfs) \Gamma \cap \gamma(\mfe,\mft) \Gamma.
$$
By symmetry, it suffices to show $\lcm_{\bfC}(\gamma(\mfe,\mfs), \gamma(\mfe,\mft)) \in \gamma(\mfe,\mfs) \Gamma$. For $s \in \mfs$, define 
$$
 \mfs(s) \defeq \mfd(s) \bbd^{-1}(\bbd(s)^{-1} \lcm(\bbd(s),\bbd(\mft))) \cap \fS.
$$
Then it is straightforward to check that $\mfs(s)$ is saturated and mcm-closed, and we claim that
$$
 \lcm_{\bfC}(\gamma(\mfe,\mfs), \gamma(\mfe,\mft)) = \gamma(\mfe,\mfs) (\coprod_{s \in \mfs} \gamma(\mff_s,\mfs(s))).
$$
For the proof, first observe that $s \mfs(s) =^* \mft(s) \bbd^{-1}(\lcm(\bbd(s),\bbd(\mft))) =^* \mcm(s,\mft)$. So $\lcm_{\bfC}(\gamma(\mfe,\mfs), \gamma(\mfe,\mft))$ and $\gamma(\mfe,\mfs) (\coprod_{s \in \mfs} \gamma(\mff_s,\mfs(s)))$ both consist of bisections of the form $ss' O$, for the same elements $ss' \in \fC$. It suffices to compare domains. Suppose that $ss' = t t'$ for some $t \in \mft$. The domain corresponding to $ss'$ in $\gamma(\mfe,\mfs) (\coprod_{s \in \mfs} \gamma(\mff_s,\mfs(s)))$ is given by $X(\mfd(s'); \mff_{s'})$, where $\mff_{s'}$ consists of the minimal elements of $\menge{g \in \fS}{s' g \in \mfs(s)}$ and $\menge{g \in \fS}{s' g \in \mcm(s',\mff_s)}$. We have $s' g \in \mfs(s)$ if and only if $s'g \in s^{-1} \mcm(s,\mft)$ if and only if $ss'g \in \mcm(s,\mft)$ if and only if $ss'g \in \mcm(ss',\mft)$. Moreover, $s'g \in \mcm(s',\mff_s)$ if and only if $s'g \in \mcm(s',s^{-1}\mfs)$ or $s'g \in \mcm(s',s^{-1}(\mcm(s,\mfe)))$ if and only if $ss'g \in \mcm(ss',\mfs)$ or $ss'g \in \mcm(ss', \mcm(s,\mfe)) = \mcm(ss', \mfe)$. The domain corresponding to $ss'$ in $\lcm_{\bfC}(\gamma(\mfe,\mfs), \gamma(\mfe,\mft))$ is given by
$$
 (ss')^{-1}(s X(\mfd(s); \mff_s) \cap t X(\mfd(t); \mff_t) \cap ss'X) = (ss')^{-1}(X(ss'; \mcm(ss', t \mff_t), \mcm(ss', s \mff_s)) = X(\mfd(s'); \mfg),
$$
where $\mfg$ consists of precisely the elements of $\fS$ which are minimal among the elements $g \in \fS$ with $ss'g \in \mcm(ss',\mft)$ or $ss'g \in \mcm(ss',\mcm(\mft,\mfe)) = \mcm(ss',\mfe)$ or $ss'g \in \mcm(ss',\mfs)$ or $ss'g \in \mcm(ss',\mcm(s,\mfe)) = \mcm(ss',\mfe)$. In comparison, we see that both domains are empty if $ss' \in \mfs$, and if $ss' \notin \mfs$, then the domains coincide, as desired. 
\eproof
\pars

Putting all this together, Theorem~\ref{thm:Fn} implies the following result.
\nopar

\btheo
\label{thm:deg}
Let $P$ be a left cancellative monoid with $P^* = \gekl{1}$. Assume that $P$ is right Noetherian and admits conditional lcms. Suppose that $S_P \subseteq P$ is a finite Garside family in $P$ with $1 \notin S_P$ and assume that $(S_P^{\leq L})^{\sharp}$ is closed under left divisors for all $L \geq 1$.
\pari

Let $\fC$ be a left cancellative small category with finite $\fC^0$ equipped with a $P$-valued degree map $\bbd$ such that $\mfv \dot{\bbd}^{-1}(p) < \infty$ for all $\mfv \in \fC^0$ and $p \in P$. Let $X \subseteq \Omega_{\infty}$ be a closed invariant subspace.

Suppose that condition (F) holds, and that condition ($\bmt < \bmd$) is satisfied. 
\pars

For all natural numbers $n$, if $\fC^*(\mfv,\mfv)$ is of type ${\rm F}_n$ for all $\mfv \in \fC^0$, then $\Q(*,*)$ is of type ${\rm F}_n$. In particular, if $* = (X(\mfv;\mfe))_{\mfv \in \fV}$ for some $\fV \in \fC^0$ and $Y = \coprod_{\mfv \in \fC^0} Y_{\mfv}$, where $Y_{\mfv} = \emptyset$ if $\mfv \notin \fV$ and $Y_{\mfv} = X(\mfv;\mfe)$ if $\mfv \in \fV$, then for all natural numbers $n$, $\bmF((I_l \ltimes X)_Y^Y)$ is of type ${\rm F}_n$ if $\fC^*(\mfv,\mfv)$ is of type ${\rm F}_n$ for all $\mfv \in \fC^0$.
\etheo
\pars

Let us now identify a situation where condition ($\bmt < \bmd$) is satisfied.
\blemma
Assume that $P$ is left reversible and that for all $\mfv \in \fC^0$ and $s \in S_P$, we have $\# \mfv \bbd^{-1}(s) \mfv \geq 2$. Then for every $\mfv \in \fC^0$ and every saturated $\mfe \subseteq \mfv \fS$, we have $\partial \Omega(\mfv) = \bigcup_{\epsilon \in \mfe} \partial \Omega(\epsilon)$.
\elemma
\nopar

\bproof
By \cite[\S~3]{Spi20} and \cite[Theorem~10.5]{Spi20}, it suffices to prove for all $c \in \mfv \fC$ that $c \fC \cap \bigcup_{\epsilon \in \mfe} \epsilon \fC \neq \emptyset$. Without loss of generality, we may assume that $\bbd(\mfe) = \gekl{e_P}$. Write $\lcm(\bbd(c),\bbd(\mfe)) = \lcm(\bbd(c),e_P) = \bbd(c) p = e_P q$. By our assumption that $\# \mfd(c) \bbd^{-1}(p) \mfd(c) \geq 2$, there exists $c' \in \fC$ with $\mft(c') = \mfd(c)$ and $\bbd(c') = p$. We conclude that $\bbd(cc') = e_P q$, so (UFP*) implies that there are $\epsilon, \epsilon' \in \fC$ with $\bbd(\epsilon) = e_P$, $\bbd(\epsilon') = q$ such that $cc' = \epsilon \epsilon'$. It follows that $\epsilon \in \mfe$. Hence we deduce $cc' \in c \fC \cap \bigcup_{\epsilon \in \mfe} \epsilon \fC$. 
\eproof
\pars

\bcor
For $X = \partial \Omega$, we always have $\partial \Omega(\mfv; \mfe) = \emptyset$ unless $\mfe = \emptyset$, in which case we obtain $\partial \Omega(\mfv)$. 
\pari

For $X = \partial \Omega$, we have $\fX = \gekl{(\partial \Omega(\mfv_i))_i}$.
\ecor
\pars

\blemma
\label{lem:deg--t<d}
Assume that $P$ is left reversible and that for all $\mfv \in \fC^0$ and $s \in S_P$, we have $\# \mfv \bbd^{-1}(s) \mfv \geq 2$. Then condition ($\bmt < \bmd$) is satisfied for $X = \partial \Omega$.
\elemma
\nopar

\bproof
Consider $\gamma(\mfe,\mfs)$ for saturated $\mfs$. It suffices to treat the case $\mfe = \emptyset$, i.e., we consider $\gamma(\mfs)$ whose target is $X(\mfv)$. By assumption, there exist $s, s' \in \mfs$ with $s \neq s'$ such that $\mfv = \mft(s) = \mft(s') = \mfd(s) = \mfd(s')$. Hence $\bmm_{\bmd(\gamma(\mfs))}(X(\mfv)) \geq 2$. This implies that condition ($\bmt < \bmd$) holds.
\eproof
\pars

We obtain the following consequence of Lemma~\ref{lem:deg--t<d} and Theorem~\ref{thm:deg}.
\nopar

\bcor
\label{cor:deg}
In the situation of Theorem~\ref{thm:deg}, suppose that $X = \partial \Omega$. If $P$ is left reversible and for all $\mfv \in \fC^0$ and $s \in S_P$, we have $\# \mfv \bbd^{-1}(s) \mfv \geq 2$, then for all natural numbers $n$, $\Q(*,*)$ is of type ${\rm F}_n$ if $\fC^*(\mfv,\mfv)$ is of type ${\rm F}_n$ for all $\mfv \in \fC^0$. In particular, if $* = (\partial \Omega(\mfv))_{\mfv \in \fV}$ for some $\fV \in \fC^0$ and $Y = \coprod_{\mfv \in \fC^0} Y_{\mfv}$, where $Y_{\mfv} = \emptyset$ if $\mfv \notin \fV$ and $Y_{\mfv} = \partial \Omega(\mfv)$ if $\mfv \in \fV$, then for all natural numbers $n$, $\bmF((I_l \ltimes X)_Y^Y)$ is of type ${\rm F}_n$ if $\fC^*(\mfv,\mfv)$ is of type ${\rm F}_n$ for all $\mfv \in \fC^0$.
\ecor
\pars

\section{Examples}

\subsection{Higher rank graphs}
\label{ss:k-graphs}

Consider the special case where $P = \Zz_{\geq 0}^k$ and $S_P = \menge{(z_j) \in P}{0 \leq z_j \leq 1}$ of the setting of \S~\ref{s:Gars-deg}. Let $\epsilon_j$ be the canonical basis vector whose only non-zero component is the $j$th component, which is $1$. Assume that $\fC$ is finite. With the same assumptions as for Theorem~\ref{thm:deg}, we obtain the following:
\btheo
\label{thm:k-graphs}
Suppose that for all $\mfv \in \fC^0$ and $1 \leq j \leq k$, we have $\# \mfv \bbd^{-1}(\epsilon_j) \mfv \geq 2$. Let $X \subseteq \Omega_{\infty}$ be a closed invariant subspace and $Y = \coprod_{\mfv \in \fC^0} Y_{\mfv}$, where $Y_{\mfv} = \emptyset$ or $Y_{\mfv} = X(\mfv; \mfv \bbd^{-1}(\mfe_P))$, for some $\mfe_P \subseteq S_P$, for all $\mfv \in \fC^0$. Then $\bmF((I_l \ltimes X)_Y^Y)$ is of type ${\rm F}_{\infty}$.
\etheo
\nopar

\bproof
This follows from Theorem~\ref{thm:deg} once we verify condition ($\bmt < \bmd$), which in turn follows from Lemma~\ref{lem:deg--desc:d} because $\mfd(s) \bbd^{-1}(\bbd(s)^{-1} \lcm(\bbd(s),\bbd(\mfe))) \cap \fS = \mfd(s) \bbd(\mfe) \cap \fS$ and $\mfd(s) \bbd^{-1}(\min \menge{t \in S_P}{\bbd(s)t \in \bbd(\mfs)}) \cap \fS$ is empty for maximal elements $s \in \mfs$.
\eproof
\pars

Higher rank graphs fit naturally here. By definition, $\fC$ is a higher rank graph if, in addition to our assumptions above, we have $\fC^* = \fC^0$. In this setting, Theorem~\ref{thm:k-graphs} specializes to the following:
\nopar

\bcor
Let $\fC$ be a finite higher rank graph. Suppose that for all $\mfv \in \fC^0$ and $1 \leq j \leq k$, we have $\# \mfv \bbd^{-1}(\epsilon_j) \mfv \geq 2$. Let $X \subseteq \Omega_{\infty}$ be a closed invariant subspace and $Y = \coprod_{\mfv \in \fC^0} Y_{\mfv}$, where $Y_{\mfv} = \emptyset$ or $Y_{\mfv} = X(\mfv; \mfv \bbd^{-1}(\mfe_P))$, for some $\mfe_P \subseteq S_P$, for all $\mfv \in \fC^0$. Then $\bmF((I_l \ltimes X)_Y^Y)$ is of type ${\rm F}_{\infty}$.
\ecor
\pars
Note that we are not allowing arbitrary higher rank graphs because of the condition that $\# \mfv \bbd^{-1}(\epsilon_j) \mfv \geq 2$.

As an even more special case, we obtain an answer to a natural question left open in \cite[\S~5.3]{Mat16}.
\nopar

\bcor
\label{cor:ProdGraphs}
Suppose $\fC_j$, $1 \leq j \leq k$, are path categories of finite graphs. Assume that for all $1 \leq j \leq k$ and $\mfv \in \fC_j^0$, we have $\# \mfv \fC_j \mfv \geq 2$. Let $\cG_j \defeq I_l(\fC_j) \ltimes \partial \Omega_{\fC_j}$. Then $\bmF(\cG_1 \times \dotso \times \cG_k)$ is of type ${\rm F}_{\infty}$. 
\pari

In particular, if $\cG_j$, $1 \leq j \leq k$, are groupoids arising from irreducible one-sided shifts of finite type as in \cite{Mat15}, then $\bmF(\cG_1 \times \dotso \times \cG_k)$ is of type ${\rm F}_{\infty}$.
\ecor
\nopar

\bproof
The first part follows from Theorem~\ref{thm:k-graphs}~(B), applied to $\fC \defeq \fC_1 \times \dotso \fC_k$, once we observe that $\cG_1 \times \dotso \times \cG_k \cong I_l(\fC) \ltimes \partial \Omega_{\fC}$. 
\pars

For the second part, let $\fC_j$, $1 \leq j \leq k$, be finite (directed) graphs describing the given irreducible one-sided shifts of finite type. The argument in \cite[\S~6.5, p. 67]{Mat15} shows that we may without loss of generality assume that for all $1 \leq j \leq k$ and $\mfv \in \fC_j^0$, we have $\# \mfv \fC_j \mfv \geq 2$. Now the second part follows from the first.
\eproof
\pars

We can also specialize to one vertex higher rank graphs, i.e., higher rank graphs with just one vertex. In this case, we are actually dealing with monoids given by particular presentations. Let $\hat{P}$ be such a monoid with a degree map $\bbd: \: \hat{P} \to P$. We denote the identity element of $\hat{P}$ by $1_{\hat{P}}$.
\nopar

\bcor
\label{cor:OneVertex}
Let $\hat{P}$ be a monoid given by a finite one vertex higher rank graph. Suppose that $\bbd^{-1}(\epsilon_j) \geq 2$ for all $1 \leq j \leq k$. Let $X \subseteq \Omega_{\infty}$ be a closed invariant subspace and $Y = X(1_{\hat{P}}; \bbd^{-1}(\mfe_P))$ for some $\mfe_P \subseteq S_P$. Then $\bmF((I_l \ltimes X)_Y^Y)$ is of type ${\rm F}_{\infty}$.
\ecor
\pars

\subsection{Zappa-Sz{\'e}p products}
\label{ss:ZS}

The next class of examples is given by Zappa-Sz{\'e}p products of a small category and a groupoid. In related contexts, such examples have been considered in more special situations in \cite{Bri05,Nek04,Nek18a,EP,Sta,LRRW14,LRRW18,LY,ABRW,BKQS,Wit19}. For our purposes, Zappa-Sz{\'e}p products allow to adjoin invertible elements to a small category while keeping the key properties which are needed to establish finiteness properties for the corresponding topological full groups of groupoids attached to left regular representations.

The setting is as follows: Suppose that a groupoid $\fG$ acts on a left cancellative small category $\fC$ via a map $\pi: \: \fC^0 \to \fG^0$. The action is such that $g.c$ is defined if $\rms(g) = \pi(\mft(c))$, and then $\pi(\mft(g.c)) = \rmr(g)$. We require the usual axioms for actions, i.e., $\mfx.c = c$ and $g.(h.c) = gh.c$ for all $\mfx \in \fG^0$, $g, h \in \fG$ and $c \in \fC$, whenever these equations make sense. Note that we do not require $\fG$ to act by automorphisms of $\fC$ as a category. 

We also require that $\mft(g.c) = g.\mft(c)$. This means that the $\fG$-action restricts to an action $\fG \curvearrowright \fC^0$ since $g. \mfv = g. \mft(\mfv) = \mft(g. \mfv) \in \fC^0$. Moreover, set $\fG {}_{\rms} \hspace{-0.25em} \times_{\mft} \fC \defeq \menge{(g,a) \in \fG \times \fC}{\rms(g) = \pi(\mft(a)}$ and let $\varphi: \: \fG {}_{\rms} \hspace{-0.25em} \times_{\mft} \fC \to \fG$ satisfy
\nopar

\begin{itemize}
\item $\rms(\varphi(g,a)) = \pi(\mfd(a))$,
\item $\varphi(gh,a) = \varphi(g,h.a) \varphi(h,a)$,
\item $\varphi(g,\mfv) = g$,
\item $\varphi(g,a). \mfd(a) = \mfd(g.a)$,
\item $g.(ab) = (g.a) (\varphi(g,a).b)$,
\item $\varphi(g,ab) = \varphi(\varphi(g,a),b)$,
\end{itemize}
for all $g, h \in \fG$, $\mfv \in \fC^0$ and $a, b \in \fC$, whenever these equations make sense.
\pars

We call an action $\fG \curvearrowright \fC$ together with a $\fG$-valued map $\varphi$ satisfying the conditions above a self-similar action of $\fG$ on $\fC$. Note that part of the conditions means that $\varphi$ is a cocycle. 

\bdefin
The Zappa-Sz{\'e}p product $\fC \bowtie \fG$ is the small category with underlying set 
$$
 \fC {}_{\mfd} \hspace{-0.25em} \times_{\rmr} \fG \defeq \menge{(a,g) \in \fC \times \fG}{\pi(\mfd(a)) = \rmr(g)},
$$set of objects $\menge{(\mfv,\pi(\mfv))}{\mfv \in \fC^0} \subseteq \fC^0 \times \fG^0$, and multiplication $(a,g)(b,h) \defeq (a (g.b),\varphi(g,b) h)$.
\edefin
\nopar

We refer to \cite{Wit19,BKQS} for more details about the construction of Zappa-Sz{\'e}p products. In the following, we denote $\fC \bowtie \fG$ by $\fD$.
\pars

\bremark
By first forming the transformation groupoid for the restricted action $\fG \curvearrowright \fC^0$ and then the Zappa-Sz{\'e}p product, we may assume without loss of generality that $\fG^0 = \fC^0$ (as in \cite{Wit19}). 
\eremark

In the following, given $a \in \fC$, we write $\ti{a} \defeq (a,\pi(\mfd(a))) \in \fD$, and identify $\fC$ with a subcategory of $\fD$ via $\fC \into \fD, \, a \ma \ti{a}$.
\blemma
\label{lem:fCfD}
\begin{enumerate}
\item[(i)] The $\fG$-action on $\fC$ restricts to an action $\fG \curvearrowright \fC^*$.
\item[(ii)] For $a, b \in \fC$ and $g \in \fG$, we have $a \prec b$ if and only if $g.a \prec g.b$.
\item[(iii)] We have $(a,g) \preceq (b,h)$ if and only if $a \preceq b$, and $(a,g) \prec (b,h)$ if and only if $a \prec b$.
\item[(iv)] If $\fC$ is left cancellative/finitely aligned/left Noetherian/right Noetherian, then so is $\fD$. 
\item[(v)] For $a, b \in \fC$, we have $\mcm(a,b) \, \ti{ } = \mcm(\ti{a},\ti{b})$ whenever $\mcm(a,b)$ exists in $\fC$. In particular, $\fD$ admits disjoint mcms if and only if $\fC$ admits disjoint mcms.
\end{enumerate}
\elemma
\nopar

\bproof
(i) If $a \in \fC^*$, then there exists $b \in \fC$ with $ab = \mfv \in \fC^0$. Then $g.\mfv = g.(ab) = (g.a)(\varphi(g,a).b)$ and $g.\mfv \in \fC^0$ imply that $g.a \in \fC^*$.
\pars

(ii) If $b = ax$, then $g.b = (g.a)(\varphi(g,a).x)$ implies that $g.a \prec g.b$.

(iii) This is easy to see.

(iv) It is easy to see that $\fD$ is left cancellative if $\fC$ is left cancellative. 
\pari

Suppose that $\fC$ is finitely aligned. If $a \fC \cap b \fC = \bigcup_i c_i \fC$, then we claim that $\ti{a} \fD \cap \ti{b} \fD = \bigcup_i \ti{c}_i \fD$. This shows that $\fD$ is finitely aligned because every principal right ideal of $\fD$ is of the form $\ti{a} \fD$. To prove the claim, note that \an{$\supseteq$} is clear, and for \an{$\subseteq$}, suppose that $(c,g) = (a,\pi(\mfd(a)))(a',g) = (b,\pi(\mfd(b)))(b',g)$. Then $c = aa' = bb'$ and thus $c = c_i x$ for some $x \in \fC$. It follows that $(c,g) = (c_i x,g) = (c_i, \pi(\mfd(c_i))) (x,g)$.

It follows from (iii) that if $\fC$ is left Noetherian, then so is $\fD$. Suppose that $\fD$ is not right Noetherian. Then we have an infinite descending chain $\dotso \tiprec (a_3,g_3) \tiprec (a_2,g_2) \tiprec (a_1,g_1)$. So $(a_1,g_1) = (b_2,h_2)(a_2,g_2) = (b_2 (h_2.a_2),\varphi(h_2,a_2) g_2)$. We deduce $h_2.a_2 \tiprec a_1$. Proceeding inductively, we obtain an infinite descending chain $\dotso \tiprec a'_3 \tiprec a'_2 \tiprec a_1$ in $\fC$. Thus $\fC$ is not right Noetherian.
\pars

(v) We have seen in (iv) that $a \fC \cap b \fC = \bigcup_i c_i \fC$ implies that $\ti{a} \fD \cap \ti{b} \fD = \bigcup_i \ti{c}_i \fD$. Now the first part of (v) follows because of (iii). The second part of (v) is a consequence of the first. 
\eproof
\pars

\bremark
If $\fC$ is right cancellative, $\fD$ does not need to be right cancellative (see \cite{LW}).
\eremark

Now suppose that $\fS$ is a subset of $\fC$. We set $\ti{\fS} \defeq \menge{(s,\pi(\mfd(s)))}{s \in \fS} \subseteq \fD$.
\nopar
\blemma
\label{lem:SC-->SD}
\begin{enumerate}
\item[(i)] $\fD^* \ti{\fS} \subseteq \ti{\fS} \fD^*$ if and only if $\fC^* \fS \subseteq \fS \fC^*$ and $\fG.\fS \subseteq \fS \fC^*$.
\item[(ii)] If $\fG.\fS \subseteq \fS \fC^*$, then $\fS^{\sharp}$ is closed under right divisors if and only if $\ti{\fS}^{\sharp}$ is closed under right divisors.
\item[(iii)] $\ti{\fS}$ is a Garside family in $\fD$ if and only if $\fS$ is a Garside family of $\fC$ and $\fG.\fS \subseteq \fS \fC^*$.
\item[(iv)] $(\ti{\fS}^{\leq L})^{\sharp}$ is closed under left divisors if and only if $(\fS^{\leq L})^{\sharp}$ is closed under left divisors.
\item[(v)] $\ti{\fS}^{\sharp}$ is closed under mcms if and only if $\fS^{\sharp}$ is closed under mcms.
\item[(vi)] $\ti{\fS}$ is $=^*$-transverse if and only if $\fS$ is $=^*$-transverse. $\ti{\fS}$ is locally bounded if and only if $\fS$ is locally bounded.
\end{enumerate}
\elemma
\bproof
(i) Given $(u,g) \in \fD^*$ and $(s,\pi(\mfd(s))) \in \ti{\fS}$, we have $(u,g)(s,\pi(\mfd(x))) = (u (g.s), \varphi(g,s))$, and that element lies in $\ti{\fS} \fD^*$ if and only if $u (g.s) \in \fS \fC^*$. The latter holds for all $u \in \fC^*$, $g \in \fG$ and $s \in \fS$ (whenever it makes sense) if and only if $\fC^* \fS \subseteq \fS \fC^*$ and $\fG.\fS \subseteq \fS \fC^*$.
\pars

(ii) Given $(a,g)$ and $(b,h)$, we have $(a (g.b),\varphi(g,b)h) = (a,g) (b,h) \in \ti{\fS} \fD^*$ if and only if $a (g.b) \in \fS \fC^*$. The latter implies $g.b \in \fS \fC^*$ if $\fS^{\sharp}$ is closed under right divisors. Thus $b = g^{-1}.(su)$ for some $s \in \fS$ and $u \in \fC^*$. We conclude that $b = g^{-1}.(su) = (g^{-1}.s) (\varphi(g^{-1},s).u) \in \fS \fC^* \fC^*$. Conversely, suppose that $\ti{\fS}^{\sharp}$ is closed under right divisors. Given $su \in \fS \fC^*$ and $a,b \in \fC$ with $ab = su$, then $\ti{a} \ti{b} = \ti{s} \ti{u} \in \ti{\fS} \fD^*$. It follows that $\ti{b} \in \ti{\fS}^{\sharp}$, so that $b \in \fS^{\sharp}$. 

(iii) This follows from the characterization of Garside families in \cite[Chapter~IV, Proposition~1.24]{Deh15}: First of all, if $\ti{\fS}$ is Garside, then we must have $\fG.\fS \subseteq \fS \fC^*$ by (i). Now $\fS^{\sharp}$ generates $\fC$ if and only if $\ti{\fS}^{\sharp}$ generates $\fD$, and if $\fG.\fS \subseteq \fS \fC^*$, then (ii) tells us that $\fS^{\sharp}$ is closed under right divisors if and only if $\ti{\fS}^{\sharp}$ is closed under right divisors. Finally, $\fS$-heads exist in $\fC$ if and only if $\ti{\fS}$-heads exist in $\fD$ because $s \preceq a$ if and only if $\ti{s} \preceq (a,g)$ for all $g \in \rmr^{-1}(\pi(\mfd(a)))$.

(iv) $(a,g)(b,h) = (a (g.b),\varphi(g,b)h)$ lies in $(\ti{\fS}^{\leq L})^{\sharp}$ if and only if $a (g.b)$ lies in $(\fS^{\leq L})^{\sharp}$, and $(a,g)$ lies in $(\ti{\fS}^{\leq L})^{\sharp}$ if and only if $a$ lies in $(\fS^{\leq L})^{\sharp}$. Now the claim follows.

(v) follows from Lemma~\ref{lem:fCfD}~(v).

(vi) We have $s \in t \fC^*$ if and only if $\ti{s} \in \ti{t} \fD^*$. The second claim follows from Lemma~\ref{lem:fCfD}~(iii).
\eproof
\pars

\bremark
Lemmas~\ref{lem:fCfD} and \ref{lem:SC-->SD} imply that if $\fC$ is a finitely aligned, left cancellative, countable small category, $\fS$ is a Garside family in $\fC$ with $\fS \cap \fC^* = \emptyset$ which is $=^*$-transverse and locally bounded, and if $\fG \curvearrowright \fC$ is a self-similar action of a countable groupoid $\fG$ with $\fG.\fS \subseteq \fS \fC^*$, then \cite[Theorem~B]{Li21a} applies to $\fD = \fC \bowtie \fG$, $\ti{\fS}$ and yields a description of all closed invariant subspaces of $I_l(\fD) \ltimes \Omega_{\fD}$.
\eremark

The following is a direct consequence of Lemma~\ref{lem:fCfD}~(v).
\nopar

\blemma
Suppose that $\fC$ is finitely aligned and left cancellative. We have an isomorphism of semilattices $\cJ_{\fC} \isom \cJ_{\fD}, \, x \fC \ma \ti{x} \fD$. This induces an $I_l(\fC)$-equivariant homeomorphism $\Omega_{\fC} \isom \Omega_{\fD}, \, \chi \ma \ti{\chi}$. In particular, if $\ti{X}$ is a closed invariant subspace of $\Omega_{\fD}$, then $X$ is a closed invariant subspace of $\Omega_{\fC}$. In this situation, we obtain an embedding $I_l(\fC) \ltimes X \into I_l(\fD) \ltimes \ti{X}, \, [c d^{-1},\chi] \ma [\ti{c} \ti{d}^{-1},\ti{\chi}]$ of $I_l(\fC) \ltimes X$ as an open subgroupoid of $I_l(\fD) \ltimes \ti{X}$. This in turn induces a map of compact open bisections, denoted by $a \ma \ti{a}$, and hence a map $\bfC \to \bfD, \, \bma \ma \ti{\bma}$.
\elemma
Note that in general, the image $\ti{X}$ of a closed invariant subspace $X$ of $\Omega_{\fC}$ is not $I_l(\fD)$-invariant.
\pars

Now assume that $\fX$, $\bfX$, $\Gamma$, $\bGamma$, $* \in \bfX$, $\C$ and $\Q$ be as in \S~\ref{s:GarsCat-TFG} (for $\fC$). Let $\ti{\fX} = \menge{\ti{U}}{U \in \fX}$, $\ti{\bfX} = \menge{\ti{\bmU}}{\bmU \in \bfX}$, $\ti{\Gamma} \defeq \menge{\ti{\gamma}}{\gamma \in \Gamma}$ and define $\ti{\bGamma}$ correspondingly. For $\ti{*} \in \ti{\bfX}$, we let $\D$ be the component of $\ti{*}$ in the smallest subcategory $\ti{\Pi} \defeq \langle \ti{\Gamma}, \bfD_{\ti{\bfX}}^* \rangle$ of $\bfD$ containing $\ti{\Gamma}$ and $\bfD_{\ti{\bfX}}^*$. Let $\Q_{\D}$ be the enveloping groupoid of $\D$.
\nopar

\blemma
Suppose that for all $g \in \fG$ and $\mfv$, $\mfe$ with $X(\mfv;\mfe) \in \fX$, we have $X(g.\mfv;g.\mfe) \in \fX$, and that for all $g \in \fG$ and $\gamma = (a_i) \in \Gamma$, we have $g.\gamma \defeq (g.a_i) \in \Gamma \C^*$.
\nopar

\begin{enumerate}
\item[(i)] If (F) holds for $\bfD$, $\ti{\fX}$, then (F) also holds for $\bfC$, $\fX$. 
\item[(ii)] If $\fX$, $\Gamma$ satisfy condition (St), then so do $\ti{\fX}$, $\ti{\Gamma}$.
\item[(iii)] We have $\D = \spkl{\ti{\Gamma}} \D^* \cup \D^*$. If $\C$ admits lcms, then so does $\D$.
\item[(iv)] If ($\bmt < \bmd$) is satisfied for $\fX$, $\Gamma$, then ($\bmt < \bmd$) is also satisfied for $\ti{\fX}$, $\ti{\Gamma}$.
\end{enumerate}
\elemma
Note that given a bisection of the form $a X(\mfv; \mfe)$, then we define $g.(a X(\mfv; \mfe)) \defeq (g.a) X(\varphi(g,a).\mfv; \varphi(g,a).\mfe)$. 
\bproof
(i) and (ii) are straightforward to check.
\pars

(iii) The first claim is straightforward to check. Now we claim that if $\gamma = \lcm_{\C}(\alpha,\beta)$, then $\ti{\gamma} = \lcm_{\D}(\ti{\alpha},\ti{\beta})$. It is clear that $\ti{\alpha} \preceq \ti{\gamma}$ and $\ti{\beta} \preceq \ti{\gamma}$. Now suppose that we are given a common multiple of $\ti{\alpha}$ and $\ti{\beta}$ of the form $\ti{\alpha} \delta = \ti{\beta} \varepsilon$. As $\D = \langle \ti{\Gamma} \rangle \D^* \cup \D^*$, we may assume that $\delta, \varepsilon \in \langle \ti{\Gamma} \rangle$, i.e., $\delta = (\alpha') \ti{ }$ and $\varepsilon = (\beta') \ti{ }$ for some $\alpha', \beta' \in \langle \Gamma \rangle$. It follows that $(\alpha \alpha') \ti{ } = (\beta \beta') \ti{ }$ and thus $\alpha \alpha' = \beta \beta'$. The last conclusion is justified because $(\ti{a},\ti{\chi}) \sim (\ti{b}, \ti{\chi})$ implies that $\ti{a} \ti{x} = \ti{b} \ti{x}$ for some $x \in \fC$, so that $ax = bx$ and thus $(a,\chi) \sim (b,\chi)$. Now it follows that $\gamma \preceq \alpha \alpha' = \beta \beta'$ and thus $\ti{\gamma} \preceq \ti{\alpha} (\alpha') \ti{ } = \ti{\beta} (\beta') \ti{ }$, as desired.

(iv) This follows from the observations that $\bmt(\ti{\gamma}) = \bmt(\gamma) \ti{ }$ and $\bmd(\ti{\gamma}) = \bmd(\gamma) \ti{ }$.
\eproof
\pars

Our findings above motivate the following terminology.
\nopar

\bdefin
We say that condition (Inv) is satisfied if $\fG.\fS \subseteq \fS \fC^*$, for all $g \in \fG$ and $\mfv$, $\mfe$ with $X(\mfv;\mfe) \in \fX$, we have $X(g.\mfv;g.\mfe) \in \fX$, and for all $g \in \fG$ and $\gamma = (a_i) \in \Gamma$, we have $g.\gamma \defeq (g.a_i) \in \Gamma \C^*$.
\edefin
\pars

With the straightforward observation that cancellation up to $=^*$ passes from $\fD$ to $\fC$, we obtain the following:
\nopar

\btheo
\label{thm:ZS}
Assume that $\fC$ is a finitely aligned left cancellative small category with finite $\fC^0$. Let $\fS \subseteq \fC$ be a locally finite Garside family such that $(\fS^{\leq L})^{\sharp}$ is closed under left divisors for all $L \geq 1$.
\pari

Let $\fG \curvearrowright \fC$ be a self-similar action and form $\fD \defeq \fC \bowtie \fG$. Let $\ti{X} \subseteq \Omega_{\fD,\infty}$ be a closed invariant subspace. Assume that condition (Inv) is satisfied, $\fD$ is right cancellative up to $=^*$, that condition (F) holds for $\fD$, that $\fX$, $\Gamma$ satisfy condition (St) and that $\C$ admits lcms. The latter is for instance the case when $\fC$ is right Noetherian, admits disjoint mcms, and if condition (LCM) holds. Further suppose that ($\bmt < \bmd$) holds for $\Gamma(*)$.
\pars

Then for all natural numbers $n$, $\Q_{\D}(\ti{*},\ti{*})$ is of type ${\rm F}_n$ if $\fD^*(\mfw,\mfw)$ is of type ${\rm F}_n$ for all $\mfw \in \fD^0$. 
\pari

In particular, if $* = (X(\mfv;\mfe))_{\mfv \in \fV}$ for some $\fV \in \fC^0$ and $Y = \coprod_{\mfv \in \fC^0} Y_{\mfv}$, where $Y_{\mfv} = \emptyset$ if $\mfv \notin \fV$ and $Y_{\mfv} = X(\mfv;\mfe)$ if $\mfv \in \fV$, then for all natural numbers $n$, $\bmF((I_l(\fD) \ltimes \ti{X})_{\ti{Y}}^{\ti{Y}})$ is of type ${\rm F}_n$ if $\fD^*(\mfw,\mfw)$ is of type ${\rm F}_n$ for all $\mfw \in \fD^0$.
\etheo
\nopar

Note that condition (F) for $\fD$ is related to faithfulness of the action of $\fG$ on $\fC$.

Let us describe the groups $\fD^*(\mfw,\mfw)$ in terms of $\fC^*$ and $\fG$. 
\nopar

\bremark
\label{rem:ZS-Stab}
We have $\fD^*(\ti{\mfv}, \ti{\mfv}) = \menge{(u,g) \in \fD}{u \in \fC^*, \, g \in \fG, \, \mfd(u) = g.\mfv, \, \mft(u) = \mfv}$. In particular, if $\fC^* = \fC^{*,0}$, then $\fD^*(\ti{\mfv}, \ti{\mfv}) = \fC^*(\mfv,\mfv) \bowtie {\rm St}(\fG,\mfv)$, where ${\rm St}(\fG,\mfv) = \{ g \in \fG_{\pi(\mfv)}^{\pi(\mfv)}: \: g.\mfv = \mfv \}$. If we even have $\fC^* = \fC^0$, then $\fD^*(\ti{\mfv}, \ti{\mfv}) = {\rm St}(\fG,\mfv)$.
\eremark
\pars

Here is a more concrete class of examples where our findings above apply.
\bex
\label{ex:ZS}
Assume that we are in the setting of \S~\ref{s:Gars-deg} of a cancellative small category $\fC$ with a degree map $\bbd: \: \fC \to P$. Assume that the self-similar action $\fG \curvearrowright \fC$ satisfies $\bbd(g.a) = \bbd(a)$ for all $g \in \fG$ and $a \in \fC$. If $P$ is right cancellative, then $\fD \defeq \fC \bowtie \fG$ is automatically right cancellative up to $=^*$, and for the Garside family $\fS$ from \S~\ref{s:Gars-deg}, we have $\fG.\fS \subseteq \fS \fC^*$. Moreover, if $\fX$ and $\Gamma$ are constructed as in \S~\ref{s:Gars-deg}, then for all $g \in \fG$ and $\mfv$, $\mfe$ with $X(\mfv;\mfe) \in \fX$, we have $X(g.\mfv;g.\mfe) \in \fX$, and that for all $g \in \fG$ and $\gamma = (a_i) \in \Gamma$, we have $g.\gamma \defeq (g.a_i) \in \Gamma \C^*$. This means that if $\fC$ satisfies the same assumptions as for Theorem~\ref{thm:deg}, and condition (F) holds, then we obtain the following conclusions, with the same notations as in Theorem~\ref{thm:ZS}:
\nopar

\begin{enumerate}
\item[(I)] If ($\bmt < \bmd$) holds for $\Gamma(*)$, then for all natural numbers $n$, $\Q_{\D}(\ti{*},\ti{*})$ is of type ${\rm F}_n$ if $\fD^*(\mfw,\mfw)$ is of type ${\rm F}_n$ for all $\mfw \in \fD^0$. In particular, if $* = (X(\mfv;\mfe))_{\mfv \in \fV}$ for some $\fV \in \fC^0$ and $Y = \coprod_{\mfv \in \fC^0} Y_{\mfv}$, where $Y_{\mfv} = \emptyset$ if $\mfv \notin \fV$ and $Y_{\mfv} = X(\mfv;\mfe)$ if $\mfv \in \fV$, then for all natural numbers $n$, $\bmF((I_l(\fD) \ltimes \ti{X})_{\ti{Y}}^{\ti{Y}})$ is of type ${\rm F}_n$ if $\fD^*(\mfw,\mfw)$ is of type ${\rm F}_n$ for all $\mfw \in \fD^0$.
\item[(II)] Suppose that $X = \partial \Omega$. If $P$ is left reversible and for all $\mfv \in \fC^0$ and $s \in S_P$, we have $\# \mfv \bbd^{-1}(s) \mfv \geq 2$, then for all natural numbers $n$, $\Q_{\D}(\ti{*},\ti{*})$ is of type ${\rm F}_n$ if $\fD^*(\mfw,\mfw)$ is of type ${\rm F}_n$ for all $\mfw \in \fD^0$. In particular, if $* = (\partial \Omega(\mfv))_{\mfv \in \fV}$ for some $\fV \in \fC^0$ and $Y = \coprod_{\mfv \in \fC^0} Y_{\mfv}$, where $Y_{\mfv} = \emptyset$ if $\mfv \notin \fV$ and $Y_{\mfv} = \partial \Omega(\mfv)$ if $\mfv \in \fV$, then for all natural numbers $n$, $\bmF((I_l(\fD) \ltimes \ti{X})_{\ti{Y}}^{\ti{Y}})$ is of type ${\rm F}_n$ if $\fD^*(\mfw,\mfw)$ is of type ${\rm F}_n$ for all $\mfw \in \fD^0$.
\item[(III)] If $P = \Zz_{\geq 0}^k$ and $S_P$, $\epsilon_j$ are as in \S~\ref{ss:k-graphs}, for all $\mfv \in \fC^0$ and $1 \leq j \leq k$, we have $\# \mfv \bbd^{-1}(\epsilon_j) \mfv \geq 2$, and if $X \subseteq \Omega_{\infty}$ is a closed invariant subspace, then for all natural numbers $n$, $\Q_{\D}(\ti{*},\ti{*})$ is of type ${\rm F}_n$ if $\fD^*(\mfw,\mfw)$ is of type ${\rm F}_n$ for all $\mfw \in \fD^0$. In particular, if $* = (X(\mfv;\mfe))_{\mfv \in \fV}$ for some $\fV \in \fC^0$ and $Y = \coprod_{\mfv \in \fC^0} Y_{\mfv}$, where $Y_{\mfv} = \emptyset$ if $\mfv \notin \fV$ and $Y_{\mfv} = X(\mfv;\mfe)$ if $\mfv \in \fV$, then for all natural numbers $n$, $\bmF((I_l(\fD) \ltimes \ti{X})_{\ti{Y}}^{\ti{Y}})$ is of type ${\rm F}_n$ if $\fD^*(\mfw,\mfw)$ is of type ${\rm F}_n$ for all $\mfw \in \fD^0$.
\end{enumerate}
In particular, (III) covers self-similar groups as in \cite{Nek04,Nek18a}, and, in combination with Remark~\ref{rem:ZS-Stab}, yields a generalization of \cite[Theorem~4.15]{SWZ}. (III) also covers self-similar actions on graphs as in \cite{EP,LRRW14,LRRW18} and self-similar actions on higher rank graphs as in \cite{LY,ABRW}. To give a concrete example, consider groupoid models for Katsura algebras as discussed in \cite[\S~18]{EP}. In this case, $\fC$ is the path category a finite graph, so we are in case (III) with $P = \Zz_{\geq 0}$, and $\fG$ is the group $\Zz$. Because of Remark~\ref{rem:ZS-Stab}, and because quotients of subgroups of $\Zz$ are of type ${\rm F}_{\infty}$, condition (F) is not needed in this case, and we obtain the same conclusions as in (III) if for all $\mfv \in \fC^0$, we have $\# \mfv \bbd^{-1}(1) \mfv \geq 2$.
\eex
\pars

\textbf{Conflict of interest statement:} There is no conflict of interest.

\textbf{Data availability statement:} This manuscript has no associated data.


\begin{thebibliography}{99}

\bibitem{ABRW} Z. \textsc{Afsar}, N. \textsc{Brownlowe}, J. \textsc{Ramagge} and M.F. \textsc{Whittaker}, \emph{$C^*$-algebras of self-similar actions of groupoids on higher-rank graphs and their equilibrium states}, preprint, arXiv:1910.02472.

\bibitem{BKQS} E. \textsc{B{\'e}dos}, S. \textsc{Kaliszewski}, J. \textsc{Quigg} and J. \textsc{Spielberg}, \emph{On finitely aligned left cancellative small categories, Zappa-Sz{\'e}p products and Exel-Pardo algebras}, Theory Appl. Categ. \emph{33} (2018), Paper No. 42, 1346--1406.

\bibitem{Bes} D. \textsc{Bessis}, \emph{Finite complex reflection arrangements are $K(\pi,1)$}, Ann. of Math. (2) \emph{181} (2015), no. 3, 809--904.


\bibitem{Bjo95} A. \textsc{Bj{\"o}rner}, \emph{Topological methods}, in Handbook of combinatorics, Vol. 2, 1819--1872, Elsevier Sci. B. V., Amsterdam, 1995.

\bibitem{Bjo03} A. \textsc{Bj{\"o}rner}, \emph{Nerves, fibers and homotopy groups}, J. Combin. Theory Ser. A \emph{102} (2003), no. 1, 88--93.



\bibitem{Bri04} M.B. \textsc{Brin}, \emph{Higher dimensional Thompson groups}, Geom. Dedicata \emph{108} (2004), 163--192.

\bibitem{Bri05} M.B. \textsc{Brin}, \emph{On the Zappa-Sz{\'e}p product}, Comm. Algebra \emph{33} (2005), no. 2, 393--424.

\bibitem{Bro} K.S. \textsc{Brown}, \emph{Finiteness properties of groups}, J. Pure Appl. Algebra \emph{44} (1987), no. 1-3, 45--75.

\bibitem{BG} K.S. \textsc{Brown} and R. \textsc{Geoghegan}, \emph{An infinite-dimensional torsion-free ${\rm FP}_{\infty}$ group}, Invent. Math. \emph{77} (1984), no. 2, 367--381.

\bibitem{BRRW} N. \textsc{Brownlowe}, J. \textsc{Ramagge}, D. \textsc{Robertson} and M.F. \textsc{Whittaker}, \emph{Zappa-Sz{\'e}p products of semigroups and their $C^*$-algebras}, J. Funct. Anal. \emph{266} (2014), no. 6, 3937--3967.


\bibitem{CFP} J.W. \textsc{Cannon}, W.J. \textsc{Floyd} and W.R. \textsc{Parry}, \emph{Introductory notes on Richard Thompson's groups}, Enseign. Math. (2) \emph{42} (1996), no. 3-4, 215--256.

\bibitem{Cun} J. \textsc{Cuntz}, \emph{Simple $C^*$-algebras generated by isometries}, Comm. Math. Phys. \emph{57} (1977), no. 2, 173--185.

\bibitem{Deh15} J. \textsc{Dehornoy}, F. \textsc{Digne}, E. \textsc{Godelle}, D. \textsc{Krammer}, J. \textsc{Michel}, \emph{Foundations of Garside theory}, EMS Tracts in Mathematics, 22, European Mathematical Society, Z{\"u}rich, 2015.

\bibitem{EP} R. \textsc{Exel} and E. \textsc{Pardo}, \emph{Self-similar graphs, a unified treatment of Katsura and Nekrashevych $C^*$-algebras}, Adv. Math. \emph{306} (2017), 1046--1129.

\bibitem{FMWZ} M.G. \textsc{Fluch}, M. \textsc{Marschler}, S. \textsc{Witzel}, M.C.B. \textsc{Zaremsky}, \emph{The Brin-Thompson groups $sV$ are of type ${\rm F}_{\infty}$}, Pacific J. Math. \emph{266} (2013), no. 2, 283--295.

%\bibitem{FS} N.J. \textsc{Fowler} and A. \textsc{Sims}, \emph{Product systems over right-angled Artin semigroups}, Trans. Amer. Math. Soc. \emph{354} (2002), no. 4, 1487--1509.

\bibitem{Geo} R. \textsc{Geoghegan}, \emph{Topological methods in group theory}, Graduate Texts in Mathematics, 243, Springer, New York, 2008.

\bibitem{Hat} A. \textsc{Hatcher}, \emph{Algebraic topology}, Cambridge University Press, Cambridge, 2002.

\bibitem{HRSW} R. \textsc{Hazlewood}, I. \textsc{Raeburn}, A. \textsc{Sims}, S.B.G. \textsc{Webster}, \emph{Remarks on some fundamental results about higher-rank graphs and their $C^*$-algebras}, Proc. Edinb. Math. Soc. (2) \emph{56} (2013), no. 2, 575--597.

\bibitem{Hig} G. \textsc{Higman}, \emph{Finitely presented infinite simple groups}, Notes on Pure Mathematics, No. 8, Department of Pure Mathematics, Department of Mathematics, I.A.S. Australian National University, Canberra, 1974.

\bibitem{HO} J. \textsc{Huang} and D. \textsc{Osajda}, \emph{Helly meets Garside and Artin}, Invent. Math. \emph{225} (2021), no. 2, 395--426.

\bibitem{JM} K. \textsc{Juschenko} and N. \textsc{Monod}, \emph{Cantor systems, piecewise translations and simple amenable groups}, Ann. of Math. (2) \emph{178} (2013), no. 2, 775--787.

\bibitem{LRRW14} M. \textsc{Laca}, I. \textsc{Raeburn}, J. \textsc{Ramagge} and M.F. \textsc{Whittaker}, \emph{Equilibrium states on the Cuntz-Pimsner algebras of self-similar actions}, J. Funct. Anal. \emph{266} (2014), no. 11, 6619--6661.

\bibitem{LRRW18} M. \textsc{Laca}, I. \textsc{Raeburn}, J. \textsc{Ramagge} and M.F. \textsc{Whittaker}, \emph{Equilibrium states on operator algebras associated to self-similar actions of groupoids on graphs}, Adv. Math. \emph{331} (2018), 268--325.

\bibitem{LV} M.V. \textsc{Lawson} and A. \textsc{Vdovina}, \emph{A generalization of higher rank graphs}, to appear in Bull. Aust. Math. Soc., arXiv:2104.09421.

\bibitem{LW} M.V. \textsc{Lawson} and A.R. \textsc{Wallis}, \emph{A correspondence between a class of monoids and self-similar group actions II}, 
Internat. J. Algebra Comput. \emph{25} (2015), no. 4, 633--668.

\bibitem{Li17} X. \textsc{Li},
  \emph{Partial transformation groupoids attached to graphs and semigroups}, Int. Math. Res. Not. \emph{2017}, no. 17, 5233--5259.

\bibitem{Li20} X. \textsc{Li}, \emph{Every classifiable simple C*-algebra has a Cartan subalgebra}, Invent. Math. \emph{219} (2020), 653--699. 

\bibitem{Li21a} X. \textsc{Li}, \emph{Left regular representations of Garside categories I. $C^*$-algebras and groupoids}, to appear in Glasg. Math. J., arXiv:2110.04501.

\bibitem{Li22} X. \textsc{Li}, \emph{Ample groupoids, topological full groups, algebraic K-theory spectra and infinite loop spaces}, preprint, arXiv:2209.08087.

%\bibitem{LN} X. \textsc{Li} and M.D. \textsc{Norling}, \emph{Independent resolutions for totally disconnected dynamical systems I: Algebraic case}, 
%J. Algebra \emph{424} (2015), 98--125.

\bibitem{LY} H. \textsc{Li} and D. \textsc{Yang}, \emph{Self-Similar $k$-Graph $C^*$-Algebras}, Int. Math. Res. Not. \emph{2021}, no. 15, 11270--11305. 

\bibitem{Mat15} H. \textsc{Matui}, \emph{Topological full groups of one-sided shifts of finite type}, J. Reine Angew. Math. \emph{705} (2015), 35--84.

\bibitem{Mat16} H. \textsc{Matui}, \emph{{\'E}tale groupoids arising from products of shifts of finite type}, Adv. Math. \emph{303} (2016), 502--548.

\bibitem{Nek04} V. \textsc{Nekrashevych}, \emph{Cuntz-Pimsner algebras of group actions}, J. Operator Theory \emph{52} (2004), no. 2, 223--249.

\bibitem{Nek18a} V. \textsc{Nekrashevych}, \emph{Finitely presented groups associated with expanding maps}, Geometric and cohomological group theory, 115--171, London Math. Soc. Lecture Note Ser., 444, Cambridge Univ. Press, Cambridge, 2018.

\bibitem{Nek18b} V. \textsc{Nekrashevych}, \emph{Palindromic subshifts and simple periodic groups of intermediate growth}, Ann. of Math. (2) \emph{187} (2018), no. 3, 667--719.

\bibitem{Nek19} V. \textsc{Nekrashevych}, \emph{Simple groups of dynamical origin}, Ergodic Theory Dynam. Systems \emph{39} (2019), no. 3, 707--732.

\bibitem{PS} G. \textsc{Paolini} and M. \textsc{Salvetti}, \emph{Proof of the $K(\pi,1)$ conjecture for affine Artin groups}, Invent. Math. \emph{224} (2021), no. 2, 487--572.

\bibitem{Par} L. \textsc{Paris}, \emph{$K(\pi,1)$ conjecture for Artin groups}, Ann. Fac. Sci. Toulouse Math. (6) \emph{23} (2014), no. 2, 361--415.

\bibitem{Rov} C.E. \textsc{R{\"o}ver}, \emph{Constructing finitely presented simple groups that contain Grigorchuk groups}, J. Algebra \emph{220} (1999), no. 1, 284--313.

\bibitem{SWZ} R. \textsc{Skipper}, S. \textsc{Witzel} and M.C.B. \textsc{Zaremsky}, \emph{Simple groups separated by finiteness properties}, Invent. Math. \emph{215} (2019), no. 2, 713--740.

\bibitem{Spi20} J. \textsc{Spielberg}, \emph{Groupoids and $C^*$-algebras for left cancellative small categories}, Indiana Univ. Math. J. \emph{69} (2020), no. 5, 1579--1626.

\bibitem{Sta} C. \textsc{Starling}, \emph{Boundary quotients of $C^*$-algebras of right LCM semigroups}, J. Funct. Anal. \emph{268} (2015), no. 11, 3326--3356.

\bibitem{Ste} M. \textsc{Stein}, \emph{Groups of piecewise linear homeomorphisms}, Trans. Amer. Math. Soc. \emph{332} (1992), no. 2, 477--514.

\bibitem{Wit19} S. \textsc{Witzel}, \emph{Classifying spaces from Ore categories with Garside families}, Algebr. Geom. Topol. \emph{19} (2019), no. 3, 1477--1524.



\end{thebibliography}
\end{document}